\documentclass{informs3} 

\OneAndAHalfSpacedXI 

\usepackage{endnotes}
\let\footnote=\endnote

\usepackage{fix-cm}
\usepackage{amsmath}
\usepackage{amsfonts}
\usepackage{graphicx}
\usepackage{geometry}
\usepackage{enumerate}
\usepackage[normalem]{ulem}
\usepackage{color}
\usepackage{float}
\usepackage{xcolor} 
\usepackage{environ}

\newcommand{\rr}[1]{{#1}}

\usepackage[final]{showkeys}
\usepackage[colorlinks=true,linkcolor=cyan,citecolor=blue,anchorcolor=blue]{hyperref}
\newcommand{\TX}{\widetilde{X}}

\newcommand{\E}{\mathbb{E}}
\newcommand{\M}{\mathcal{M}}
\newcommand{\R}{\mathbb{R}}

\newcommand{\Prob}{\mathbb{P}}

\newcommand{\lipone}{\text{\rm Lip(1)}}
\providecommand{\abs}[1]{\left\lvert#1\right\rvert}
\providecommand{\norm}[1]{\lVert#1\rVert}
\newcommand{\startproof}{\proof}
\newcommand{\finishproof}{\hfill $\square$\endproof}

\usepackage{natbib}
\PassOptionsToPackage{caption=false}{subfig}
\usepackage{subfig}

\usepackage{natbib}
 \bibpunct[, ]{(}{)}{,}{a}{}{,}%
 \def\bibfont{\small}%
\TheoremsNumberedThrough     
\ECRepeatTheorems

\EquationsNumberedThrough    

\begin{document}


\RUNAUTHOR{Braverman, Scully}

\RUNTITLE{Stein's method and general clocks}

\TITLE{Diffusion approximation error for queueing systems with general primitives}

\ARTICLEAUTHORS{%
\AUTHOR{Anton Braverman}
\AFF{Kellogg School of Management, Northwestern University, Evanston, IL 60208, \EMAIL{anton.braverman@kellogg.northwestern.edu}}  
\AUTHOR{Ziv Scully\thanks{Supported by NSF grant no.\ CMMI-2307008}}
\AFF{School of Operations Research and Information Engineering,  Cornell University, Ithaca, NY 14853, \EMAIL{zivscully@cornell.edu}}
} 

\ABSTRACT{%
We investigate the steady-state diffusion-approximation error for continuous-time queueing systems  with
generally distributed primitives. \rr{A common picture emerges after analyzing a number of canonical systems:} the error decomposes into interior and boundary terms. The former are simpler to handle and can be bounded using only low-order moments  of the system's primitives --- when the  approximation error is measured using the Wasserstein distance, three moments suffice. The boundary terms are inherently more delicate: while   crude bounds are easy to obtain, sharper (e.g., order optimal) bounds require   deeper, model specific, insights.

Methodologically, we extend the generator comparison approach of Stein's method to piecewise-deterministic Markov processes (PDMPs). The discontinuous nature of the PDMP at jump times necessitates using   the  basic adjoint relationship (BAR), instead of the infinitesimal generator, to characterize the stationary distribution. A second-order Taylor expansion of the BAR’s jump terms, coupled with a Palm-inversion step that converts event-averaged quantities into time averages, yields the candidate diffusion generator and a transparent interior/boundary error decomposition. In parallel, we show how the prelimit generator approach --- working with the Poisson equation of the queueing system instead of the diffusion process --- offers a promising avenue for bounding the challenging boundary terms.
}%



\maketitle

%


\section{Introduction}
Diffusion approximations play a central role in the analysis of queueing systems. They depend on the first two moments of the system’s primitives and can yield valuable insights into the performance of models that are otherwise intractable. These approximations are often justified through heavy-traffic limit theorems, which ensure that the original system converges to a diffusion limit as its parameters are taken to an appropriate asymptotic regime. Yet a key limitation of such limit theorems is that they provide no information about the accuracy of the approximation when the system operates outside this asymptotic regime. Over the past decade, this gap has been addressed by the generator approach of Stein’s method \citep{Barb1988,Barb1990,Gotz1991}, which offers a systematic way to bound approximation errors even in non-asymptotic settings.

A notable contribution of the generator approach has been in establishing the universality of diffusion approximations --- that is, demonstrating the accuracy of the diffusion model simultaneously across multiple heavy-traffic regimes, such as quality-driven, quality- and efficiency-driven, and efficiency-driven regimes (see, e.g., \cite{GlynWard2003,GurvHuanMand2014}). This is due to the diffusion model's dependence on the finite-parameter specification of the queueing system (in contrast to the diffusion limit, which depends on the limit of the parameters). In this paper, we explore a different dimension of universality: for continuous-time queueing models with generally distributed primitives (e.g., general interarrival and service times), we investigate how the diffusion approximation error behaves across primitive distributions.

The assumption of exponentially distributed primitives is often adopted for analytical tractability; in such models, the queue-length process naturally admits a continuous-time Markov chain (CTMC) representation. Models with general primitives are important not only because they better reflect real-world variability, but also because a system’s behavior can change qualitatively when exponential assumptions are relaxed. For example,  generalized Jackson networks no longer have a product-form stationary distribution when interarrivals and service times are not exponentially distributed; for further examples see \cite{Whit1986a, DaiHaseVand2004, BasaRand2010}. In such cases, maintaining a Markovian representation requires augmenting the state space with residual inter-event times; for instance,  the $G/G/1$ queue-length process includes both the residual interarrival time and the residual service time. Phase-type primitives, which can approximate any positive-valued distribution arbitrarily well, also admit a CTMC representation, but the enlarged state space significantly complicates the error analysis \citep{Gurv2014,BravDai2017}.

In what follows, we introduce the terminology of a \emph{general clock} to denote a residual inter-event time when the corresponding  distribution is general (non-exponential).  For instance, in the Markovian representation of the $G/G/1$ queue length,  both the residual interarrival and  service times are general clocks, whereas in the $G/M/\infty$ system where service times are exponentially distributed, only the residual interarrival time is a general clock.

Our main goal is to establish bounds on the distance between the steady-state distribution of the queueing system and its diffusion approximation; in this paper, we focus on the Wasserstein distance. We seek error bounds that (a) hold across a broad class of primitive distributions, (b) depend on those distributions only through low-order statistics (e.g., moments of the distribution), and (c) converge to zero at the optimal rate. The reassurance brought by such distributional universality   is particularly valuable when calibrating model primitives from data: while accurately estimating the exact distribution is often infeasible, low-order statistics --- such as the first three moments --- can  be estimated with much greater precision.

We analyze four canonical queueing systems: the $G/G/1$ queue, the
join-the-shortest-queue system (JSQ), the $G/M/\infty$ queue, and the tandem queue. Across these models, a common structure
emerges: the diffusion-approximation error naturally decomposes into \emph{interior} and
\emph{boundary} terms. \rr{Provided that the Stein factor bounds for the approximating diffusion are known (this is true for all but the tandem queue example, for which Stein factor bounds remain an open problem),} establishing properties (a)--(c) for the interior terms is relatively
straightforward --- we can bound them using only the first three moments of the primitive
distributions, and we expect this to hold more broadly beyond the cases studied here. The boundary
terms, by contrast, present a greater challenge. A useful analogy is a random walk reflected at the
origin to remain on $[0,\infty)$: away from the boundary the dynamics match those of the unconstrained
walk, while at the reflection point the behavior becomes markedly more intricate.

To highlight this added difficulty, we analyze the boundary terms in a representative case: the
workload process of the $G/G/1$ queue. A simple bound can be obtained using H\"older’s inequality,
which ensures properties (a) and (b) but yields a suboptimal convergence rate; this technique
extends to other models as well. Achieving a bound that also satisfies (c) is substantially harder.
Despite considerable effort, we were only able to establish such a bound for a restricted class of
distributions with nonincreasing hazard rates (see
Theorem~\ref{thm:prelimit:idle:bound}). Although these partial results do not rule out the existence
of simple, universally applicable bounds, they underscore the significant challenge of attaining
such generality.

Achieving our stated objective  requires the development of new methodological tools, which form the second main contribution of this work. We extend the generator approach of Stein’s method, which, to date, has only been applied to CTMCs and discrete-time systems, to models with generally distributed primitives. Such models naturally admit a piecewise-deterministic Markov process (PDMP) representation \citep{Davi1984}. In contrast to \cite{Davi1984}, which formulates PDMPs via hazard rates, we adopt a residual-time representation that tracks the remaining time to the next event, or jump, of the process. Extending Stein’s method to such PDMPs requires new tools. In particular, we replace the infinitesimal generator with the basic adjoint relationship (BAR) --- the stationary balance equation for the process --- because the PDMP does not have a well-defined infinitesimal generator at the boundaries where jumps occur.

Starting from the BAR, we use a Taylor expansion to bridge the stationary equation of the queueing model to the generator of its diffusion approximation, mirroring the approach in the CTMC setting. The novelty here lies in the additional jump terms in the BAR due to the general renewal processes --- these terms can be naturally interpreted as Palm expectations. To extract the diffusion generator from these jump terms, we rely on the Palm inversion formula, which links Palm expectations (event-averaged quantities) to time-averaged expectations, because the event-based nature of these terms prevents a direct Taylor expansion. We emphasize that, although our analysis has a close connection to Palm calculus, all results are derived directly from first principles, without invoking the broader machinery of that theory.

Returning to the boundary-term difficulty noted earlier, one promising way to circumvent it is the prelimit generator approach \citep{Brav2022}. This method works with the Poisson equation of the queueing system itself, rather than that of the diffusion approximation. In effect, it reshuffles the problem’s difficulty: instead of bounding complex boundary terms, the challenge is transferred to obtaining Stein factor bounds --- bounds on the second- and third-order derivatives of the Poisson equation’s solution.

For the $G/G/1$ workload process, this reshuffling works in our favor and directly addresses one of the most stubborn boundary terms: the expected equilibrium idle-period length. This term has been extensively studied \citep{Sieg1979,LiOu1995,WolfWang2003}, and although \cite{BlanGlyn2006} provide an asymptotic expansion for it, the only ``clean'' upper bound we know of is due to \cite{Koll1976}. In contrast, while the Stein factor bounds remain technically involved, they lead to novel and explicit bounds on this quantity.

Finally, we note that extending the prelimit approach to the PDMP setting requires tools distinct from those used for our extension of the classical generator approach. We are able to carry out this extension in settings with a single general clock, but extending it to models with multiple clocks remains an open problem. Taken together, our results broaden the scope of both the generator and prelimit approaches beyond CTMCs and discrete-time models, and yield new insights into the quality of diffusion approximations of queueing models.

\subsection{Literature review}
\label{sec:litrev}

The generator comparison approach of Stein's method provides a principled framework for comparing stationary
distributions of Markov processes. While Stein’s method originates with \cite{Stei1972}, the link to
Markov processes was developed by \cite{Barb1988,Barb1990} and \cite{Gotz1991}. Queueing theory has
been a particularly fertile domain: early applications to birth–death processes appear in
\cite{BrowXia2001}. Subsequent work analyzed more complex systems --- see, for example,
\cite{GurvHuanMand2014,Stol2015,Gurv2014} --- which, although not framed explicitly in Stein’s
terminology, employ the core generator-comparison ideas. The approach was later systematized and
popularized in queueing settings by \cite{BravDaiFeng2016,BravDai2017} for diffusion approximations,
and by \cite{ying2016,Ying2017,Gast2017} for mean-field models. Since then, the generator approach
within Stein’s method has remained an active line of research in the queueing community.

A related line of work analyzes approximation errors for \emph{discrete-time} queueing systems, both
via the drift method \citep{EryiSrik2012,MaguSrik2016,MaguBurleSrik2018,HurtMagu2022a} and via
Stein’s method \citep{FengShi2018,ZhouShro2020,HurtMagu2022}. Discrete time affords
simplifications not available in our continuous-time setting; nevertheless, the approximation error
in these works likewise decomposes into interior and boundary components. For example, one  common boundary term is related to the unused service  within a slot, and is commonly controlled by imposing a bound on the
number of service completions per slot.

 Another relevant line of work is the BAR approach \citep{Miya2015,Miya2017,BravDaiMiya2017,BravDaiMiya2024,DaiXu2024,DaiGlynXu2025}, which develops heavy-traffic steady-state limit theorems for queueing systems. Its core device is an  exponential test function, engineered so that all but an asymptotically negligible portion of the BAR's jump terms vanish, which  facilitates establishing heavy-traffic limit theorems via the convergence of moment generating functions. In contrast, our test functions come from the   Poisson equation, which enables us to establish upper bounds on the Wasserstein distance.

Our $G/G/1$ example has a long history. The seminal work of \cite{King1961a,King1962} initiated an
extensive literature on upper and lower bounds for the expected waiting time; see
\cite{Daleetal1992,WolfWang2003} for surveys. Because the expected waiting time and expected workload
are tightly linked \citep[Corollary~X.3.5]{Asmu2003}, bounds for one translate directly to bounds
for the other.

While many classical bounds are asymptotically tight as $\rho\uparrow 1$, they often do not quantify
the nonasymptotic gap to the true expectation. Notable exceptions include \citep[Eq.~(7)]{Boonetal2023},
who use transform methods to obtain moment formulas and an approximation to the expected waiting time with an
$O(1/\abs{\log(1-\rho)})$ error bound, albeit with an implicit constant. By contrast, our result provides
an $O(1)$ error bound --- uniform in $\rho$ as $\rho\to 1$ --- with an explicit constant. Another
approximation that is $o(1)$ accurate (error vanishing as $\rho\to 1$) is stated following
\citep[Theorem~3.3]{BertGama2022}. A related line of research studies \emph{extremal} queues
\citep{ChenWhit2020,ChenWhit2021,ChenWhit2022,ChenWhit2022a}, identifying interarrival and
service-time distributions for which waiting-time bounds are attained with equality.

Beyond the many bounds for the expected waiting time, there is a growing body of work applying Stein’s method to
single-server queues. Using equilibrium couplings, \cite{GaunWalt2020} obtain exponential
approximation bounds for the $G/G/1$ waiting time --- building on the approaches of
\cite{PekoRoll2011,Ross2011} and exploiting the representation of the waiting time as a geometric sum of i.i.d.\
increments. They also derive error bounds for the $M/G/1$ system via the generator approach. In line
with our Theorem~\ref{thm:workload:1}, their $G/G/1$ bound depends on the expected equilibrium idle
period.

In \cite{GurvHuan2018}, the generator comparison approach is applied to the workload of the
$M/G/1{+}GI$ system (Poisson arrivals with a general patience-time distribution). The focus there is
on diffusion-approximation error bounds that are universal across patience-time distributions and load
regimes, together with bounds for arbitrary moments of the stationary workload.

Finally, \cite{Brav2022} gives a simple prelimit treatment of the $M/M/1$ queue length. At the
process level, \cite{Besaetal2020} establish rates of convergence to diffusion limits for the
$M/M/1$ and $M/M/\infty$ systems. More recently, \cite{Barbetal2023} develop Gaussian-process
approximations and apply them to the $G/G/\infty$ system.

\subsection{Notation}
\label{sec:notation}
We write $f^{(k)}(\xi)$ to denote a generic  $k$th-order Taylor expansion remainder term that may change from expression to expression. For a function $f: \R^{d} \to \R$ and integer $k$, we let $\partial_{x_i}^{k} f(x) = \partial^{k} f(x)/\partial x^{k}_{i}$. We let   $\lipone$ denote  the set of Lipschitz-$1$ functions and let $\M_{k}$ denote the set of all Lipschitz-$1$ functions whose derivatives  up to the $(k-1)$st order are also Lipschitz-$1$. 
For a counting process $\{N(t) : t \geq 0\}$ with event times $\{\tau_{m}\}_{m=1}^{\infty}$, we define 
\begin{align*}
    \int_{0}^{t} f(Z(s)) d N(s) = \sum_{m=1}^{\infty} f(Z(\tau_{m})) 1(\tau_{m} \leq t),
\end{align*}
and we let  $\Delta f(Z(t-)) = f(Z(t)) - f(Z(t-))$, where $f(Z(t-)) = \lim_{s \uparrow t} f(Z(s))$.  Given independent random variables $X\in\R^n$ and $Y\in\R^m$ and a measurable
function $f:\R^{n+m}\to\R$, we let $\E^X f(X,Y)$ denote the random variable $\int_{\R^{n}} f(x,Y) d \mu_{X}(x)$, where $\mu_{X}$ is the law of $X$.

\subsection{Outline for the rest of the paper}
\label{sec:outline}
\rr{
Using the workload process of the $G/G/1$ system as a simple introductory example, Section~\ref{sec:workload} illustrates how to extend the generator approach to PDMPs, while Section~\ref{sec:prelimit} develops the parallel prelimit generator approach. Section~\ref{sec:stein:factor:bounds} establishes the Stein factor bounds for the workload process needed in the prelimit analysis. Section~\ref{sec:jsq} analyzes the total customer count in the JSQ system, demonstrating how the classical generator approach extends to a model with multiple general clocks and highlighting the more complicated boundary terms that arise in such systems. Section~\ref{sec:conclusion} concludes with open research directions.
}

\rr{
In the  appendix, we analyze two additional examples: the $G/M/\infty$ queue and the tandem queue. Appendix~\ref{sec:gminfinity} studies the $G/M/\infty$ queue, which illustrates the BAR/Palm-inversion machinery in a setting with both general and exponential clocks, but where the approximation error contains only interior terms. Appendix~\ref{app:tandem} studies the tandem queue, showing how the  approach extends to a multidimensional diffusion approximation while also highlighting the current lack of Stein factor bounds for the corresponding reflected Brownian motion.
}

\section{The G/G/1 workload} 
\label{sec:workload}
The workload process of the single-server queue has only one general clock, which tracks arrivals to the system, making it a natural introductory example. After defining the model, we derive the BAR and use it to establish several useful identities for the system. We then extract a diffusion-approximation generator from the BAR via a Taylor expansion and the Palm inversion formula. 

\rr{Using the Poisson equation for the limiting (exponential) distribution, we compare the diffusion generator with the expanded $G/G/1$ BAR and obtain the resulting approximation error in Theorem~\ref{thm:workload:1}. We follow the theorem with a discussion of the associated error terms, in particular the boundary term involving the system idle period. This discussion motivates the prelimit approach as an alternative method for handling what is otherwise a historically difficult error term. Before concluding the section, we state Theorem~\ref{thm:prelimit:idle:bound}, which gives the analogous bound obtained via the prelimit approach, and provide a detailed comparison of the two methods, including the advantages and limitations of each.}   
\subsection{The model}
\label{sec:model:workload} 
Consider a single-server queue operating under any non-idling policy. Arrivals follow a renewal process with rate $\lambda$ and interarrival distribution $U$, and the service time distribution $S$ has rate $\mu$. Let $\rho = \lambda/\mu$, $c_{U}^2 = \lambda^2 \text{Var}(U)$, and $c_{S}^2 = \mu^2 \text{Var}(S)$, and assume that 
\begin{align}
    \E U^3 < \infty \quad \text{ and } \quad \E S^3 < \infty. \label{eq:moments:workload}
\end{align}
For simplicity, we assume that simultaneous events (arrivals/departures)  do not occur, with probability one. Let $V(t)$ and $R_{a}(t)$ be the remaining workload in the system and  residual interarrival time (time until the next arrival), respectively, at time $t \geq 0$. 

Let  $\delta = 1-\rho$, $X(t) = \delta V(t)$, and define    the scaled  workload process 
\begin{align*}
    \{Z(t) = (X(t),R_{a}(t)) :  t \geq 0\},
\end{align*}
which is a right-continuous with left limits (RCLL) PDMP taking values in 
  \begin{align*}
  \mathbb{S} = \{(x,r) \in \R^{2}_+ : x \geq 0, r > 0\}.
  \end{align*}
Note that $(x,0)\notin\mathbb{S}$ for $x\geq 0$: since the process is RCLL, whenever the residual interarrival time hits zero an arrival occurs instantaneously, and the process immediately jumps to a state with strictly positive residual time.

 Let  $A(t)$ be  the number of arrivals on $[0,t]$.  If an arrival occurs at time $t$, let $U(t)$ and $S(t)$   be the subsequent interarrival time and the workload brought by the arriving customer, respectively. If no arrival occurs at $t$, set $U(t)=S(t)=0$. 
Conditional on an arrival at $t$, $U(t)\stackrel{d}{=}U$ and $S(t)\stackrel{d}{=}S$, and both are independent of the process history over $[0,t)$.

The workload process is a regenerative process, with regeneration happening at those instances when a customer arrives to an empty system. Going forward we assume that 
\begin{align}
\rho <1  \text{ and $U$ is nonlattice}, \label{eq:stability}
\end{align} 
which ensures the existence of a limiting steady-state distribution \citep[Corollary~X.3.3]{Asmu2003}. We let $Z = (X,R_{a}) \in \mathbb{S}$ have this distribution. 
\subsection{The BAR and some basic identities} \label{sec:firstBAR}
\rr{We begin by deriving the BAR for the workload process, which characterizes the stationary distribution of $Z$. The derivation relies on two assumptions: pathwise fundamental theorem of calculus (FTC) conditions, which justify decomposing the sample-path dynamics, and an integrability condition, which justifies taking expectations and interchanging limits and integrals. These assumptions have natural analogues in the other queueing models considered in later sections, so this section also serves as a template for deriving the BAR in those settings. }

\rr{Let $\tau_0=0$, and suppose that arrivals occur at times $\{\tau_m\}_{m=1}^{\infty}$. More generally, beyond the workload process, one may interpret $\tau_m$ as the time of the $m$th event in the system, where an event may be, for example, an arrival or a departure. Recall that, by assumption, simultaneous events do not occur, with probability one. This can be relaxed with a minor modification of the arguments that follow; see also the discussion around (3.14) in \cite{BravDaiMiya2017}.}

\rr{Given any $f:\mathbb S\to\mathbb R$ and initial condition $Z(0)\in\mathbb S$, we decompose the evolution of the workload on $[0,t]$ into its jump contributions and interarrival intervals by writing
\begin{align*}
    &f(Z(t)) - f(Z(0)) \\
    =&\  \big(f(Z(t)) - f(Z(\tau_{A(t)}))\big) + \sum_{m=1}^{A(t)} \big(f(Z( \tau_{m}-))  - f(Z(  \tau_{m-1}))\big)   + \int_{0}^{t} \Delta f(Z(s-))d A(s).
\end{align*} 
For each sample path of the workload process, we say that the FTC conditions are satisfied if the function $s\mapsto f(Z(s))$ 
\begin{align*} 
&(i) \text{  is absolutely continuous on each interval $[\tau_{m-1},\tau_{m}]$, $1 \leq m \leq A(t)$, when evaluated} \notag\\
    &\hspace{.75cm} \text{as $f(Z(\tau_{m-1}))$ at the left endpoint and $f(Z(\tau_{m}-))$  at the right endpoint,} \notag\\
    &(i') \text{  is absolutely continuous on $[\tau_{A(t)},t]$, when evaluated} \notag\\
    &\hspace{.75cm} \text{as $f(Z(\tau_{A(t)}))$ at the left endpoint and $f(Z(t))$  at the right endpoint,} \notag\\
    &(ii) \text{ satisfies the multivariate chain rule for those $s$ at which $d f(Z(s))/ds$ exists.}  
\end{align*} 
Assume these conditions are satisfied. Then,  whenever $d f(Z(s))/ds$  exists, it follows that
\begin{align*}
    \frac{d}{ds}f(Z(s)) = -\delta 1(X(s)>0)\partial_{x}f(Z(s)) - \partial_{r_{a}} f(Z(s)).
\end{align*}}
\rr{Applying the FTC  to the  telescoping series yields 
\begin{align*}
    &f(Z(t)) - f(Z(0)) \\
    =&\    \int_{0}^{t} \big( -\delta 1(X(s)>0) \partial_{x} f(Z(s)) - \partial_{r_{a}} f(Z(s)) \big) ds +    \int_{0}^{t} \Delta  f(Z(s-)) dA(s).
\end{align*} 
 Let $Z(0) \sim Z$ and   assume that 
\begin{align}
    \E \abs{f(Z)},\ \E \abs{\partial_{x} f(Z)},\ \E \abs{\partial_{r_{a}} f(Z)},\ \E \int_{0}^{t} \abs{\Delta  f(Z(s-))} dA(s) < \infty. \label{eq:ev}
\end{align} 
We refer to \eqref{eq:ev} as the integrability condition. When it is satisfied, we can apply the Fubini-Tonelli theorem to interchange the expectation with the integrals to arrive at the following lemma.
\begin{lemma}
\label{lem:full:BAR:workload}
Fix $f: \mathbb{S} \to \R$ and initialize $Z(0) \sim Z$. If $f(Z(s))$ satisfies the FTC conditions with probability one under $Z(0) \sim Z$ and if the integrability condition \eqref{eq:ev} holds, then 
    \begin{align}
    0 =   -\delta  \E \big(1(X>0) \partial_{x} f(Z)\big) - \E \big(\partial_{r_{a}} f(Z) \big) + \frac{1}{t} \E  \int_{0}^{t} \Delta  f(Z(s-)) dA(s). \label{eq:full:BAR:workload}
\end{align}
\end{lemma} 
We point out that Lemma~\ref{lem:full:BAR:workload} does not require any assumptions beyond the FTC and integrability conditions. In particular, it continues to hold when the interarrival distribution $U$ has atoms, or heavy tails. Rather than giving an exhaustive characterization of when the FTC conditions hold, we highlight two important cases that will be used repeatedly throughout the paper.
}

\rr{
First, observe that $\{Z(t)\}$ is piecewise linear, with finitely many linear pieces between any two arrivals. Moreover, on any interval $[\tau_{m-1},\tau_m)$, one of two cases occurs. If $R_{a}(\tau_{m-1}) \leq V(\tau_{m-1})$, then the workload remains in $\operatorname{int}(\mathbb S)$ throughout the interval. If instead $R_{a}(\tau_{m-1}) > V(\tau_{m-1})$, then the workload starts in $\operatorname{int}(\mathbb S)$, eventually hits $\{0\}\times(0,\infty)$, which lies on the boundary of $\mathbb S$, and remains there until the next arrival at $\tau_m$. The same reasoning applies on the final interval $[\tau_{A(t)},t]$.
}

\rr{
It follows that if $f \in C^1(\operatorname{int}(\mathbb S))$ and the boundary restriction $f(0,r_a)$ belongs to $C^1(0,\infty)$ as a function of $r_a$, then $f(Z(s))$ satisfies the FTC conditions for every initial condition $Z(0)\in\mathbb S$. In most applications, $f(z)$ will be the solution to a Poisson equation and hence will typically be at least twice continuously differentiable. 
}

\rr{
Second, we will frequently work with truncated functions. These are useful because they avoid the need to verify $\E \abs{f(Z)}<\infty$ directly: truncation makes $f(Z)$ bounded, and hence integrable. The tradeoff is that truncation typically introduces points at which the function is not differentiable. This does not create a problem for the FTC conditions, provided that the workload process spends zero Lebesgue time at those nondifferentiability points. Two examples that we will use repeatedly are
\begin{align*}
    f(x,r_a)=x \wedge M
    \qquad\text{and}\qquad
    f(x,r_a)=r_a \wedge M.
\end{align*}
Both satisfy the FTC conditions, since the workload process crosses the truncation level at most once between any two arrivals. We do not attempt to give a general definition of all admissible truncated functions. Rather, in each instance where such a function is used, the FTC conditions are straightforward to verify directly.
}

Going forward we fix $t = 1$. The BAR grants access to the properties of $Z$. The  following identities can be established by using the BAR with $f(z) = r_{a} \wedge M$,   $f(z) = x \wedge M$, and  $f(z) = r_{a}^m \wedge M$, and taking $M \to \infty$.
\begin{lemma}
    \label{lem:relationships:workload}
    In steady state,  $\E A(1) = \lambda$,  $\Prob(X > 0) = \rho$, and $\E R_{a}^{m-1} = \lambda\,\E U^{m}/m$.
\end{lemma}

\rr{In the next section we show how to extract the approximating diffusion generator from the BAR.}
\subsection{Deriving the diffusion generator}
When $f(Z)=f(X,R_{a})$  is a function of $X$ alone, the BAR \eqref{eq:full:BAR:workload} reduces to 
\begin{align}
    0 =   -\delta  \E \big(1(X>0) f'(X)\big)   + \E  \int_{0}^{1} \big( f(X(t-)+\delta S(t))- f(X(t-)) \big) dA(t). \label{eq:BAR:uncompensated}
\end{align}
Deriving the diffusion generator requires expanding the jump term on the right-hand side. In the special case of Poisson arrivals,   the ``PASTA'' property (Theorem~3.3.1 of \cite{BaccBrem2003}; see also their equation (1.8.8)) implies
\begin{align*}
    &  -\delta  \E \big(1(X>0) f'(X)\big)   + \E  \int_{0}^{1} \big( f(X(t-)+\delta S(t))- f(X(t-)) \big) dA(t) \\
    =&\   -\delta  \E \big(1(X>0) f'(X)\big)   +   \lambda \E  \big( f(X+\delta S)- f(X) \big).
\end{align*}
The right-hand side coincides with the generator of the $M/G/1$ workload process; see \cite{GurvHuan2018} and \cite{GaunWalt2020} for analyses of the $M/G/1$ workload using the generator approach.

For general arrival processes, we consider the  \emph{compensated} workload
\begin{align}
    \TX(t) = X(t)-\delta \rho R_{a}(t) \quad \text{ and } \quad \TX = X - \delta \rho R_{a}. \label{eq:compensated:workload}
\end{align}
The defining property of the compensated workload is that at all arrival instances $t$, 
\begin{align*}
    \E \Delta \TX(t-) = \delta(\E S(t) - \rho \E U(t)) = \delta (\E S - \rho \E U) = 0.
\end{align*}
Specializing the BAR \eqref{eq:full:BAR:workload} to functions of the form $f(z) = f(x-\delta \rho r_{a})$ yields
\begin{align}
  & 0 = \delta  \E  \big((\rho -  1(X>0))f'(\TX )\big) +  \E \int_{0}^{1}   \Delta f(\TX(t-))    d A(t).\label{eq:BAR2:workload}
\end{align}  
With the compensated process,  the  first-order term in the jump-term expansion of \eqref{eq:BAR2:workload} equals zero. We elaborate on why this matters at the end of Section~\ref{sec:extracting:workload} after proving the following key proposition, which expands \eqref{eq:BAR2:workload} into a second-order differential operator plus an error term.  Recall that $f^{(k)}(\xi)$ denotes a generic $k$th-order Taylor remainder term, whose value may change from line to line. When the expansion point depends on time, we write $f^{(k)}(\xi(t))$ to make this dependence explicit.
\begin{proposition}
\label{prop:extraction:workload}
If  $f\in C^{2}(\R)$ with $\norm{f''} < \infty$ and $f''(x)$ absolutely continuous, then, provided that all expectations are well defined, 
\begin{align}
    \delta  \E  \big(\rho - 1(X>0)\big)f'(\TX) =&\ -\delta^2 \E f'(X) + \delta^2 f'(0) + \epsilon_{0}(f), \label{eq:e0:workload} \\
        \E \int_{0}^{1} \Delta  f(\TX(t-))  d A(t) =&\ \frac{1}{2}\delta^2 \lambda   \E (S -\rho  U)^2 \E f''(X) + \epsilon_{A}(f)  \label{eq:eA:workload}
    \end{align}
    where 
    \begin{align*}
        \epsilon_{0}(f) =&\ \delta^3 \rho \E(R_{a} f''(\xi)) - \delta^2 \rho \E (1(X=0)R_{a} f''(\xi)),\\
        \epsilon_{A}(f) =&\ \frac{1}{6} \delta^3\E \int_{0}^{1} (S(t) - \rho U(t))^3 f'''(\xi(t))d A(t) \\
        &- \frac{1}{2} \delta^2\lambda \E (S -\rho  U)^2  \E \int_{0}^{1} \int_{0}^{U(t)} \big( X(t+u)-X(t-)\big) f'''(\xi(t+u))du d A(t).
    \end{align*}
\end{proposition}
Note that although the left-hand sides of \eqref{eq:e0:workload} and \eqref{eq:eA:workload} are in terms of $\TX(t)$, the right-hand sides (excluding the error terms) are in terms of the uncompensated $X$.

The terms $-\delta^2 \E f'(X)$ and $\delta^2 f'(0)$ in \eqref{eq:e0:workload} capture the drift of the diffusion and its reflection at the boundary, while the term involving $\E f''(X)$ in \eqref{eq:eA:workload} captures its variability. Combining these terms and noting  that 
\begin{align*}
  \lambda \E (S -\rho  U)^2 = \lambda \text{Var}(S-\rho U) = \lambda \text{Var}(S) + \lambda \rho^2 \text{Var}(U)  =  \rho \E S (c_{U}^{2}+c_{S}^{2}),
\end{align*} 
we arrive at the diffusion generator 
\begin{align*}
    G_{Y} f(x) = -\delta^2 f'(x) + \frac{1}{2}\delta^2 \rho \E S (c_{U}^{2}+c_{S}^{2}) f''(x) + \delta^2 f'(0), \quad x \geq 0.
\end{align*}
This generator corresponds to the one-dimensional reflected Brownian motion (RBM)  \citep{HarrReim1981} and the corresponding stationary distribution is exponential with mean $\rho \E S (c_{U}^{2}+c_{S}^{2})/2$. We let $Y$ be the random variable having this distribution. 

 \rr{A critical tool in the proof of Proposition~\ref{prop:extraction:workload} is a relationship between event-average and time-average expectations. We state this relationship next, and then use it to prove the proposition.  }

\subsubsection{Palm inversion and the proof of Proposition~\ref{prop:extraction:workload}}
\label{sec:extracting:workload}
The following is a special case of the  Palm inversion formula \citep[Equation~(1.2.25)]{BaccBrem2003}. We prove it  in Appendix~\ref{app:workload:inversion}. 
Note that the conditions of Lemma~\ref{lem:palm:inversion:2} are  trivially satisfied when $f(x)$ is bounded.  
\begin{lemma}
    \label{lem:palm:inversion:2}
   Initialize $Z(0) \sim Z$ and fix  $f: \R \to \R$ with
    \begin{align*}
     \E\abs{f(X)}< \infty \quad \text{ and  } \quad \E \abs{\int_{0}^{R_{a}(0)} f(X(t)) dt} < \infty.
 \end{align*}
 Then 
    \begin{align}
        \E f(X) =&\ \E   \int_{0}^{1} \int_{0}^{U(t)} f(X(t+u))du  d A(t).  \label{eq:palm:inversion:arrive:workload}
    \end{align}
\end{lemma}
We now use Lemma~\ref{lem:palm:inversion:2} to expand the BAR for the compensated workload.
\startproof{Proof of Proposition~\ref{prop:extraction:workload}}
We first prove \eqref{eq:e0:workload}, which follows from some basic algebra:
\begin{align*}
    &\delta  \E  \big(\rho -1(X>0)\big)f'(X-\delta \rho R_{a})  \\
    =&\ -\delta^2 \E f'(X-\delta \rho R_{a}) + \delta \E \big( 1(X=0)f'(-\delta \rho R_{a})\big) \\
    =&\ -\delta^2 \E \big(f'(X)+ \delta(-\rho R_{a})f''(\xi)\big) + \delta \E \big( 1(X=0) \big( f'(0) - \delta(\rho R_{a})f''(\xi)\big)\big) \\
    =&\ -\delta^2 \E f'(X)+ \delta^2 f'(0)  +\delta^3 \rho \E(R_{a} f''(\xi)) - \delta^2 \rho \E (1(X=0)R_{a} f''(\xi)),
\end{align*} 
where in the last equality we used $\Prob(X=0)=1-\rho$.
Next we prove \eqref{eq:eA:workload}. Since $\E \Delta \TX(t-) =0$ at all arrival instants $t$, it follows that 
\begin{align}
    \E \int_{0}^{1} \Delta  f(\TX(t-))  d A(t)  =&\   \frac{1}{2} \delta^2\E (S -\rho  U)^2  \E \int_{0}^{1}f''(\TX(t-))d A(t) \notag \\
    &+ \frac{1}{6} \delta^3\E \int_{0}^{1} (S(t) - \rho U(t))^3 f'''(\xi(t))d A(t). \label{eq:analogous:2}
\end{align}
Noting that $\TX(t-) = X(t-)$ at jump times $t$,  to conclude \eqref{eq:eA:workload} it suffices to  show that 
\begin{align*}
    &\E \int_{0}^{1}  f''(X(t-))  d A(t)\\
    =&\ \lambda \E f''(X) - \lambda \E  \int_{0}^{1} \int_{0}^{U(t)} \big(X(t+u)-X(t-)\big)f'''(\xi(t+u)) du  d A(t),
\end{align*}
which follows from expanding the right-hand side of \eqref{eq:palm:inversion:arrive:workload} in  Lemma~\ref{lem:palm:inversion:2} as follows:
\begin{align*}
    \E f''(X) =&\   \E  \int_{0}^{1} U(t) f''(X(t-))  d A(t) \\ 
    &+ \E  \int_{0}^{1} \int_{0}^{U(t)} \big(X(t+u)-X(t-)\big)f'''(\xi(t+u)) du  d A(t).
\end{align*} 
\finishproof
\begin{remark}
\label{rem:compensated}
The first-order term in the expansion of $\E \int_{0}^{1} \Delta f(\TX(t-)) d A(t)$ in \eqref{eq:analogous:2} is zero because the compensated workload satisfies $ \E \Delta \TX(t-)=0$ at jump times $t$. In contrast,  the jump-term expansion in the BAR  for the uncompensated workload \eqref{eq:BAR:uncompensated} contains a first-order term that equals $\delta \E S\E \int_{0}^{1} f'(X(t-)) dA(t)$. Similar to the final step in the proof of Proposition~\ref{prop:extraction:workload}, we can relate this term to $\delta \E S\E f'(X)$ by expanding the right-hand side of \eqref{eq:palm:inversion:arrive:workload} in Lemma~\ref{lem:palm:inversion:2}. Namely,  
\begin{align*}
   \delta \E S\E  \int_{0}^{1}  f'(X(t-))   d A(t) =&\ \lambda\delta \E S  \E f'(X)   \\
    &-\lambda \delta \E S  \E  \int_{0}^{1} f''(X(t-)) \int_{0}^{U(t)} \big(X(t+u)-X(t-)\big) du  d A(t) \\
    &-  \frac{1}{2} \lambda \delta \E S \E  \int_{0}^{1} \int_{0}^{U(t)} \big(X(t+u)-X(t-)\big)^2f'''(\xi(t+u)) du  d A(t).
\end{align*}
Let us informally discuss the last two terms on the right-hand side.  In the following sections we will work with $f(x)$ such that $\norm{f''},\norm{f'''} = O(\delta^{-2})$. Note that 
\begin{align*}
    X(t+u)-X(t-) = (X(t-) + \delta S(t)- \delta u)^{+} - X(t-)
\end{align*}
is of order $\delta$, so that   the third term is   an  $O(\delta)$ error term. The second term however is  $O(1)$ and  cannot be treated as error. We must use the Palm inversion formula a second time to separate
\begin{align*}
    \E  \int_{0}^{1} f''(X(t-)) \int_{0}^{U(t)} \big(X(t+u)-X(t-)\big) du  d A(t)
\end{align*}
into a term involving  $\E f''(X)$ (that we include in the diffusion generator), plus a third-order error term. However, this task requires computing $\E \int_{0}^{U(t)} \big(X(t+u)-X(t-)\big) du$. This is relatively simple (though algebraically involved) for the workload process, because no other jumps happen at times $u \in [t,t+U(t))$. However, if we were working with the $G/G/1$ queue-length process, computing $\E \int_{0}^{U(t)} \big(X(t+u)-X(t-)\big) du$ would require knowing the expected number of departures during $[t,t+U(t))$, but this quantity is not straightforward.

Thus, working with the compensated process accomplishes two things. First, it reduces the volume of calculations required because we only need to use the Palm inversion formula once. Second, it eliminates the need to compute $\E \int_{0}^{U(t)} \big(X(t+u)-X(t-)\big) du$, which is generally intractable.
\end{remark}
\rr{We now introduce the final major ingredient of our approach: the Poisson equation.}


\subsection{The Poisson equation}
\label{sec:stein:step:workload}
In order to compare the workload to its diffusion approximation, we require the Poisson equation. Given  $\theta,\sigma^2 > 0$ and $h : \R \to \R$ with $\E \abs{h(Y)} < \infty$, consider the Poisson equation for the exponential distribution 
 \begin{align}
       -\theta f_h'(x)+\frac{1}{2}\sigma^2 f_h''(x) =&\ \E h(Y) - h(x), \quad x\geq 0, \notag \\
       f_h'(0)=&\ 0. \label{eq:stein:exponential}
 \end{align}
 One may verify that the solution satisfies
\begin{align*}
    f_h'(x) = - e^{\frac{2 \theta}{\sigma^2} x} \int_{x}^{\infty} \frac{2}{\sigma^2} \big( \E h(Y) - h(y) \big) e^{-\frac{2 \theta}{\sigma^2} y} dy, \quad x \geq 0.
\end{align*}
Since the integral representation is well defined for all $x \in \R$, we may extend $f_h'(x), f_h''(x)$, and the Poisson equation \eqref{eq:stein:exponential} to the entire real line.

  Setting $x = X$ and taking expectations allows us to compare $\E h(Y)$ to $\E h(X)$ by comparing the left-hand side to the expanded BAR in Proposition~\ref{prop:extraction:workload}. To do so however, we require bounds on the derivatives of $f_h(x)$.  Since the solution to \eqref{eq:stein:exponential} is unique up to an additive constant,  we assume, without loss of generality, that $f_h(0)=0$.  The following lemma is a rescaled version of \cite[Lemma~4.1]{PekoRoll2011}. 
\begin{lemma}
\label{lem:exp:stein:factors}
Extend \eqref{eq:stein:exponential} to all $x \in \R$ and let $f_h: \R \to \R $ be the  solution  with $f_h(0)=0$. Then  $f_h'(0)=0$. Moreover, if $h(x)$ is Lipschitz, then $f_h'''(x)$ is absolutely continuous and 
    \begin{align*}
       \norm{f_h''} \leq 1/\theta \quad \text{ and } \quad \norm{f_h'''} \leq 4/\sigma^2.
    \end{align*}
As a consequence, for all $x \in \R$, both $\abs{f_h'(x)} \leq \abs{x}/\theta$ and $\abs{f_h(x)} \leq  x^2/(2\theta)$ (obtained by integrating the bound on $f_h''(x)$). 
\end{lemma}

\subsection{Putting everything together}
Let $d_{W}(X,Y)$ denote the Wasserstein distance between $X$ and $Y$, which is defined as 
\begin{align*}
    d_{W}(X,Y) = \sup_{h \in \lipone} \abs{\E h(X) - \E h(Y)}.
\end{align*}
Combining Proposition~\ref{prop:extraction:workload} with the Stein factor bounds for the exponential distribution yields the following bound; we leave the detailed proof to Appendix~\ref{app:workload:theorem}. 
\begin{theorem}
\label{thm:workload:1}
 Assume that \eqref{eq:moments:workload} and \eqref{eq:stability} hold. Then  
\begin{align}
    d_{W}(X,Y) \leq  \abs{\epsilon_{0}(f_h)} + \abs{\epsilon_{A}(f_h)}, \label{eq:dw:workload}
\end{align} 
where $\epsilon_{0}(f_h)$ and $\epsilon_{A}(f_h)$ are as in Proposition~\ref{prop:extraction:workload}. Furthermore, 
    \begin{align*}
                    \abs{\epsilon_{0}(f_h)} \leq&\ \delta \rho \lambda \E U^2/2 + \rho \E ( 1(X=0) R_{a})  , \\
        \abs{\epsilon_{A}(f_h)} \leq&\ \frac{2\delta \E \abs{S-\rho U}^{3}}{ 3  (\E S)^2 (c_{U}^{2}+c_{S}^{2})}  + 2 \delta (\E S + \lambda \E U^2) .
    \end{align*} 
\end{theorem}
The workload $X$ is measured in units of time, and the bounds on
$\lvert\epsilon_{0}(f_h)\rvert$ and $\lvert\epsilon_{A}(f_h)\rvert$  carry the same units. For example,
$\delta\rho\,\lambda\,\mathbb{E}U^{2}/2$ combines the scale-free factors $\delta$ and $\rho$ with $\lambda$ (time$^{-1}$) and  $\mathbb{E}U^{2}$ (time$^{2}$),
yielding time units. Let us say a few words about the boundary term  $\E ( 1(X=0) R_{a})$; the rest of the terms are of order $\delta = (1-\rho)$ and depend only on the first three moments of $U$ and $S$.  

\subsubsection*{A connection to the idle period.}
Note that $\E (R_{a} 1(X=0)) = (1-\rho) \E(R_{a} | X=0)$. Since  the workload process   regenerates at those instances when a customer arrives to an empty system, the workload cycles between busy and idle periods; let $\bar I$ be the length of such an idle period. Thus, conditioned on $X=0$, the distribution of  $R_{a}$ is the same as the distribution of the remaining idle time (which equals the equilibrium distribution of $\bar I$), so that
\begin{align*}
    \E (R_{a} 1(X=0)) = (1-\rho) \E(R_{a} | X=0) = (1-\rho) \frac{\E \bar I^2}{2 \E \bar I}.
\end{align*}

\rr{The idle-period term in the error bound of Theorem~\ref{thm:workload:1} is not merely an artifact of the methodology. To see this, we apply two well known identities (Equation~X.2.5 and Corollary~X.3.5 of \cite{Asmu2003}) to get 
\begin{align}
    \frac{\E \bar I^2}{2 \E \bar I} =&\     \frac{ \lambda \E (S-U)^2}{2(1-\rho)}  - \frac{1}{\rho (1-\rho)} \E X +   \frac{\E S^2}{2\E S}.\label{eq:idle:expression}
\end{align}
Recall that $\E Y = \lambda \E (S-\rho U)^2/2$. By writing $ \lambda \E (S-U)^2 =  \lambda \E (S- \rho U - (1-\rho)U)^2$ and expanding the right-hand side, one can show that
\begin{align}
    \E Y - \E X = \rho(1-\rho) \Big( \frac{\E \bar I^2}{2 \E \bar I}+ \frac{1}{\rho}\E Y -   \frac{1}{2} \lambda (1-\rho)  \E U^2 -  \rho\E U (  \lambda^2 \E U^2-1) -  \frac{\E S^2}{2\E S}  \Big).
\end{align}
Thus, the idle-period term is an unavoidable component of the approximation error itself.}

Significant effort has gone into understanding $\E \bar I^2/(2 \E \bar I)$, which is related to the first descending ladder height \citep[Proposition X.1.5]{Asmu2003}. \cite{LiOu1995} characterize the  distribution function of $\bar I$ as a solution to a certain (complicated) non-linear integral equation.  \cite{WolfWang2003} characterize the asymptotic behavior of the equilibrium distribution of $\bar I$ in a $G/M/1$ system as $\rho \to 1$. \cite{BlanGlyn2006} derive a full asymptotic expansion of $\E \bar I^2 /(2 \E \bar I)$ as $\rho \to 1$ (see also \cite{Sieg1979}). However, none of these papers offer an upper bound on $\E \bar I^2 /(2 \E \bar I)$ in terms of simple model primitives. 

The only  bounds that we are aware of are found in   equations (16) and (17) of \cite{Koll1976}. Letting  $r = (S-U)^{-}$, these are 
\begin{align*}
    \frac{\E \bar I^2}{2 \E \bar I} \leq \frac{\E r^3 / 3 + a \E r^2}{2 \E r \Prob(r \geq a) a} \quad \text{ and } \quad \frac{\E \bar I^2}{2 \E \bar I} \leq \frac{\E r^3 / 3 + a \E r^2}{2 \E (r-a)^{+}   a},
\end{align*}
for any $a > 0$ such that $\Prob(r \geq a)>0$ and $\E (r-a)^{+} > 0$, respectively. Although these bounds require only three finite moments of $U$ and $S$, they are still not ``simple'' as they depend on the entire distribution function of $r$. 

\subsubsection*{A crude bound.}
For a crude but simple bound on the boundary term we can use the fact that $\Prob(X=0)=(1-\rho)$ together with  H\"older's inequality to get 
\begin{align}
      \E (R_{a} 1(X=0)) \leq&\ \big(\Prob(X=0)\big)^{1/p} \big(\E (R_{a}^q )\big)^{1/q} = (1-\rho)^{1/p} \big( \lambda \E U^{q+1}/q \big)^{1/q},  \label{eq:crude:workload}
\end{align}
for any $p,q \geq 1$ with $1/p+1/q = 1$.   Setting  $p = q =2$, we recover  a bound depending only on $\lambda$ and the third moment of $U$,  but the rate of convergence is only $\sqrt{1-\rho}$. We can get a rate of convergence of $(1-\rho)^{1-\epsilon}$ for any $\epsilon > 0$ by taking $q$ large enough, but this comes at the expense of assuming that  higher moments of $U$ exist.

\rr{
It is natural to ask whether the boundary term $\E(R_a1(X=0))$ is also  of order $(1-\rho)$, like the remaining error terms in Theorem~\ref{thm:workload:1}. This is indeed true in several familiar settings. For example, if $U$ is exponentially distributed, then the idle period $\bar I$ is also exponentially distributed by the memoryless property, and hence
\begin{align*}
    \E(R_a1(X=0))
    =
    (1-\rho)\frac{\E \bar I^2}{2\E \bar I}
    =
    \frac{1-\rho}{\lambda}.
\end{align*}
More generally, if $U$ has decreasing mean residual life, then
\begin{align*}
     \E(R_a\mid X=0)\leq \E U=\frac{1}{\lambda},
\end{align*}
and therefore $\E(R_a1(X=0))\leq (1-\rho)/\lambda$. Theorem~\ref{thm:prelimit:idle:bound}, stated below, identifies another class of interarrival distributions for which the boundary term is of order (at most) $1-\rho$. At present, however, we are not aware of a result that guarantees this property for general interarrival distributions, nor of a clean counterexample in the literature showing that this property fails.
}

\subsubsection*{Idle period bounds via the prelimit approach.}  
A key contribution of this paper is to demonstrate how the prelimit generator approach can be used to bound the distance between $X$ and $Y$ without explicitly handling complicated boundary terms such as $\E \bar I^2 /(2 \E \bar I)$. Rather than analyzing these boundary terms directly, the difficulty is transferred to bounding Stein factors.  

The next theorem illustrates the type of bounds obtained under this approach. The main technical challenge lies in establishing universal (across distributions) bounds for the third derivative of the Poisson equation solution. Although an explicit expression is available, obtaining the bounds is delicate and requires additional structure on the interarrival distribution. One tractable family consists of interarrival distributions with a nonincreasing hazard rate, as assumed in our theorem. Another admissible family includes distributions whose hazard rate is uniformly bounded away from zero. In that case, however, $\E \bar I^2 / (2\E \bar I)$ can be bounded directly using the hazard-rate lower bound, and the prelimit-based bounds are strictly weaker. For this reason, we do not include that case in the theorem, although we show how to handle it in the proof.  

In summary, while our machinery extends beyond these two specific distributional families, it relies essentially on exploiting structural properties of $U$ in order to control the third-order Stein factors.   After developing the prelimit approach for models with a single general clock in Section~\ref{sec:prelimit}, we prove the following theorem in Section~\ref{sec:stein:factor:bounds}. 
\begin{theorem}
\label{thm:prelimit:idle:bound}
Assume that \eqref{eq:moments:workload} and \eqref{eq:stability} hold,  and let $Y_{2}$ be exponentially distributed with mean $\lambda \E (S-U)^2/2$. Further assume that $U$ has a nonincreasing hazard rate that is bounded from above by $\overline{\eta} < \infty$. Then for $h(x) = x$,  
    \begin{align*}
        &\abs{\E h(Y_{2}) -  \E h(X)} \leq \delta( K_1 + K_2 + K_3 + K_4+K_5), \\
        &K_1 =   ( \lambda \E U^2/2 +\E S), \\
    &K_2  =   \rho \E U^2/2  + \delta (1/\E Y_2) \lambda (\E U^2/2)(\E S^2 + (\E S)^2)( \delta + 2 \overline{\eta}  \delta \E \bar B ), \\
    &K_3 =    \big(  (1/\E Y_2) \big( \rho \E U^2 + 5\lambda \E U^3/3 \big) +  \lambda \E U^3/3 \big) \big( \delta + (\lambda + \overline{\eta}) \delta \E \bar B  \big)  \\
    &K_{4} = \lambda    \E  S^3  (3   + (\E U + \E \bar I) \lambda  (1 + \rho + \lambda^2 \E U^2))   (1/\E Y_2) \overline{\eta} \delta \E \bar B, \\
    &K_{5} =  \E (S^2 - \rho U^2 + \lambda \E U^3/3)  (3   + (\E U + \E \bar I) \lambda  (1 + \rho + \lambda^2 \E U^2))   (1/\E Y_2) \overline{\eta} \delta \E \bar B.
    \end{align*}
    Furthermore,  $\E \bar I = \delta \E \bar B/\rho$ and 
    \begin{align*}
        &\E \bar B \leq \min\Bigg\{ \frac{\rho\text{Var}(S-U)}{2 (1-\rho) \E ( S - U)^{+}}, \\
        &\hspace{2.5cm}0.9 \rho \frac{\sqrt{\text{Var}(U-S)}}{1-\rho} \exp\bigg(   5.4\frac{\E \abs{U-S}^{3} }{ (\text{Var}(U-S) )^{3/2}} + 0.8 \frac{\E (U-S) }{ \sqrt{\text{Var}(U-S) }}\bigg)\Bigg\}.
    \end{align*}
\end{theorem} 
\rr{Theorem~\ref{thm:prelimit:idle:bound} assumes that the hazard rate of $U$ is nonincreasing.  The technical role of this assumption is highlighted following the statement of Lemma~\ref{lem:mixing:time}.  This assumption covers a broad class of useful interarrival distributions. For instance, mixtures of exponentials have nonincreasing hazard rates, as well as heavy-tailed  Pareto-type interarrival times with
\begin{align*}
    \Prob(U>x)=\left(1+\frac{x}{\beta}\right)^{-\alpha},
\end{align*}
whose hazard rate is $\eta(x)=\alpha/(\beta+x)$, which is nonincreasing and bounded above.
}

\rr{Recalling the definition of $\M_{3}$ from Section~\ref{sec:notation}, Theorem~\ref{thm:prelimit:idle:bound} can be extended to all $h \in \M_{3}$, at the cost of more complicated constants in the upper bound. Indeed, the machinery in Section~\ref{sec:stein:factor:bounds} is developed for generic test functions, except for the third-order Stein factor bounds in Lemma~\ref{lem:3bounds:stein:factors}. In the interest of space, we specialize that lemma to $h(x)=x$ to reduce the  algebra required. Nevertheless, Section~\ref{sec:stein:factor:bounds} contains all the ingredients needed to establish the corresponding bounds for general $h \in \M_{3}$.}

\rr{The restriction to $h \in \M_{3}$, rather than $h \in \lipone$ as in Theorem~\ref{thm:workload:1}, is not intrinsic to the prelimit approach. It stems from the way we establish the Stein factor bounds. In particular, our synchronous-coupling argument requires more regularity of $h(x)$ than the exponential Stein factor bounds in Lemma~\ref{lem:exp:stein:factors}, which are not proved using couplings. Such additional smoothness requirements are common in coupling-based Stein factor arguments; see, for example, \cite{GorhMack2016}.}

\rr{
A second reason for specializing Theorem~\ref{thm:prelimit:idle:bound} to $h(x)=x$, specific to the $G/G/1$ system, is that this choice lets us use the classical identity \eqref{eq:idle:expression} to bound the difficult boundary term from Theorem~\ref{thm:workload:1}, namely $\E \bar I^2/(2\E \bar I)$, in terms of simple model primitives. In this way, the two theorems reinforce one another: Theorem~\ref{thm:prelimit:idle:bound} controls the idle-period term, and Theorem~\ref{thm:workload:1} then converts that control into an explicit Wasserstein bound that applies to all $h\in\lipone$.
}

\rr{
Indeed, using the first equality in \eqref{eq:idle:expression} and the identity $\E Y_{2}=\lambda \E(S-U)^2/2$, we obtain, under the assumptions of Theorem~\ref{thm:prelimit:idle:bound},
\begin{align*}
    \frac{\E \bar I^2}{2 \E \bar I}    = \frac{ 1}{1-\rho}\E Y_2  - \frac{1}{\rho(1-\rho)} \E X +   \frac{\E S^2}{2\E S}\leq&\  \frac{1}{1-\rho}  ( \E Y_{2} -  \E X) +  \frac{\E S^2}{2\E S}  \leq \sum_{i=1}^{5} K_{i} +  \frac{\E S^2}{2\E S}.
\end{align*}
Since $\E \bar I^2/(2\E \bar I)=\E(R_a1(X=0))/(1-\rho)$, the bound on $\abs{\epsilon_{0}(f_{h})}$ in Theorem~\ref{thm:workload:1} can therefore be strengthened to
\begin{align*}
    \abs{\epsilon_{0}(f_h)} \leq&\ \delta \rho \lambda \E U^2/2 + \rho \delta \Big(\sum_{i=1}^{5} K_{i} +  \frac{\E S^2}{2\E S} \Big), \quad h \in \lipone. 
\end{align*}  
}

\subsubsection*{Comparing the two approaches.}
\rr{
Determining which of the bounds in Theorems~\ref{thm:workload:1} and~\ref{thm:prelimit:idle:bound} is tighter is not straightforward, especially because the prelimit bound is algebraically involved. Rather than ranking the bounds across parameter regimes, it is more useful to compare the qualitative information they provide. Theorem~\ref{thm:workload:1} gives a relatively clean bound, but leaves a boundary term that is difficult to control; analogous boundary terms in more complex models may be even harder to bound. By contrast, Theorem~\ref{thm:prelimit:idle:bound} gives a fully explicit bound under structural assumptions on the interarrival distribution. Although its constants are more complicated, they depend only on primitive quantities such as low-order moments and an upper bound on the interarrival hazard rate.
}
 
\rr{
The classical generator approach follows a relatively standard sequence of steps. Although it uses Palm inversion, that step is largely mechanical and changes little across examples. The resulting algebra can be lengthy, but it is not expected to be the main obstacle in applying the method to other models. Another advantage is that the approach relies on the Poisson equation of the limiting diffusion, whose Stein factor bounds may already be available. This is especially true for one-dimensional approximations, such as the exponential limits in the workload and JSQ examples, or the normal limit in the $G/M/\infty$ example. For multidimensional approximations, however, Stein factor bounds are much harder to obtain; for instance, no Stein factor bounds are currently available for the two-dimensional reflected Brownian motion that approximates the tandem queue in Appendix~\ref{app:tandem}.
}
 
\rr{
The prelimit approach has the opposite profile: it bypasses boundary terms, but requires Stein factor bounds for the prelimit Poisson equation. These bounds are much more model-specific than their diffusion counterparts and, as the workload example suggests, require a detailed understanding of the underlying queueing dynamics. Moreover, while the approach extends naturally to models with a single general clock, extending it to models with multiple general clocks remains open.
}

\rr{
The choice between the two approaches therefore depends on where the main difficulty lies for the model at hand: controlling the boundary terms or establishing the prelimit Stein factor bounds. The $G/G/1$ workload example illustrates precisely this tradeoff.
}

\section{The G/G/1 workload: the prelimit approach}
\label{sec:prelimit}


\rr{This section develops the prelimit generator approach for the $G/G/1$ workload and uses it to prove Theorem~\ref{thm:prelimit:idle:bound}.   After introducing some additional notation and providing a roadmap of the technical approach, the section proceeds in two main steps. First, we define the workload Poisson equation and derive a BAR-like identity suitable for generator comparison. Second, we expand this identity and compare it with the generator of the exponential approximation. The resulting expression, Lemma~\ref{lem:taylor:final}, writes the approximation error in terms of second- and third-order derivatives of the prelimit Poisson equation solution. The corresponding Stein factor bounds are developed separately in Section~\ref{sec:stein:factor:bounds}.  In Section~\ref{sec:conclusion} we discuss why extending the prelimit approach to models with multiple general clocks is substantially more difficult, and leave this extension as an open problem.
}
\subsection{Additional notation}
In addition to the  notation introduced in Section~\ref{sec:model:workload}, let $G(x) = \Prob(U \leq x)$ be the interarrival distribution function and  assume for  simplicity that $\Prob(U>0)=1$. Let 
\begin{align}
    &\bar B = \text{the duration of a busy period  initialized by an arrival to an empty system, and } \notag \\
      &\bar I = \text{the duration of the subsequent idle period}.  \label{eq:BI}
\end{align}
We  let $B_{0}$ be the length   of the initial busy period, with the convention that $B_{0} = 0$ if $X(0) = 0$, and let  $B_{n}$, $n \geq 1$,   be the lengths of the subsequent busy periods, which are i.i.d.\   $\bar B$. We also let $I_{n}$, $n \geq 0$, be the duration of the idle period following $B_{n}$, and note that $I_{n}$, $n \geq 1$, are i.i.d.\  $\bar I$. 

When $\rho < 1$,  \citep[Propositions~X.1.3 and X.3.1]{Asmu2003} say that
\begin{align}
    \E \bar I =  \frac{1-\rho}{\rho} \E \bar B < \infty. \label{eq:bari:barb}
\end{align} 
Furthermore, note that 
\begin{align*}
     \E \bar B = \lim_{\epsilon \downarrow 0} \E_{0,\epsilon} B_{0} = \E(\E_{\delta S, U} B_{0}),
\end{align*}
where $\E_{x,r_{a}}(\cdot) = \E (\cdot | Z(0)=(x,r_{a}))$ and the outer expectation on the right-hand side is with respect to the  distributions of $U$ and $S$.

 Lastly, we prove the following lemma in Appendix~\ref{app:prelim:proofs}. It shows  that the expected time to regeneration given any initial condition $(x,r_a)\in\mathbb S$ is finite under the stability condition
\eqref{eq:stability}. 
\begin{lemma}
    \label{lem:finite:busy:period}
If \eqref{eq:stability} holds then 
\begin{align}
\E_{x,r_a} B_0
\leq
\E \big(\E_{x+\delta S,U} B_0\big)
<\infty .
\label{eq:busy:finite}
\end{align}
\end{lemma}

\subsection{A technical roadmap}
\rr{We begin with an \emph{informal} outline of the prelimit approach and its key differences from the approach in Section~\ref{sec:workload}. Given $h : \R \to \R$ such that $\E|h(X)| < \infty$ and
\begin{align}
    \int_{0}^{\infty} \abs{\E_z h(X(t)) - \E h(X)} dt < \infty, \quad z \in \mathbb{S}, \label{eq:integrability}
\end{align}
define  
\begin{align*}
    F_h(z) = \int_{0}^{\infty} \big( \E_z h(X(t)) - \E h(X) \big) dt, \quad z \in \mathbb{S}.
\end{align*}
We immediately draw the reader’s attention to the boundary behavior of this function:
\begin{align}
    \lim_{\epsilon\to0} F_h(x + \delta\epsilon, r + \epsilon) = \E F_h(x + \delta S, U) 
    \quad \text{and} \quad
    \lim_{\epsilon\to0} F_h(\delta\epsilon, r + \epsilon) = \E F_h(\delta S, U). \label{eq:jumpzero}
\end{align}
For example, the first equality follows from
\begin{align*}
    F_h(x + \delta\epsilon, r + \epsilon)
    =
    \int_{0}^{\epsilon} \big( \E_z h(X(t)) - \E h(X) \big)dt
    +
    \E F_h(x+\delta S,U).
\end{align*}
Thus, if $F_h(z)$ satisfies the conditions of Lemma~\ref{lem:full:BAR:workload}, one may show that the boundary term in the BAR vanishes by \eqref{eq:jumpzero}, and the BAR reduces to
\begin{align*}
    0
    =
    -\delta \E \big(1(X>0)\partial_x F_h(Z)\big)
    -
    \E \big(\partial_{r_a}F_h(Z)\big).
\end{align*}
In this sense, $F_h(Z)$ is characterized entirely by the interior behavior of the workload process.
}

\rr{
In the remainder of this section, we outline the technical implementation of our approach. The idea is to vary the test function $h(x)$ over a sufficiently rich class so that the resulting boundary-term-free BAR fully characterizes the law of $Z$. We accomplish this by first showing that  $F_h(z)$ solves the Poisson equation
\begin{align*}
    -\delta 1(x>0)\partial_xF_h(z)-\partial_{r_a}F_h(z)
    =
    \E h(X)-h(x), \quad z=(x,r_a)\in\mathbb S.
\end{align*}
The left-hand side depends on both $x$ and $r_a$, whereas the right-hand side depends only on $x$. Thus, it is not immediate how to extract the diffusion generator, which approximates only the workload coordinate $X$, by a Taylor expansion. The absence of $r_a$ from the right-hand side also suggests some freedom in choosing the residual-time coordinate.
}

\rr{
Our approach is to set $r_a=R_a$ and take expectations. As we show, this gives
\begin{align}
    -\delta\partial_x\E F_h(x,R_a)
    +
    \lambda\E\big(F_h(x+\delta S,U)-F_h(x,U)\big)
    =
    \E h(X)-h(x), \quad x\geq 0. \label{eq:roadmap:poisson}
\end{align}
Before expanding $F_h(x+\delta S,U)-F_h(x,U)$ to derive the diffusion generator, we must relate expectations involving $U$ to those involving $R_a$. This step is unnecessary when all clocks are exponential, since then $U\stackrel{d}{=}R_a$. With general clocks, however, this relationship is the main new technical issue, and it becomes the key obstacle in extending the prelimit approach to models with two or more general clocks.
}

\rr{
It is useful to contrast this strategy with those of \cite{BravDaiMiya2017} and \cite{BravDaiMiya2024}. Those papers also rely on carefully engineered test functions that make the boundary terms in the BAR vanish. The resulting boundary-term-free BAR is then used to characterize the heavy-traffic limit of the corresponding queueing model.
}

\rr{
The key difference lies in the choice of test function. In \cite{BravDaiMiya2017} and \cite{BravDaiMiya2024}, the authors use an exponential test function, augmented so that the jump terms vanish. In contrast, our test function $F_h(z)$ is intrinsic to the queueing model: it is the solution to the prelimit Poisson equation. Moreover, as $h(x)$ varies, $F_h(z)$ gives an entire class of test functions rather than a single specially tailored one.
 }

 \subsection{The workload Poisson equation.}  
\label{sec:mhorizon}
In this section we prove that \eqref{eq:roadmap:poisson} holds for all $h \in \lipone$, which paves the way for the generator extraction and comparison step performed in Section~\ref{sec:taylor}. We begin with a discussion of how to verify \eqref{eq:integrability}. 

As one option, \cite{DaiMeyn1995}  show  that  when $\E U^{p+1} < \infty$ and $\E S^{p+1}<\infty$,  then  $\abs{ \E_{z}  Q(t)  - \E Q }$ decays at a rate of $1/t^{p-1}$, where $Q(t)$ and $Q$ are the customer count at time $t$ and in steady state, respectively. Their result holds for  queueing-network models beyond the $G/G/1$ system, but we wish to avoid using their complex machinery (in particular we would need to convert their bound in $\abs{ \E_{z}  Q(t)  - \E Q }$ into a bound on $\abs{ \E_{z}  X(t)  - \E X }$). 

A second option is to  verify \eqref{eq:integrability}  by noticing that 
\begin{align*}
    \int_{0}^{\infty} \abs{ \E_{z}  h(X(t)) - \E h(X) }dt \leq  \int_{0}^{\infty} \E \big(\abs{ \E_{z}  h(X(t)) - \E_{Z} h(X(t)) } \big)  dt,
\end{align*}
where the outer expectation is taken with respect to the stationary distribution of $Z$. To bound the right-hand side we would need to couple the workload process starting from $z$ with one starting from $Z$, and bound the expected coupling time (in terms of $z$ and $Z$). Constructing such a coupling is complicated by the fact that   $z$ and $Z$ may differ both in the initial workload  and residual interarrival time. 

In this paper we illustrate a third option that bypasses the need to verify \eqref{eq:integrability} directly. It involves   using the truncated
\begin{align*}
F^{M}_h(z) =&\ \int_{0}^{M} \big( \E_{z}  h(X(t)) - \E h(X) \big)  dt,\quad z \in \mathbb{S},
\end{align*}
and may be useful in other models where verifying \eqref{eq:integrability} directly is challenging. Our starting point is the following proposition, which is  proved in Appendix~\ref{app:mhorproofs}.
\begin{proposition}
\label{prop:Mpoisson}
For any $h \in \lipone$ and almost all $M > 0$,
\begin{align}
  -\delta \E \partial_{x} F_{h}^{M}(x,R_{a}) + \lambda \E \big(F_{h}^{M}(x+\delta S,U) -  F_{h}^{M}(x,U) \big) =&\ \E \big( \E_{x,R_{a}} h(X(M)) - h(x) \big). \label{eq:Mpoisson}
\end{align}
\end{proposition} 
We wish to take $M \to \infty$ in \eqref{eq:Mpoisson} and recover \eqref{eq:roadmap:poisson}, but  since we do not assume that $F_{h}(z)$ is well defined,  we first need to specify what we mean by both $\partial_{x} \E  F_h(x,R_{a})$ and $\E \big(F_h(x + \delta S,U) - F_h(x,U) \big)$. 

We first define 
\begin{align}
F_h(x+\epsilon, r_{a}) - F_h(x,r_{a}) =&\ \int_{0}^{\infty} \big( \E_{x+\epsilon,r_{a}}  h(X(t)) - \E_{x,r_{a}}  h(X(t)) \big)dt, \quad x \geq 0,\ \epsilon,r_{a} > 0. \label{eq:abuse1}  
\end{align} 
To argue that the right-hand side is well defined, we now introduce a simple  synchronous coupling  of the workload process that differs only in the initial workload and not in the initial residual time. This coupling also plays a central role in  Section~\ref{sec:stein:factor:bounds}.

Given $\epsilon > 0$, let $\{Z^{(\epsilon)}(t) = (X^{(\epsilon)}(t),R_{a}(t)) : t \geq 0\}$ be a coupling of $\{Z(t)= (X(t),R_{a}(t)): t \geq 0\}$ with  initial condition $ X^{(\epsilon)}(0) = X(0)+\epsilon$. Both systems share the same arrival process, and the service time of each arriving customer is identical in both systems. Similar to $B_0$, we define $B_{0}^{(\epsilon)} = \inf\{t \geq 0: X^{(\epsilon)}(t) = 0\}$.  It follows that for every sample path, 
\begin{align}
 \frac{d}{d t} \big( X^{(\epsilon)}(t) - X(t)\big) =&\ -\delta 1(X(t)=0)\,1(t\le B_0^{(\epsilon)}), \notag \\
Z^{(\epsilon)}(t) =&\ Z(t) \text{ for $t > B_{0}^{(\epsilon)}$}.  \label{eq:couplinggap} 
\end{align} 
 We adopt the convention that $\E_{x,r_{a}}(\cdot)$  is the expected value conditional on $Z(0) = (x,r_{a})$, even if the quantity inside the parentheses is a function of $Z^{(\epsilon)}(t)$. 
 
Returning to \eqref{eq:abuse1}, our synchronous coupling yields
\begin{align}
  \int_{0}^{\infty} \big(\E_{x+\epsilon,r_{a}} h(X(t)) - \E_{x,r_{a}} h(X(t)) \big) dt  = \E_{x,r_{a}} \int_{0}^{B^{(\epsilon)}_{0}    }  \big( h(X^{(\epsilon)}(t)) - h(X(t)) \big) dt. \label{eq:findiff}
\end{align}
The right-hand side is well defined because $\abs{h(X^{(\epsilon)}(t)) - h(X(t))} \leq \epsilon \norm{h'}$ and $\E_{x,r_{a}}B^{(\epsilon)}_{0} = \E_{x+\epsilon,r_{a}}B_{0} < \infty$ by \eqref{eq:busy:finite}. A similar line of reasoning yields the following two lemmas. The detailed proofs are found in  Appendix~\ref{app:poissonderive}.
\begin{lemma}
    \label{lem:abuse:d1}
    Let $T \geq 0$ be any random variable and define $ \partial_{x}  \E  F_{h}(x,T)$ by 
\begin{align}
 \partial_{x}  \E  F_{h}(x,T) =&\  \lim_{\epsilon \to 0} \frac{1}{\epsilon} \E \big( F_h(x+\epsilon, T) - F_h(x,T) \big), \quad x \geq 0, \label{eq:abuse3}
\end{align}
with the convention that $\partial_{x}  \E  F_{h}(x,T)  = \partial_{x}   F_{h}(x,r_{a}) $ when $T=r_{a}$ is deterministic. Then for any $h \in \lipone$, 
\begin{align}
 \partial_{x}  \E  F_{h}(x,T) =&\   \E \Big(\E_{x,T} \int_{0}^{B_{0}}    h'(X(t))  dt\Big), \quad x \geq 0, \label{eq:d1}
\end{align}
where the outer expectation is with respect to $T$. 
In the special case that $T=r_{a}$ is deterministic, \eqref{eq:d1} yields an expression for $\partial_{x}   F_{h}(x,r_{a})$. In particular, this implies that  for any random $T \geq 0$,
\begin{align}
    \partial_{x}  \E  F_{h}(x,T) =  \E  \partial_{x}  F_{h}(x,T), \quad x \geq 0. \label{eq:d1:interchange}
\end{align} 
\end{lemma}

\begin{lemma}
\label{lem:infpoisson}
For any   $h \in \lipone$ and  $x \geq 0$,  
\begin{align}
\lim_{M \to \infty} \E \big( \E_{x,R_{a}} h(X(M)) \big) =&\ \E h(X),  \label{eq:lim1}  \\
\lim_{M \to \infty}  \E \partial_{x} F_{h}^{M}(x,R_{a})  =&\  \E \partial_{x} F_{h}(x,R_{a}), \label{eq:lim2}\\
\lim_{M \to \infty}  \E \big(F_{h}^{M}(x+\delta S,U) -  F_{h}^{M}(x,U) \big) =&\ \E \big(F_{h}(x+\delta S,U) -  F_{h}(x,U) \big).\label{eq:lim3} 
\end{align}  
\end{lemma}
Starting from \eqref{eq:Mpoisson} of Proposition~\ref{prop:Mpoisson}, we apply  Lemma~\ref{lem:infpoisson}   to take $M \to \infty$ there. Since   $\partial_{x}  \E  F_{h}(x,R_{a}) =  \E  \partial_{x}  F_{h}(x,R_{a})$  and $\partial_{x}  \E  F_{h}(0,R_{a}) = 0$ by Lemma~\ref{lem:abuse:d1}  and \eqref{eq:d1}, respectively, we arrive at 
\begin{align}
    -\delta    \partial_{x} \E  F_h(x,R_{a}) + \lambda \E \big(F_h(x + \delta S,U) - F_h(x,U) \big) =  \E h(X) - h(x), \quad x \geq 0.  \label{eq:the:poisson:abuse} 
\end{align}
In the following section, we replace $\E \big(F_h(x + \delta S,U) - F_h(x,U) \big)$ by a term where the expectation is taken over $R_{a}$ instead of $U$. We then perform a Taylor expansion to compare the left-hand side with the diffusion generator.

\subsection{Generator expansion and comparison}
\label{sec:taylor}
 Let $S'$ be an independent copy of $S$ and introduce the random variable 
 \begin{align*}
J(x,r_{a})  = - (x \wedge \delta r_{a}) + \delta S', \quad (x,r_{a}) \in \mathbb{S}. 
\end{align*}
Conditioned on $Z(0)=(x,r_{a})$,  $J$ represents the jump in workload between time $t=0$  and $t = r_{a}$ (the first arrival time). 
We present the  following lemma, which is proved in Appendix~\ref{app:taylorproofs}.
\begin{lemma}
\label{lem:relate}
For any $(x,r_{a}) \in \mathbb{S}$, $s \geq 0$, and  $h \in \lipone$,
\begin{align*}
F_h(x+ \delta s,r_{a}) - F_h(x,r_{a}) =&\ \E \big( F_h(x+\delta s + J(x,r_{a}),U) - F_h(x+J(x,r_{a}),U) \big) + \epsilon(x,r_{a},s).
\end{align*}
The $U$ on the right-hand side is independent of $S'$ (and thus $J(x,r_{a})$) and 
\begin{align*}
\epsilon(x,r_{a},s) =&\  \E^{S'} \Big(  \int_{-x \wedge(\delta r_{a})}^{-(x+\delta s) \wedge(\delta r_{a})} \partial_{x}   \E^{U} F_h(x+\delta s +v +\delta S',U) dv\Big) \\
&+ \int_{0}^{r_{a}} \Big(  h\big((x + \delta s-\delta t)^{+}\big) -  h\big((x  -\delta t)^{+}\big) \Big) dt,
\end{align*}
where we recall that $\E^{U}(\cdot)$ and $\E^{S'}(\cdot)$ denote  expectations with respect to $U$ only and $S'$ only, respectively. 
\end{lemma}
\noindent To simplify notation, we define
\begin{align}
\bar F_h'(x) = \partial_{x}  \E F_h(x,R_{a}). \label{eq:fbar}
\end{align}
Replacing $x$ by $x + J(x,R_{a})$ in the Poisson equation \eqref{eq:the:poisson:abuse} and taking expectations  yields 
\begin{align}
& \E h(X) - \E h(x + J(x,R_{a}))  \notag \\
=&\ -\delta \E \bar F_h'(x + J(x,R_{a})) + \lambda \E \big(F_h(x+ \delta S,R_{a}) - F_h(x,R_{a}) \big) - \E \big(\epsilon(x,R_{a},S)\big). \label{eq:backpoisson}
\end{align}
We emphasize that   $\E \bar F_h'(x + J(x,R_{a}))$ is actually  $\E^{R_{a}}  \partial_{x} \E^{R_{a}'} F_h(x + J(x,R_{a}),R_{a}')$, where $R_{a}$ and $R_{a}'$ are independent copies.

The following lemma is the result of  expanding the right-hand side of \eqref{eq:backpoisson} and taking expected values with respect to the steady-state of the approximating diffusion.  It is proved in Appendix~\ref{app:taylorlemma}.
\begin{lemma}
\label{lem:taylor:final}
Let $Y_{2}$ be exponentially distributed with mean $\lambda \E(U-S)^2/2$.      Fix $h \in \lipone$ and assume that  $\bar F_{h}''(x)$ and $\bar F_{h}'''(x)$ exist for all $x \geq 0$, and that  $\E \abs{\bar  F_h'(Y_{2})}, \E \abs{\bar F_h''(Y_{2})} < \infty$. Then 
    \begin{align*}
        \E h(Y_{2}) - \E h(X)  
        =&\  \E  \big(    h(Y_{2}) -  h(Y_{2} + J(Y_{2},R_{a}))\big)   + \E \big(\epsilon(Y_{2},R_{a},S)\big) \notag  \\
 &+\delta \E \Big( 1(\delta R_{a} \geq Y_{2}) \big( \bar F_{h}'(\delta S )  - \bar F_h'(Y_{2}) - \delta(S-R_{a})\bar F_h''(Y_{2})\big) \Big) \\
 & +\delta \E \Big(  1(\delta R_{a}  < Y_{2}) \int_{0}^{\delta(S -R_{a} )} \int_{0}^{v} \bar F_h'''(Y_{2}+u) du  dv \Big) \\
 &- \lambda \E \int_{0}^{\delta S} (\delta S - v)\int_{0}^{v} \bar F_h'''(Y_{2}+u) du dv. 
    \end{align*} 
    \end{lemma}
Lemma~\ref{lem:taylor:final} expresses $\E h(Y_{2}) - \E h(X)$ in terms of the difference of the diffusion and prelimit generators, and  is the prelimit-approach analog of Proposition~\ref{prop:extraction:workload}. 


\section{Stein factor bounds for the G/G/1 workload process}
\label{sec:stein:factor:bounds}
In this section we derive second- and third-order Stein factor bounds and then combine them with Lemma~\ref{lem:taylor:final} to prove Theorem~\ref{thm:prelimit:idle:bound}. 
Although we focus on the case $h(x)=x$, the same approach extends to Lipschitz test functions with Lipschitz first and second derivatives.

Throughout, we assume that $U$ is absolutely continuous with density $G'(x)$ and hazard rate
\begin{align*}
    \eta(x) = \frac{G'(x)}{1-G(x)}, \quad x \geq 0,
\end{align*}
and write $\overline{\eta} = \sup_{x \geq 0} \eta(x)$. 
Extending the analysis to interarrival distributions with point masses would require additional work; see the discussion following Lemma~\ref{lem:stein2} in Section~\ref{sec:2bounds}.

Our proof of Theorem~\ref{thm:prelimit:idle:bound} proceeds in stages. 
We first establish auxiliary results containing Stein factor bounds for the workload process. 
These appear in two parts: Lemma~\ref{lem:2bounds:stein:factors}, proved in Section~\ref{sec:2bounds}, and Lemma~\ref{lem:3bounds:stein:factors}, proved in Section~\ref{sec:3bounds}. 
Only the latter requires $\eta(x)$ to be nonincreasing; the proof can be adapted to more general hazard rate structures, but at the expense of added complexity.
Throughout, recall that $\bar B$ denotes the busy-period duration as defined in \eqref{eq:BI}.

\rr{\begin{lemma}
    \label{lem:2bounds:stein:factors}
    Assume that $\overline{\eta} < \infty$. The following bounds hold for   $x \geq 0$. For any $h \in \lipone$, 
    \begin{align}
        \abs{\delta\bar F_{h}'(x)} \leq   x( 1 + (\lambda + \overline{\eta}) \E \bar B ) \quad \text{ and } \quad \abs{ \delta\partial_{x} \E F_{h}(x,U)} \leq  x( 1 +  2\overline{\eta} \E \bar B  ),   \label{eq:1bounds}
    \end{align}
    and for any $h \in \M_{2}$,    
    \begin{align}
        \abs{\delta \bar F_{h}''(x)} \leq   (1+x)\big( 1 + (\lambda + \overline{\eta}) \E \bar B  \big)  \quad \text{ and } \quad  \abs{\delta \partial_{x}^2  \E  F_{h}(x,U)} \leq (1+x)(1 + 2\overline{\eta} \E \bar B). \label{eq:2bounds} 
    \end{align}
  In the special case that $h(x) = x$,
    \begin{align}
          \abs{\delta \bar F_{h}''(x)} \leq    1 + (\lambda + \overline{\eta}) \E \bar B      \quad \text{ and }   \quad    \abs{\delta \partial_{x}^2  \E  F_{h}(x,U)} \leq  1 + 2\overline{\eta} \E \bar B. \label{eq:2boundsx} 
    \end{align} 
\end{lemma}}

\begin{lemma}
    \label{lem:3bounds:stein:factors}
    Suppose that $\overline{\eta} < \infty$, that $\eta(x)$ is nonincreasing, and fix $h(x) =x$. Then for any $x \geq 0$,
\begin{align*}
    \abs{\delta^2 \bar F_{h}'''(x)} \leq&\    \Prob(\delta U > x)  3\lambda \overline{\eta} \E \bar B   +  \Prob(\delta U < x  < \delta\bar I +\delta U) \lambda \overline{\eta} (1 + \rho + \lambda^2 \E U^2)  \E \bar B,
\end{align*}
where $U$ and $\bar I$ are independent.
\end{lemma}
\rr{Unlike the bounds in Lemma~\ref{lem:2bounds:stein:factors}, the third-order bounds in Lemma~\ref{lem:3bounds:stein:factors} are stated only for the special case $h(x)=x$. This restriction is made for brevity and to limit the length of the paper. In the proof of the lemma (Section~\ref{sec:3bounds}), we indicate how the argument can be extended to general test functions $h \in \M_{3}$. Doing so requires only additional algebraic manipulations, which are straightforward but would significantly increase the length of the exposition.}

The Stein factor bounds in Lemmas~\ref{lem:2bounds:stein:factors} and \ref{lem:3bounds:stein:factors} both involve the busy-period duration $\E \bar B$. 
Thus, in order to make the bounds explicit we require an upper bound on $\E \bar B$ in terms of the interarrival and service distributions. 
A classical result due to \cite{Loul1978} provides such a bound:
\begin{align}
    \E \bar B \leq  0.9 \rho \frac{\sqrt{\text{Var}(U-S)}}{1-\rho} 
    \exp\bigg(   5.4\frac{\E \abs{U-S}^{3} }{ (\text{Var}(U-S) )^{3/2}} 
    + 0.8 \frac{\E (U-S) }{ \sqrt{\text{Var}(U-S) }}\bigg).  
    \label{eq:loulou:bound}
\end{align}
Notably, the right-hand side depends only on the first three moments of $S$ and $U$. 
We next present an alternative bound, derived in Appendix~\ref{sec:busy:period:bound}. 
\begin{lemma}
\label{lem:loulou}
Recall that $I_1, I_2, \ldots$ are i.i.d.\ $\bar I$. For any $\rho < 1$, 
\begin{align}
\E V = \lambda \frac{\E S^2}{2} +  \lambda\E \bar B\Big[  \E ( S - U)^{+} +   \sum_{k=2}^{\infty}    \E  \Big( S - U -   \sum_{i=1}^{k-1} I_{i} \Big)^{+}\Big].  \label{eq:loulou}
\end{align}
Let $W$ denote  the steady-state customer waiting time. A consequence of \eqref{eq:loulou} is that  
\begin{align}
    \E \bar B = \frac{\E W \E S}{ \E ( S - U)^{+} +   \sum_{k=2}^{\infty}    \E  \Big( S - U -   \sum_{i=1}^{k-1} I_{i} \Big)^{+}} \leq \frac{\rho\text{Var}(S-U)}{2 (1-\rho) \E ( S - U)^{+}}. \label{eq:barBbound}
\end{align}
\end{lemma}  
\begin{remark}
    The bound  in \eqref{eq:barBbound} is a result of the upper bound on $\E W$ due to \cite{King1962}. Tighter bounds on $\E W$ have been established since then \citep{Daleetal1992}, and any one of them could be used instead. 
\end{remark}
We require one final auxiliary lemma, proved in Appendix~\ref{app:final:aux:proof},  which uses the Stein factor bounds  from Lemmas~\ref{lem:2bounds:stein:factors} and \ref{lem:3bounds:stein:factors} to bound the right-hand side of the Taylor expansion in Lemma~\ref{lem:taylor:final}.  We recall from Lemma~\ref{lem:relationships:workload} that $\E R_{a}^{k} = \lambda \E U^{k+1}/(k+1)$.
\begin{lemma}
    \label{lem:final:aux}
    Let $h(x) = x$. Then 
    \begin{align}
    &\E  \abs{ h(Y_2 + J(Y_2,R_{a})) - h(Y_2) } \leq \delta ( \lambda \E U^2/2 +\E S), \label{eq:ineq:1}\\
    &\E \abs{\epsilon(Y_2,R_{a},S)} \leq   \delta \big(\rho \E U^2/2  + \delta (1/\E Y_2) \lambda (\E U^2/2)(\E S^2 + (\E S)^2)( \delta + 2 \overline{\eta}  \delta \E \bar B )\big), \label{eq:ineq:2} \\
    &\delta \E \abs{ 1(\delta R_{a} \geq Y_2) \big( \bar F_{h}'(\delta S )  - \bar F_h'(Y_2) - \delta(S-R_{a})\bar F_h''(Y_2)\big)}  \notag \\
    &\hspace{2cm} \leq \delta \big(  (1/\E Y_2) \big( \rho \E U^2 + 5\lambda \E U^3/3 \big) +  \lambda \E U^3/3 \big) \big( \delta + (\lambda + \overline{\eta}) \delta \E \bar B  \big) ,\label{eq:ineq:3}  \\
    & \lambda  \Big| \E \int_{0}^{\delta S} (\delta S - v)\int_{0}^{v}  \bar F_{h}'''(Y_2+u)      du dv \Big| \leq   \delta  \lambda    \E  S^3 C  , \label{eq:aineq:1}   \\
     &\delta \Big| \E  1(\delta R_{a}  < Y_2) \int_{0}^{\delta(S -R_{a} )} \int_{0}^{v} \bar F_h'''(Y_2+u) du  dv \Big|  \leq  \delta \E (S^2 - \rho U^2 + \lambda \E U^3/3)  C, \label{eq:aineq:2}  
\end{align} 
where $C = (3   + (\E U + \E \bar I) \lambda  (1 + \rho + \lambda^2 \E U^2))   (1/\E Y_2) \overline{\eta} \delta \E \bar B$.
\end{lemma}

 \startproof{Proof of Theorem~\ref{thm:prelimit:idle:bound}}
We apply the bounds in Lemma~\ref{lem:final:aux} to the  expression for $\E h(X) - \E h(Y_2)$ in Lemma~\ref{lem:taylor:final}. The constants $K_1$ through $K_5$ correspond exactly to the upper bounds in \eqref{eq:ineq:1}--\eqref{eq:aineq:2}, respectively. The relationship between $\E \bar I$ and $\E \bar B$ was previously mentioned in \eqref{eq:bari:barb} and the bound on $\E \bar B$ follows from using the lesser of the two bounds in \eqref{eq:loulou:bound} and Lemma~\ref{lem:loulou}. 
\finishproof

\subsection{Second-derivative bound}
\label{sec:2bounds}
We now establish Lemma~\ref{lem:2bounds:stein:factors}.  
Let $\alpha(t)$ denote the age of the interarrival process (the backward recurrence time) at time $t$.  
Our first step is to identify an expression for the second derivative of $\E F_h(x,T)$ with respect to $x$, obtained by differentiating $\partial_x F_h(x,r_a)$.   We prove it in Appendix~\ref{app:stein2:proofs}.

\rr{\begin{lemma} \label{lem:stein2}
Suppose that $\overline{\eta} < \infty$. For any $h \in \M_{2}$, any absolutely continuous random variable $T \geq 0$ with bounded density $\theta(x)$,  and any $x \geq 0$, 
\begin{align*}
 \partial_{x}^2  \E  F_{h}(x,T) =&\  \partial_{x} \E F_{h'}(x,T) +  \frac{1}{\delta} h'(0) \\
 &+ \frac{1}{\delta} \Big( \theta(x/\delta)  + \E \big(1(T<x/\delta)\E_{x,T} \eta\big(\alpha(B_{0})\big) \big) \Big) \E \big( \partial_{x} F_{h}(\delta S, U) \big),
\end{align*}
where $\partial_{x} \E F_{h'}(x,T)$ is as in Lemma~\ref{lem:abuse:d1} but with $h'(x)$ instead of $h(x)$. In the special case that $h(x)=x$, 
\begin{align*}
    \partial_{x}^2  \E  F_{h}(x,T) =&\    \frac{1}{\delta}  + \frac{1}{\delta} \big( \theta(x/\delta)  + \E \big(1(T<x/\delta)\E_{x,T} \eta(\alpha(B_{0})) \big) \big) \E \bar B.
\end{align*}
\end{lemma} }
 
\begin{remark}
\label{rem:4}
In the proof of Lemma~\ref{lem:stein2}, the quantity 
\begin{align*}
    \lim_{\epsilon \to 0} \frac{1}{\epsilon}\Prob_{x,r_a}(R_{a}(B_{0}) < \epsilon/\delta)
\end{align*}
arises. Since $R_{a}(B_{0}) = I_{0}$, this is the density at zero of the idle period $I_{0}$.  
When $U$ has a density, the latter equals $\E_{x,r_a} \eta(\alpha(B_{0}))$, which will play an important role in the third-derivative bounds.  
If $U$ has point masses  one cannot use the hazard-rate connection, which makes the analysis more challenging. 
\end{remark}

\rr{
\startproof{Proof of Lemma~\ref{lem:2bounds:stein:factors}}
Setting $\widehat{h}(x) =x$, let us first prove \eqref{eq:2boundsx}. Let $x \geq 0$.  Since the densities of $R_{a}$ and $U$ are bounded by $\lambda$ and $\overline{\eta}$, respectively, the expression for $\partial_{x}^2  \E  F_{\widehat{h}}(x,T)$ in Lemma~\ref{lem:stein2} yields 
\begin{align*}
  \delta \abs{\partial_{x}^2  \E  F_{\widehat{h}}(x,R_{a})} \leq&\    1 +  \big( \lambda  +  \overline{\eta} \big) \E \bar B  \quad \text{ and }  \quad  \delta \abs{\partial_{x}^2  \E  F_{\widehat{h}}(x,U)} \leq    1 +  2\overline{\eta}  \E \bar B,
\end{align*}
which proves \eqref{eq:2boundsx}. Additionally, observe that
\begin{align}
    &\delta \E\big(\E_{x,R_{a}}B_{0}\big) = \delta  \partial_{x} \E F_{\widehat{h}}(x,R_{a}) = \delta \int_{0}^{x} \partial_{x}^2 \E F_{\widehat{h}}(u,R_{a})du \leq x \big(1 +  \big( \lambda  +  \overline{\eta} \big) \E \bar B\big), \notag \\
    &\delta \E\big(\E_{x,U}B_{0}\big) = \delta  \partial_{x} \E F_{\widehat{h}}(x,U) = \delta \int_{0}^{x} \partial_{x}^2 \E F_{\widehat{h}}(u,U)du \leq x \big(1 +  2\overline{\eta}  \E \bar B\big).  \label{eq:keyb0}
\end{align}
In both lines, the first  equality is due to \eqref{eq:d1} of Lemma~\ref{lem:abuse:d1}, and the second equality follows from $\partial_{x} \E F_{\widehat{h}}(0,R_{a})  = \partial_{x} \E F_{\widehat{h}}(0,U)  = 0$ together with  the fundamental theorem of calculus.
\\
The remaining bounds follow almost immediately. For any $h \in \lipone$, by \eqref{eq:d1} of Lemma~\ref{lem:abuse:d1}  we have
\begin{align*}
  \abs{\delta\bar F_{h}'(x)} = \abs{ \delta \partial_{x} \E F_{h}(x,R_{a})} \leq \delta \E \big( \E_{x,R_{a}} B_{0}\big) \quad \text{ and }  \quad  \abs{ \delta \partial_{x} \E F_{h}(x,U)} \leq \delta \E \big( \E_{x,U} B_{0}\big),
\end{align*}
which, when combined with  \eqref{eq:keyb0}, yields \eqref{eq:1bounds}. For the bounds in \eqref{eq:2bounds}, observe that  \eqref{eq:d1} of Lemma~\ref{lem:abuse:d1}  and our assumption that  $h \in \M_{2}$ yields
\begin{align}
    \abs{\partial_{x} \E F_{h'}(x,R_{a})} \leq \norm{h''} \E \big(\E_{x,R_{a}} B_0\big)  \leq \E \big(\E_{x,R_{a}} B_0\big). \label{eq:fh'}
\end{align}
Applying the bounds in  \eqref{eq:fh'} and \eqref{eq:keyb0} to the form of $\partial^{2}_{x} \E F_{h}(x,T)$ in Lemma~\ref{lem:stein2} yields \eqref{eq:2bounds}. 
\finishproof 
}

\subsection{Third-derivative bound}
\label{sec:3bounds}
In this section we prove Lemma~\ref{lem:3bounds:stein:factors}. Differentiating  both sides of the Poisson equation \eqref{eq:the:poisson:abuse} twice  yields
\begin{align}
   \delta \bar F_{h}'''(x) =&\ \lambda \partial_{x}^{2} \E \big(F_h(x + \delta S,U) - F_h(x,U) \big) + h''(x).  \label{eq:d3:initial}
\end{align}
 The following lemma is proved in Appendix~\ref{app:d3:form:proof}. It follows directly from Lemma~\ref{lem:stein2} after verifying that $\partial_{x}^{2} \E \big(F_h(x + \delta S,U) - F_h(x,U) \big) = \E^{S}\big(\partial_{x}^{2} \E^{U}  F_h(x + \delta S,U) - \partial_{x}^{2} \E^{U}  F_h(x,U)\big)$. 

 \rr{\begin{lemma}
    \label{lem:d3:form}
    Suppose that $U$ has a bounded density. Then for any $h \in \M_{3}$ and any $x \geq 0$, 
    \begin{align*}
         \delta^2 \bar F_{h}'''(x)  =&\ \lambda\delta \E \Big(\int_{0}^{\delta S} \partial_{x}^2 \E F_{h'}(x+y,U) dy \Big)  +  \lambda\E\big( G'(x/\delta +  S) - G'(x/\delta)\big)\E \big( \partial_{x} F_{h}(\delta S, U) \big) \\
 &+   \lambda\E \Big(1(x/\delta < U <x/\delta +  S)\E_{x+ \delta S,U} \eta\big(\alpha(B_{0})\big)\Big) \E \big( \partial_{x} F_{h}(\delta S, U) \big) \\
 &+ \lambda\E \Big( 1(U < x/\delta) \Big(\E_{x+ \delta S,U} \eta\big(\alpha(B_{0})\big) - \E_{x,U} \eta\big(\alpha(B_{0})\big)\Big) \Big) \E \big( \partial_{x} F_{h}(\delta S, U) \big) + \delta  h''(x). 
    \end{align*}
\end{lemma}}
\rr{To prove Lemma~\ref{lem:3bounds:stein:factors}, we make the simplifying assumption that $h(x) = x$. In this case, \eqref{eq:d1} of Lemma~\ref{lem:abuse:d1} implies that $\partial_{x}\E F_{h'}(\cdot,U)=0$ and $\E ( \partial_{x} F_{h}(\delta S, U) ) = \E \bar B$. As a result, the expression for $\delta^2 \bar F_{h}'''(x)$ simplifies to 
\begin{align}
         \delta^2 \bar F_{h}'''(x)  =&\   \lambda\E\big( G'(x/\delta +  S) - G'(x/\delta)\big)\E  \bar B \notag \\
 &+   \lambda\E \Big(1(x/\delta < U <x/\delta +  S)\E_{x+ \delta S,U} \eta\big(\alpha(B_{0})\big)\Big) \E  \bar B  \notag \\
 &+ \lambda\E \Big( 1(U < x/\delta) \Big(\E_{x+ \delta S,U} \eta\big(\alpha(B_{0})\big) - \E_{x,U} \eta\big(\alpha(B_{0})\big)\Big) \Big) \E  \bar B. \label{eq:spec:d3}
    \end{align}
One can use the generic expression for $\delta^2 \bar F_{h}'''(x)$ in Lemma~\ref{lem:d3:form} to extend the following arguments to any  $h \in \M_{3}$, at the cost of additional algebraic manipulations that are straightforward but would substantially increase the length of the exposition.}
     
 The first two terms in \eqref{eq:spec:d3} are straightforward to bound, whereas the third term
\begin{align*}
    \E \Big( 1(U < x/\delta) \Big(\E_{x+ \delta S,U} \eta\big(\alpha(B_{0})\big) - \E_{x,U} \eta\big(\alpha(B_{0})\big)\Big) \Big) \E \bar B 
\end{align*} 
requires substantially more effort. Naively bounding it by $\overline{\eta} \E \bar B$    is insufficient, since we need an additional factor of $\delta$ in the bound.
The following lemma is proved in Appendix~\ref{sec:lem4}, where we also describe how to establish an analogous result when the hazard rate is uniformly bounded away from zero.
\begin{lemma}
\label{lem:mixing:time}
Assume that $\overline{\eta} < \infty$ and that $\eta(x)$ is nonincreasing. For any   $x,s \geq 0$, and  $r_{a} < x/\delta$,
        \begin{align*}
            \abs{\E_{x+ \delta s,r_{a}} \eta\big(\alpha(B_{0})\big) - \E_{x,r_{a}} \eta\big(\alpha(B_{0})\big)} \leq \overline{\eta} (1 +  \lambda s + \lambda^2 \E U^2) \Prob(\bar I > x/\delta-r_{a}).
        \end{align*}
\end{lemma}
\rr{
The central quantity on the left-hand side is $\eta(\alpha(B_0))$: the hazard rate of $U$ evaluated at the interarrival age when the workload first hits zero. To prove the lemma, we relate $\alpha(B_0)$ to a renewal process driven by the idle periods of the workload process. Intuitively, the initial workload $V(0)$ is cleared once the cumulative idle time equals $V(0)$, so the difference between $\E_{x+\delta s,r_a}\eta(\alpha(B_0))$ and $\E_{x,r_a}\eta(\alpha(B_0))$ is governed by the mixing behavior of this idle-period renewal process. The assumption that $\eta(x)$ is nonincreasing is used precisely to control this mixing behavior: in the proof of Lemma~\ref{lem:mixing:time}, it yields an upper bound on the relevant coupling probability (see \eqref{eq:prob:couple}), which leads to the stated estimate. The lemma can be extended to other interarrival distributions whenever comparable control of the relevant mixing time is available.
}
\startproof{Proof of Lemma~\ref{lem:3bounds:stein:factors}}
We bound, one by one, each term in the expression  of $\delta^2 \bar F_{h}'''(x)$ in Lemma~\ref{lem:d3:form}. Observe that  
\begin{align*}
      &\lambda \abs{\E\big( G'(x/\delta +  S) - G'(x/\delta)\big)}\E \bar B  \leq   2\lambda \overline{\eta}  \Prob(U > x/\delta) \E \bar B, \\
    &\lambda\E \Big(1(x/\delta < U <x/\delta +  S)\E_{x+ \delta S,U} \eta\big(\alpha(B_{0})\big)\Big) \E \bar B \leq  \lambda \overline{\eta}\Prob(U > x/\delta) \E \bar B,
\end{align*}
where in the first inequality we used $G'(x) = \eta(x) \Prob(U > x) \leq \overline{\eta}\Prob(U > x) $.  Lastly we invoke Lemma~\ref{lem:mixing:time}, which implies that 
\begin{align*}
    &\lambda\E \Big( 1(U < x/\delta) \abs{\E_{x+ \delta S,U} \eta\big(\alpha(B_{0})\big) - \E_{x,U} \eta\big(\alpha(B_{0})\big)} \Big) \E \bar B \\
    \leq&\ \lambda \overline{\eta} (1 + \lambda \E S + \lambda^2 \E U^2) \E^{U} \big(1(U < x/\delta) \E^{\bar I} (1(\bar I > x/\delta-U)) \big) \E \bar B \\
    =&\  \lambda \overline{\eta} (1 + \rho + \lambda^2 \E U^2) \Prob(U < x/\delta < \bar I + U)\E \bar B,
\end{align*}
   where $U$ and $\bar I$ are independent.
\finishproof

\rr{Having developed the methodology for the one-dimensional $G/G/1$ workload, we now turn to an example with  multiple general clocks: the JSQ system.}

\section{The JSQ total customer count}
\label{sec:jsq}
Consider a parallel-server system with $n$ identical servers, each with its own buffer, operating under a first-come-first-served policy. Customer arrivals form a renewal process with generic interarrival time distribution $U$   and customer service times are i.i.d. $S$.  Upon arrival, customers are routed to the server with the shortest queue.  We  assume for   simplicity  that simultaneous arrivals/departures do not occur with probability one, and that  ties are broken uniformly at random (other tie-breaking rules can be handled as well). We set 
\begin{align*}
    &\lambda = 1/\E U, \quad \mu = 1/\E S, \quad \rho = \lambda/n\mu, \\
    &c_{U}^2 = \lambda^2 \text{Var}(U), \quad c_{S}^2 = \mu^2 \text{Var}(S),
\end{align*}
 and assume that  $\rho < 1$, $\E U^3 < \infty$, and $\E S^3 < \infty$. 
 
Let $Q_{i}(t)$ be the number of customers (both in service and waiting) at server $i$ at time $t \geq 0$. Also let $R_{a}(t)$ be the time until the next arrival and let $R_{s,i}(t)$ be the remaining service time of the customer being served by server $i$; when server $i$ is idle, $R_{s,i}(t)$ denotes the service time of the next customer to arrive to server $i$. 
Define   $Q(t) = (Q_{1}(t),\ldots, Q_{n}(t))$ and $R_{s}(t) = (R_{s,1}(t),\ldots, R_{s,n}(t))$, and let 
\begin{align*} 
    Z(t) =&\ (Q(t), R_{a}(t),R_{s}(t) ).  
\end{align*}
We assume that $\{Z(t): t \geq 0\}$ is positive Harris recurrent  (see \cite{Bram2011a} for sufficient conditions) and let $Z = (Q,R_a,R_s)$ be the vector having its stationary distribution. Define the (scaled) total customer count
\begin{align*}
    X(t) = \delta \sum_{i=1}^{n} Q_{i}(t) \quad \text{ and } \quad X = \delta \sum_{i=1}^{n} Q_{i},
\end{align*} 
where $\delta = (1-\rho)$. 

In addition to positive Harris recurrence, we assume that $\mathbb{E}X^{2}<\infty$. 
While one can establish this via a combination of the Palm inversion formula and the BAR --- first by applying the test function 
$f(z) = \sum_{i=1}^{n} (q_i \wedge M + \mu r_{s,i})^{2}$ and letting $M \to \infty$ to obtain $\mathbb{E}X < \infty$, 
and then combining this with the test function 
$f(z) = \sum_{i=1}^{n} (q_i \wedge M + \mu r_{s,i})^{3}$ to deduce $\mathbb{E}X^{2}<\infty$ --- the calculations are lengthy and would distract from the paper’s main focus, so we take $\E X^{2}<\infty$ as an assumption. The following is our main result. 
\begin{theorem}
    \label{thm:jsq:main}
    Assume that $\E X^2 < \infty$ and let $Y$ be exponentially distributed with mean $( \rho c_{U}^{2} +   c_{S}^{2})/2$. Then
    \begin{align*}
        d_{W}(X,Y) \leq&\  \frac{1}{n} \sum_{i=1}^{n} \E \Big(1(Q_i=0)   \Big(\sum_{j=1}^{n} Q_j + \lambda R_{a} + \sum_{j=1}^{n} \mu R_{s,j}\Big)   \Big) +   \delta (K_{0} + K_{A} + K_{D}),
    \end{align*}
    where
    \begin{align*}
        &K_{0} = \frac{1}{2}  \big(\lambda^2 \E U^2/2 +   ( \lambda \mu  \E S^2/2 + \mu \delta\E S)\big),\\
        &K_{A} = \frac{4}{\lambda c_{U}^{2} + n \mu c_{S}^{2}} \Big( \frac{1}{3} \lambda \E\abs{1-\lambda U}^{3} +   c_{U}^{2}  n\mu \lambda(\rho (\lambda/n) \E S^2/2 + \delta \E S+ \lambda \E U^2/2) \\
    & \hspace{8cm} + \lambda c_{U}^{2}   (1 + n + \lambda \mu \E U^2 +  \mu^2 \E S^2 )\Big), \\
    K_{D} =&\ \frac{4n}{\lambda c_{U}^{2} + n \mu c_{S}^{2}}\bigg( \frac{1}{3}   \frac{\lambda}{n} \E\abs{1 - \mu S}^{3} \\
    & +    c_{S}^{2}   \big(   (\lambda^2 \mu \rho + 2(n-1)\mu^2 \rho \lambda)   \E S^2/2n    + \delta(\lambda \mu + 2(n-1)\mu^2)\E S  + \frac{1}{2} \lambda^2 \mu \E U^2 \big)\\ 
    &+   \mu c_{S}^{2}     n \rho  + \mu c_{S}^{2}      (\lambda /n) \big( \E S+\lambda \E S^2 + \lambda^2 \E U^2 \E S \big)  \\
    &+  c_{S}^{2}  \lambda \big( 1 +    (\lambda \mu \E U^2/2 + 1/\rho)  + 2\mu^2 \E S^2  \big)\bigg) +\frac{1}{2}    c_{S}^{2} \lambda/n.
    \end{align*} 
\end{theorem}
Theorem~\ref{thm:jsq:main} is proved using the same machinery as in the $G/G/1$ example, so we defer its proof to Appendix~\ref{app:jsq}. The only term  we left unbounded  in Theorem~\ref{thm:jsq:main} is  the boundary term 
\begin{align*}
    \frac{1}{n} \sum_{i=1}^{n} \E \Big(1(Q_i=0)   \Big(\sum_{j=1}^{n} Q_j + \lambda R_{a} + \sum_{j=1}^{n} \mu R_{s,j}\Big)   \Big).
\end{align*}
Based on the difficulty of bounding $\E R_{a} 1(X=0)$ in the $G/G/1$ workload model, we anticipate that bounding $\E 1(Q_i=0) R_{a}$ and $\E 1(Q_i=0) R_{s,j}$ (for $j \neq i$) is also  challenging, though we can get crude (sub-$(1-\rho)$) bounds using H\"older's inequality as in \eqref{eq:crude:workload}. The 
\begin{align*}
    \E 1(Q_i=0)   \big(\sum_{j=1}^{n} Q_j\big)
\end{align*}
term can be interpreted as a state-space collapse (SSC) term that we expect to be of order $(1-\rho)$; see \cite[pg.  375]{HurtMagu2022} for a treatment of this term in the discrete-time model. Since the JSQ routing scheme seeks to equalize queue lengths, we would expect that when $Q_i = 0$ (by Lemma~\ref{lem:relationships:JSQ} we know that this happens with probability $(1-\rho)$), all other queue lengths are similar to $Q_i$. For recent work on the JSQ model with general interarrival and service times, see \cite{DaiGuanXu2024}.

\section{Future directions}
\label{sec:conclusion} 
Extending the generator approach of Stein's method to PDMPs allowed us to analyze diffusion approximation error for queueing systems with generally distributed primitives. We conclude by highlighting several directions for future research.

\subsubsection*{Extending the prelimit approach to multiple   clocks.} 
We modify the $G/G/1$ system from Section~\ref{sec:workload} to have two superimposed independent renewal streams. Denote the interarrival distribution and the residual interarrival time of stream $i = 1,2$, by $U_{i}$ and $R_{a,i}(t)$, respectively and let $\lambda_{i} = 1/\E U_{i}$.  The   state descriptor for the workload process is
\begin{align*}
    Z(t)=(V(t), R_{a,1}(t), R_{a,2}(t)),
\end{align*}
and the  Poisson equation, derived as in \eqref{eq:the:poisson:abuse}, is 
\begin{align*}
     &-  \partial_{v} \E  F_h(v,R_{a,1},R_{a,2}) \\
     &+ \lambda_{1} \E \big(F_h(v+S,U_1,R_{a,2}) - F_h(v,U_1,R_{a,2}) \big)\\
    &+ \lambda_{2} \E \big(F_h(v+S,R_{a,1},U_{2}) - F_h(v,R_{a,1},U_{2}) \big)  
    =  \E h(V) - h(v), \quad v \geq 0.  
\end{align*}
However,  relating the mixed expectations  
\begin{align*}
    \E \big(F_h(v+S,U_1,R_{a,2}) - F_h(v,U_1,R_{a,2}) \big) \quad \text{ and } \quad F_h(v+S,R_{a,1},U_{2}) - F_h(v,R_{a,1},U_{2})
\end{align*}
to the “all-residuals” quantity 
\begin{align*}
    \E (F_h(v+S,R_{a,1},R_{a,2}) - F_h(v,R_{a,1},R_{a,2}) )
\end{align*}
is substantially more delicate in the presence of multiple general clocks. 

The single-clock argument
underlying Lemma~\ref{lem:relate}  ---  which conditions on a fixed residual time $r_a$ and “waits” until
the next arrival  ---  no longer cleanly applies: the path segment from time $0$ to the next type-$1$ arrival has a complicated dependence
on the entire trajectory of the type-$2$ clock. In particular, type-$2$ arrivals may occur before
$r_{a,1}$, repeatedly resetting $R_{a,2}(t)$ and altering the workload trajectory, which introduces
nontrivial cross-terms that we do not know how to handle.

 In
models with multiple general clocks, the cross-dependencies  preclude a
straightforward extension of Lemma~\ref{lem:relate} and appear to require genuinely new ideas. We therefore leave this extension as an open problem.

\subsubsection*{Stein factor bounds for semimartingale reflected Brownian motion.}  
Semimartingale reflected Brownian motion (SRBM) is a fundamental stochastic process that arises as the heavy-traffic limit of many important queueing models, including generalized Jackson networks and certain multiclass queueing networks. Given the prominence of SRBMs in queueing theory, Stein factor bounds for their stationary distributions would be highly valuable, as they would provide a way to quantify the accuracy of SRBM approximations. At present, however, the only such bounds we are aware of are those in Lemma~\ref{lem:exp:stein:factors} for the one-dimensional RBM, whose stationary distribution is exponential.

In general, the geometry of SRBMs can be quite complex, so a natural starting point is a relatively simple example. One candidate is the SRBM corresponding to the tandem queue in Appendix~\ref{app:tandem}; the associated Poisson equation is displayed at the end of that appendix. The objective would be to bound the partial derivatives of its solution up to third order. One possible approach is to use coupling arguments, in the spirit of the synchronous coupling developed in Section~\ref{sec:prelimit}.



\bibliographystyle{informs2014} 
\bibliography{dai20190911} 

\def\cprime{$'$} \def\cprime{$'$} \def\cprime{$'$} \def\cprime{$'$}
  \def\cprime{$'$} \def\cprime{$'$} \def\cprime{$'$}
\begin{thebibliography}{64}
\providecommand{\natexlab}[1]{#1}
\providecommand{\url}[1]{\texttt{#1}}
\providecommand{\urlprefix}{URL }

\bibitem[{Asmussen(2003)}]{Asmu2003}
Asmussen S (2003) \emph{Applied probability and queues}, volume~51 of
  \emph{Applications of Mathematics (New York)} (New York: Springer-Verlag),
  second edition, ISBN 0-387-00211-1, {S}tochastic Modelling and Applied
  Probability.

\bibitem[{Baccelli \protect\BIBand{} Br{\'e}maud(2003)}]{BaccBrem2003}
Baccelli F, Br{\'e}maud P (2003) \emph{Elements of queueing theory: Palm
  martingale calculus and stochastic recurrences}, volume~26 of
  \emph{Applications of Mathematics} (Berlin: Springer), 2nd edition.

\bibitem[{Barbour(1990)}]{Barb1990}
Barbour A (1990) Stein's method for diffusion approximations. \emph{Probab.
  Theory and Related Fields} 84(3):297--322, ISSN 0178-8051,
  \urlprefix\url{http://dx.doi.org/10.1007/BF01197887}.

\bibitem[{Barbour et~al.(2023)Barbour, Ross, \protect\BIBand{}
  Zheng}]{Barbetal2023}
Barbour A, Ross N, Zheng G (2023) Stein’s method, gaussian processes and palm
  measures, with applications to queueing. \emph{The Annals of Applied
  Probability} 33(5):3835--3871.

\bibitem[{Barbour(1988)}]{Barb1988}
Barbour AD (1988) Stein's method and {P}oisson process convergence.
  \emph{Journal of Appl. Probab.} 25:175--184, ISSN 00219002,
  \urlprefix\url{http://www.jstor.org/stable/3214155}.

\bibitem[{Bassamboo \protect\BIBand{} Randhawa(2010)}]{BasaRand2010}
Bassamboo A, Randhawa RS (2010) On the accuracy of fluid models for capacity
  sizing in queueing systems with impatient customers. \emph{Operations
  Research} 58(5):1398--1413,
  \urlprefix\url{http://dx.doi.org/10.1287/opre.1100.0815}.

\bibitem[{Bertsimas \protect\BIBand{} Gamarnik(2022)}]{BertGama2022}
Bertsimas D, Gamarnik D (2022) \emph{Queueing theory: classical and modern
  methods} (Belmont, MA: Dynamic Ideas LLC), 1st edition.

\bibitem[{Besan{\c{c}}on et~al.(2020)Besan{\c{c}}on, Decreusefond,
  \protect\BIBand{} Moyal}]{Besaetal2020}
Besan{\c{c}}on E, Decreusefond L, Moyal P (2020) Stein’s method for diffusive
  limits of queueing processes. \emph{Queueing Systems} 95:173--201.

\bibitem[{Blanchet \protect\BIBand{} Glynn(2006)}]{BlanGlyn2006}
Blanchet J, Glynn P (2006) Complete corrected diffusion approximations for the
  maximum of a random walk. \emph{Annals of Applied Probability}
  16(2):951--983.

\bibitem[{Boon et~al.(2023)Boon, Janssen, \protect\BIBand{} van
  Leeuwaarden}]{Boonetal2023}
Boon M, Janssen A, van Leeuwaarden J (2023) Heavy-traffic single-server queues
  and the transform method. \emph{Indagationes Mathematicae} 34(5):1014--1037.

\bibitem[{Bramson(2011)}]{Bram2011a}
Bramson M (2011) Stability of join the shortest queue networks. \emph{Ann.
  Appl. Probab.} 21(4):1568--1625,
  \urlprefix\url{http://dx.doi.org/10.1214/10-AAP726}.

\bibitem[{Braverman(2022)}]{Brav2022}
Braverman A (2022) The prelimit generator comparison approach of stein’s
  method. \emph{Stochastic Systems} 12(2):181--204,
  \urlprefix\url{http://dx.doi.org/10.1287/stsy.2021.0085}.

\bibitem[{Braverman et~al.(2017)Braverman, Dai, \protect\BIBand{}
  Miyazawa}]{BravDaiMiya2017}
Braverman A, Dai J, Miyazawa M (2017) Heavy traffic approximation for the
  stationary distribution of a generalized {Jackson} network: the {BAR}
  approach. \emph{Stochastic Systems} 7(1):143--196,
  \urlprefix\url{http://projecteuclid.org/euclid.ssy/1495785619}.

\bibitem[{Braverman et~al.(2024)Braverman, Dai, \protect\BIBand{}
  Miyazawa}]{BravDaiMiya2024}
Braverman A, Dai J, Miyazawa M (2024) The {BAR} approach for multiclass
  queueing networks with {SBP} service policies. \emph{Stochastic Systems} .

\bibitem[{Braverman \protect\BIBand{} Dai(2017)}]{BravDai2017}
Braverman A, Dai JG (2017) Stein's method for steady-state diffusion
  approximations of ${M}/\mathit{Ph}/n+{M}$ systems. \emph{Ann. of Appl.
  Probab.} 27(1):550--581, ISSN 1050-5164,
  \urlprefix\url{http://dx.doi.org/10.1214/16-AAP1211}.

\bibitem[{{Braverman} et~al.(2016){Braverman}, {Dai}, \protect\BIBand{}
  {Feng}}]{BravDaiFeng2016}
{Braverman} A, {Dai} JG, {Feng} J (2016) {Stein's method for steady-state
  diffusion approximations: An introduction through the Erlang-A and Erlang-C
  models}. \emph{Stoch. Syst.} 6:301--366,
  \urlprefix\url{http://www.i-journals.org/ssy/viewarticle.php?id=212&layout=abstract}.

\bibitem[{Brown \protect\BIBand{} Xia(2001)}]{BrowXia2001}
Brown TC, Xia A (2001) Stein's method and birth-death processes. \emph{Ann.
  Probab.} 29(3):1373--1403,
  \urlprefix\url{http://dx.doi.org/10.1214/aop/1015345606}.

\bibitem[{Chen et~al.(2011)Chen, Goldstein, \protect\BIBand{}
  Shao}]{ChenGoldShao2011}
Chen LHY, Goldstein L, Shao QM (2011) \emph{Normal approximation by {S}tein's
  method}. Probability and its Applications (New York) (Springer, Heidelberg),
  ISBN 978-3-642-15006-7,
  \urlprefix\url{http://dx.doi.org/10.1007/978-3-642-15007-4}.

\bibitem[{Chen \protect\BIBand{} Whitt(2020)}]{ChenWhit2020}
Chen Y, Whitt W (2020) Extremal models for the gi/ gi/ k waiting-time
  tail-probability decay rate. \emph{Operations Research Letters}
  48(6):770--776.

\bibitem[{Chen \protect\BIBand{} Whitt(2021)}]{ChenWhit2021}
Chen Y, Whitt W (2021) Extremal gi/gi/1 queues given two moments: exploiting
  tchebycheff systems. \emph{Queueing Systems} 97(1):101--124.

\bibitem[{Chen \protect\BIBand{} Whitt(2022{\natexlab{a}})}]{ChenWhit2022}
Chen Y, Whitt W (2022{\natexlab{a}}) Correction to: Extremal gi/gi/1 queues
  given two moments: exploiting tchebycheff systems. \emph{Queueing Systems}
  102(3):553--556.

\bibitem[{Chen \protect\BIBand{} Whitt(2022{\natexlab{b}})}]{ChenWhit2022a}
Chen Y, Whitt W (2022{\natexlab{b}}) Set-valued performance approximations for
  the queue given partial information. \emph{Probability in the Engineering and
  Informational Sciences} 36(2):378--400.

\bibitem[{Dai et~al.(2025)Dai, Glynn, \protect\BIBand{} Xu}]{DaiGlynXu2025}
Dai JG, Glynn P, Xu Y (2025) Asymptotic product-form steady-state for
  generalized jackson networks in multi-scale heavy traffic.
  \urlprefix\url{https://arxiv.org/abs/2304.01499}.

\bibitem[{Dai et~al.(2024)Dai, Guang, \protect\BIBand{} Xu}]{DaiGuanXu2024}
Dai JG, Guang J, Xu Y (2024) Steady-state convergence of the continuous-time
  jsq system with general distributions in heavy traffic. \emph{SIGMETRICS
  Perform. Eval. Rev.} 52(2):39–41, ISSN 0163-5999,
  \urlprefix\url{http://dx.doi.org/10.1145/3695411.3695426}.

\bibitem[{Dai et~al.(2004)Dai, Hasenbein, \protect\BIBand{}
  Vande~Vate}]{DaiHaseVand2004}
Dai JG, Hasenbein JJ, Vande~Vate JH (2004) Stability and instability of a
  two-station queueing network. \emph{Annals of Applied Probability}
  14:326--377, ISSN 1050-5164.

\bibitem[{Dai \protect\BIBand{} Meyn(1995)}]{DaiMeyn1995}
Dai JG, Meyn SP (1995) Stability and convergence of moments for multiclass
  queueing networks via fluid limit models. \emph{IEEE Transactions on
  Automatic Control} 40:1889--1904.

\bibitem[{Dai \protect\BIBand{} Xu(2024)}]{DaiXu2024}
Dai JG, Xu Y (2024) Explicit steady-state approximations for parallel server
  systems with heterogeneous servers.
  \urlprefix\url{https://arxiv.org/abs/2406.04203}.

\bibitem[{Daley et~al.(1992)Daley, Kreinin, \protect\BIBand{}
  Trengove}]{Daleetal1992}
Daley D, Kreinin AY, Trengove C (1992) Inequalities concerning the waiting time
  in single-server queues: A survey. \emph{Oxford Statistical Science Series}
  1(9):177--177.

\bibitem[{Davis(1984)}]{Davi1984}
Davis MHA (1984) Piecewise deterministic {M}arkov processes: a general class of
  non-diffusion stochastic models. \emph{Journal of Royal Statist.\ Soc.\,
  series B} 46:353--388.

\bibitem[{Eryilmaz \protect\BIBand{} Srikant(2012)}]{EryiSrik2012}
Eryilmaz A, Srikant R (2012) Asymptotically tight steady-state queue length
  bounds implied by drift conditions. \emph{Queueing Systems} 72(3-4):311--359,
  ISSN 0257-0130, \urlprefix\url{http://dx.doi.org/10.1007/s11134-012-9305-y}.

\bibitem[{Feng \protect\BIBand{} Shi(2018)}]{FengShi2018}
Feng J, Shi P (2018) Steady-state diffusion approximations for discrete-time
  queue in hospital inpatient flow management. \emph{Naval Research Logistics
  (NRL)} 65(1):26--65, \urlprefix\url{http://dx.doi.org/10.1002/nav.21787}.

\bibitem[{Gast(2017)}]{Gast2017}
Gast N (2017) Expected values estimated via mean-field approximation are
  1/n-accurate. \emph{Proceedings of the ACM on Measurement and Analysis of
  Computing Systems} 1(1):1--26.

\bibitem[{Gaunt \protect\BIBand{} Walton(2020)}]{GaunWalt2020}
Gaunt RE, Walton N (2020) Stein’s method for the single server queue in heavy
  traffic. \emph{Statistics \& Probability Letters} 156:108566, ISSN 0167-7152,
  \urlprefix\url{http://dx.doi.org/https://doi.org/10.1016/j.spl.2019.108566}.

\bibitem[{G{\"o}tze(1991)}]{Gotz1991}
G{\"o}tze F (1991) On the rate of convergence in the multivariate {CLT}.
  \emph{Ann. Probab.} 19(2):724--739,
  \urlprefix\url{http://dx.doi.org/10.1214/aop/1176990448}.

\bibitem[{Gurvich(2014)}]{Gurv2014}
Gurvich I (2014) Diffusion models and steady-state approximations for
  exponentially ergodic {M}arkovian queues. \emph{Ann. Appl. Probab.}
  24(6):2527--2559, \urlprefix\url{http://dx.doi.org/10.1214/13-AAP984}.

\bibitem[{Gurvich et~al.(2014)Gurvich, Huang, \protect\BIBand{}
  Mandelbaum}]{GurvHuanMand2014}
Gurvich I, Huang J, Mandelbaum A (2014) Excursion-based universal
  approximations for the {Erlang-A} queue in steady-state. \emph{Math. Oper.
  Res.} 39(2):325--373,
  \urlprefix\url{http://dx.doi.org/10.1287/moor.2013.0606}.

\bibitem[{Gut(1974)}]{Gut1974}
Gut A (1974) On the moments and limit distributions of some first passage
  times. \emph{The Annals of Probability} 277--308.

\bibitem[{Harrison \protect\BIBand{} Reiman(1981)}]{HarrReim1981}
Harrison JM, Reiman MI (1981) Reflected {B}rownian motion on an orthant.
  \emph{Ann. Probab.} 9(2):302--308, ISSN 0091-1798,
  \urlprefix\url{http://links.jstor.org/sici?sici=0091-1798(198104)9:2<302:RBMOAO>2.0.CO;2-P&origin=MSN}.

\bibitem[{Huang \protect\BIBand{} Gurvich(2018)}]{GurvHuan2018}
Huang J, Gurvich I (2018) Beyond heavy-traffic regimes: Universal bounds and
  controls for the single-server queue. \emph{Oper. Res.} 66(4):1168--1188,
  \urlprefix\url{http://dx.doi.org/10.1287/opre.2017.1715}.

\bibitem[{Hurtado-Lange \protect\BIBand{} Maguluri(2022)}]{HurtMagu2022}
Hurtado-Lange D, Maguluri ST (2022) A load balancing system in the many-server
  heavy-traffic asymptotics. \emph{Queueing Systems} 101(3):353--391.

\bibitem[{Hurtado~Lange \protect\BIBand{} Maguluri(2022)}]{HurtMagu2022a}
Hurtado~Lange DA, Maguluri ST (2022) Heavy-traffic analysis of queueing systems
  with no complete resource pooling. \emph{Mathematics of Operations Research}
  47(4):3129--3155, \urlprefix\url{http://dx.doi.org/10.1287/moor.2021.1248}.

\bibitem[{Kingman(1961)}]{King1961a}
Kingman JFC (1961) The single server queue in heavy traffic. \emph{Mathematical
  Proceedings of the Cambridge Philosophical Society} 57:902--904,
  \urlprefix\url{http://dx.doi.org/10.1017/S0305004100036094}.

\bibitem[{Kingman(1962)}]{King1962}
Kingman JFC (1962) On queues in heavy traffic. \emph{J. Roy. Statist. Soc. Ser.
  B} 24:383--392, ISSN 0035-9246,
  \urlprefix\url{http://links.jstor.org/sici?sici=0035-9246(1962)24:2<383:OQIHT>2.0.CO;2-A&origin=MSN}.

\bibitem[{Köllerström(1976)}]{Koll1976}
Köllerström J (1976) Stochastic bounds for the single-server queue.
  \emph{Mathematical Proceedings of the Cambridge Philosophical Society}
  80(3):521–525, \urlprefix\url{http://dx.doi.org/10.1017/S0305004100053135}.

\bibitem[{Li \protect\BIBand{} Ou(1995)}]{LiOu1995}
Li J, Ou J (1995) Characterizing the idle-period distribution of gi/g/1 queues.
  \emph{Journal of applied probability} 32(1):247--255.

\bibitem[{Lieberman(2013)}]{Lieb2013}
Lieberman GM (2013) \emph{Oblique Derivative Problems for Elliptic Equations}
  (WORLD SCIENTIFIC), \urlprefix\url{http://dx.doi.org/10.1142/8679}.

\bibitem[{Lorden(1970)}]{Lord1970}
Lorden G (1970) {On excess over the boundary}. \emph{The Annals of Mathematical
  Statistics} 41(2):520 -- 527,
  \urlprefix\url{http://dx.doi.org/10.1214/aoms/1177697092}.

\bibitem[{Loulou(1978)}]{Loul1978}
Loulou R (1978) An explicit upper bound for the mean busy period in a
  {$GI/G/1$} queue. \emph{Journal of Applied Probability} 15(2):452--455,
  \urlprefix\url{http://dx.doi.org/10.2307/3213419}.

\bibitem[{Mackey \protect\BIBand{} Gorham(2016)}]{GorhMack2016}
Mackey L, Gorham J (2016) Multivariate {S}tein factors for a class of strongly
  log-concave distributions. \emph{Electron. Commun. Probab.} 21:14,
  \urlprefix\url{http://dx.doi.org/10.1214/16-ECP15}.

\bibitem[{Maguluri et~al.(2018)Maguluri, Burle, \protect\BIBand{}
  Srikant}]{MaguBurleSrik2018}
Maguluri ST, Burle SK, Srikant R (2018) Optimal heavy-traffic queue length
  scaling in an incompletely saturated switch. \emph{Queueing Systems}
  88(3):279--309.

\bibitem[{Maguluri \protect\BIBand{} Srikant(2016)}]{MaguSrik2016}
Maguluri ST, Srikant R (2016) Heavy traffic queue length behavior in a switch
  under the maxweight algorithm. \emph{Stochastic Systems} 6(1):211--250,
  \urlprefix\url{http://dx.doi.org/10.1214/15-SSY193}.

\bibitem[{Miyazawa(2015)}]{Miya2015}
Miyazawa M (2015) Diffusion approximation for stationary analysis of queues and
  their networks: a review. \emph{J. Oper. Res. Soc. Japan} 58(1):104--148,
  ISSN 0453-4514, \urlprefix\url{http://dx.doi.org/10.15807/jorsj.58.104}.

\bibitem[{Miyazawa(2017)}]{Miya2017}
Miyazawa M (2017) A unified approach for large queue asymptotics in a
  heterogeneous multiserver queue. \emph{Adv. in Appl. Probab.} 49(1):182--220,
  ISSN 0001-8678, \urlprefix\url{http://dx.doi.org/10.1017/apr.2016.84}.

\bibitem[{Peköz \protect\BIBand{} Röllin(2011)}]{PekoRoll2011}
Peköz EA, Röllin A (2011) {New rates for exponential approximation and the
  theorems of {R}ényi and {Y}aglom}. \emph{The Annals of Probability}
  39(2):587 -- 608, \urlprefix\url{http://dx.doi.org/10.1214/10-AOP559}.

\bibitem[{Ross(2011)}]{Ross2011}
Ross N (2011) Fundamentals of {S}tein's method. \emph{Probab. Surv.}
  8:210--293, ISSN 1549-5787,
  \urlprefix\url{http://dx.doi.org/10.1214/11-PS182}.

\bibitem[{Siegmund(1979)}]{Sieg1979}
Siegmund D (1979) Corrected diffusion approximations in certain random walk
  problems. \emph{Advances in Applied Probability} 11(4):701--719, ISSN
  00018678, \urlprefix\url{http://www.jstor.org/stable/1426855}.

\bibitem[{Stein(1972)}]{Stei1972}
Stein C (1972) A bound for the error in the normal approximation to the
  distribution of a sum of dependent random variables. \emph{Proceedings of the
  Sixth Berkeley Symposium on Mathematical Statistics and Probability, Volume
  2: Probability Theory}, 583--602 (Berkeley, Calif.: University of California
  Press), \urlprefix\url{http://projecteuclid.org/euclid.bsmsp/1200514239}.

\bibitem[{Stolyar(2015)}]{Stol2015}
Stolyar AL (2015) Tightness of stationary distributions of a flexible-server
  system in the {H}alfin-{W}hitt asymptotic regime. \emph{Stoch. Syst.}
  5(2):239--267, \urlprefix\url{http://dx.doi.org/10.1214/14-SSY139}.

\bibitem[{Ward \protect\BIBand{} Glynn(2003)}]{GlynWard2003}
Ward A, Glynn P (2003) A diffusion approximation for a markovian queue with
  reneging. \emph{Queueing Systems} 43(1-2):103--128, ISSN 0257-0130,
  \urlprefix\url{http://dx.doi.org/10.1023/A:1021804515162}.

\bibitem[{Whitt(1986)}]{Whit1986a}
Whitt W (1986) Deciding which queue to join: Some counterexamples. \emph{Oper.
  Res.} 34(1):55–62, ISSN 0030-364X.

\bibitem[{Wolff \protect\BIBand{} Wang(2003)}]{WolfWang2003}
Wolff RW, Wang CL (2003) Idle period approximations and bounds for the gi/g/1
  queue. \emph{Advances in Applied Probability} 35(3):773--792.

\bibitem[{Ying(2016)}]{ying2016}
Ying L (2016) On the approximation error of mean-field models.
  \emph{Proceedings of the 2016 ACM SIGMETRICS International Conference on
  Measurement and Modeling of Computer Science}, 285--297 (Antibes
  Juan-les-Pins, France: ACM),
  \urlprefix\url{http://dx.doi.org/10.1145/2964791.2901463}.

\bibitem[{Ying(2017)}]{Ying2017}
Ying L (2017) Stein's method for mean field approximations in light and heavy
  traffic regimes. \emph{Proc. ACM Meas. Anal. Comput. Syst.} 1(1):12:1--12:27,
  ISSN 2476-1249, \urlprefix\url{http://dx.doi.org/10.1145/3084449}.

\bibitem[{Zhou \protect\BIBand{} Shroff(2020)}]{ZhouShro2020}
Zhou X, Shroff N (2020) A note on stein's method for heavy-traffic analysis.
  \emph{arXiv preprint arXiv:2003.06454} .

\end{thebibliography}



%
%

\begin{APPENDICES} 

\section{The $G/M/\infty$ queue length --- a boundary-free example}
\label{sec:gminfinity}
The $G/M/\infty$ system provides an instructive example that features both general and exponential clocks, and in contrast to the $G/G/1$ setting, the diffusion error for the infinite-server queue is composed solely of interior terms and no  boundary terms. As a result, obtaining a three-moment error bound  is considerably simpler.

The $G/M/\infty$ system is an infinite-server queue where arrivals form a renewal process with generic interarrival time distribution $U$ and exponentially distributed service times with rate $\mu$. Let $R = \lambda/\mu$ be the offered load. We assume that $\E U^3 < \infty$ and set 
\begin{align*}
    &\lambda = 1/\E U, \quad c_{U}^2 = \lambda^2 \text{Var}(U).
\end{align*}
Let $Q(t)$ be the customer count, let $R_{a}(t)$ be the residual interarrival time at $t \geq 0$, and assume that
\begin{align*}
    \{Z(t) = (Q(t),R_{a}(t))\}
\end{align*}
is positive Harris recurrent. Let $Z=(Q,R_a)$ denote the vector having the stationary distribution. Let $\delta = 1/\sqrt{R}$ and  define the normalized customer count 
\begin{align*}
    X(t) = \delta(Q(t)-R) \quad \text{ and  } \quad X = \delta(Q-R).
\end{align*}
The following is the main result of this section. 
\begin{theorem}
    \label{thm:gminf:main}
    Let $Y \sim N(0,\sigma^2)$ with $\sigma^2 = (1+  c_{U}^{2})$. Then  
    \begin{align*}
        &d_{W}(X,Y) \leq  \delta(K_A + K_D), \\
        &K_A =   \frac{1}{3(1+  c_{U}^{2})}  \E \abs{1-\lambda U}^{3} +  \frac{c_{U}^{2}}{1+  c_{U}^{2}}  \Big(1 + \frac{1}{2} \lambda^2\E U^2  ( \delta^4\lambda^3 \E U^3/3 +  \delta^2\lambda^2 \E U^2  + 1  ) \Big), \\
        &K_D  = \sqrt{\frac{2/\pi}{(1+  c_{U}^{2})}} \sqrt{ \E X^2} \sqrt{\lambda^3\E U^3/3} +\frac{1}{(1+  c_{U}^{2})} \Big( \delta^2  \lambda^3 \E U^3/3 +   \lambda^2 \E U^2 + \frac{1}{3} \Big).
    \end{align*}
\end{theorem} 

The proof of Theorem~\ref{thm:gminf:main} follows the same structure as the proof of Theorem~\ref{thm:workload:1}. Namely, we use the BAR for the compensated customer count in the infinite-server queue to extract a diffusion generator and corresponding error terms in Proposition~\ref{prop:extraction:gminf}. Bounding the error terms requires both Stein factor bounds and a bound  on the second moment of the centered queue length (unlike our $G/G/1$ example, which did not require moment bounds). The former are well known   and the latter is obtained  from the BAR by using a quadratic test function. After presenting said Stein factor and moment bounds,   we combine everything to prove Theorem~\ref{thm:gminf:main} at the end of the section. 

For the BAR, let $A(t)$ and $U(t)$ be defined as in Section~\ref{sec:model:workload}, and let $D(t)$ denote the number of departures on $[0,t]$. Given a function $f(Z(t))$, partition $[0,t]$ according to event times and consider the FTC conditions analogous to those in Section~\ref{sec:firstBAR}, as well as the integrability condition 
\begin{align*}
    \E \abs{f(Z)}, \  \E \abs{\partial_{r_{a}} f(Z)},\  \E \int_{0}^{t} \abs{\Delta  f(Z(s-))} dD(s),\  \E \int_{0}^{t} \abs{\Delta  f(Z(s-))} dA(s) < \infty. 
\end{align*} 
Arguing as in Section~\ref{sec:firstBAR} yields the BAR.
\begin{lemma}
\label{lem:full:BAR:gminf}
Initialize $Z(0) \sim Z$. If $f(Z(s))$ satisfies the FTC conditions with probability one under $Z(0) \sim Z$ and if the integrability condition  holds, then   for all $t \geq 0$,
\begin{align}
    0 =&\ - t \E (\partial_{r_a} f(Z)) +\E  \int_{0}^{t} \Delta f(Z(s-)) d D(s)  +  \E  \int_{0}^{t} \Delta f(Z(s-)) d A(s)  \notag \\
    =&\ - t \E (\partial_{r_a} f(Z)) + t \mu \E  Q(f(Q-1,R_a) -f(Q,R_a))  +  \E  \int_{0}^{t} \Delta f(Z(s-)) d A(s),\label{eq:full:BAR:gminf}
\end{align}
provided that all expectations on the right-hand side are well defined (in a manner similar to \eqref{eq:ev}).
\end{lemma}
\startproof{Proof of Lemma~\ref{lem:full:BAR:gminf}}
The first equality follows as in Lemma~\ref{lem:full:BAR:workload} and the second is due to Theorem~3.3.1 of \cite{BaccBrem2003}. 
\finishproof
\begin{remark}
\label{rem:tsmall}
For convenience, we will work with the right-hand side of \eqref{eq:full:BAR:gminf} with $t = 1$ there. For an approach that avoids using the machinery of \cite{BaccBrem2003} to establish the second equality in \eqref{eq:full:BAR:gminf}, consider the first equality in \eqref{eq:full:BAR:gminf}, divide both sides by $t$, and let $t \to 0$ to get 
    \begin{align*}
        0 =&\ - \E (\partial_{r_a} f(Z)) + \lim_{t \to 0}\frac{1}{t} \E  \int_{0}^{t} \Delta f(Z(s-)) d D(s)  +  \lim_{t \to 0}\frac{1}{t}  \E  \int_{0}^{t} \Delta f(Z(s-)) d A(s) \\
        =&\ - \E (\partial_{r_a} f(Z)) + \mu \E  Q(f(Q-1,R_a) -f(Q,R_a))  +  \lim_{t \to 0}\frac{1}{t}  \E  \int_{0}^{t} \Delta f(Z(s-)) d A(s).
    \end{align*}
    The second equality is argued using the fact that $\Prob(D(t) = 1) = t \mu Q(0) + o(t)$ and $\Prob(D(t) > 1) = o(t)$ for $t$. Since $\E A(t) = \lambda t$, the presence of the limit as $t \to 0$ in the third term does not add further complexity (compared to the $t = 1$ case) to our diffusion generator extraction procedure (Proposition~\ref{prop:extraction:gminf}). 
\end{remark}
Defining the compensated customer count 
\begin{align}
    \widetilde{X}(t) = X(t) - \delta \lambda R_{a}(t) \quad \text{ and } \quad \widetilde{X} = X - \delta \lambda R_{a},
\end{align}
and specializing the BAR \eqref{eq:full:BAR:gminf} to  $\TX$ yields
\begin{align}
     \delta  \lambda  \E  f'(\TX) + \mu \E \big( Q \big( f(\TX -\delta)-f(\TX)\big)\big) +  \E  \int_{0}^{1} \Delta f(\TX(t-)) d A(t)    = 0. \label{eq:BAR3:mginf}
\end{align}
The following expansion of \eqref{eq:BAR3:mginf} is analogous to Proposition~\ref{prop:extraction:workload} and is proved in  Appendix~\ref{app:gminf}.
\begin{proposition}
\label{prop:extraction:gminf}
If  $f\in C^{2}(\R)$ with $f''(x)$ absolutely continuous and $\norm{f''},\norm{f'''} < \infty$, then, provided that all expectations are well defined,
\begin{align}
        \mu \E \big( Q \big( f(\TX -\delta)-f(\TX)\big)\big)  =&\ -\delta R f'(\TX)  -\mu \E X f'(X) + \frac{1}{2} \mu \E f''(X) + \epsilon_{D}(f)  , \label{eq:eD:gminf} \\
        \E \int_{0}^{1} \Delta  f(\TX(t-))  d A(t) =&\ \frac{1}{2}\mu c_{U}^{2} \E f''(X) + \epsilon_{A}(f)  \label{eq:eA:gminf} 
    \end{align}
    where 
\begin{align*}
    \abs{\epsilon_{A}(f)} \leq&\      \frac{1}{6} \delta \mu \norm{f'''} \E \abs{1-\lambda U}^{3}  + \frac{1}{2} \delta \mu c_{U}^{2}\norm{f'''}  \Big(1 + \frac{1}{2} \lambda^2\E U^2  ( \E (Q-R)^2/R^2 + 1  ) \Big),\\
    \abs{\epsilon_{D}(f)} \leq&\ \mu \delta \norm{f''} \sqrt{ \E X^2} \sqrt{\lambda^3\E U^3/3} +  \frac{1}{2}   \mu \delta \norm{f'''} \delta^2  \lambda \E  (Q    R_{a})   + \frac{1}{6}  \mu \delta \norm{f'''}.
\end{align*}
 \end{proposition} 

The upper bounds on $\abs{\epsilon_{A}(f_h)}$ and $\abs{\epsilon_{D}(f_h)}$ contain terms involving $\E(Q-R)^2$ and $\E  (Q    R_{a})$. The following lemma  presents some useful identities, including a bound on these quantities. It is  proved in Appendix~\ref{sec:relationships:proof:gminf}.
\begin{lemma}
    \label{lem:relationships:gminf}
    Recall that $R = \lambda/\mu$. For any $m > 1$,
    \begin{align}
    & \E A(1) = \lambda,  \quad  \E R_{a}^{m-1} =  \lambda \E U^{m}/m, \quad \E Q = R, \label{eq:standard:gminf} \\ 
    & \E (Q-R)^2 = \lambda \E Q R_{a} \leq \lambda^3 \E U^3/3 +  R \lambda^2 \E U^2,  \label{eq:mombound:gminf}\\
    & \E \int_{0}^{1}Q(t-) d A(t) = \mu \E (Q-R)^2 + \lambda R -\lambda. \label{eq:qt-:gminf}
    \end{align}
\end{lemma}
 
Applying the expansions in \eqref{eq:eD:gminf}--\eqref{eq:eA:gminf} to the BAR \eqref{eq:BAR3:mginf} suggests the diffusion approximation with  generator 
\begin{align*}
  G_{Y} f(x) =  -\mu x f'(x) + \frac{1}{2} \mu (1+  c_{U}^{2}) f''(x), \quad x \in \R,
\end{align*}
which corresponds to the $N(0,\sigma^2)$ distribution with $\sigma^2 = (1+  c_{U}^{2})$. The following lemma  is  a rescaled version  of \citep[Lemma 2.4]{ChenGoldShao2011}, and is proved in Appendix~\ref{app:stein:normal}.
\begin{lemma}
    \label{lem:normal:stein:factors}
Let $Y \sim N(0,\sigma^2)$ and let $f_{h,\sigma}(x)$ be  the solution to  
\begin{align}
 - x f_{h,\sigma}'(x) + \frac{1}{2}  \sigma^2 f_{h,\sigma}''(x) = \frac{1}{\mu} (\E h(Y) - h(x)), \quad x \in \R, \label{eq:poisson:equation:normal}  
\end{align}
with $f_{h,\sigma}(0)=0$. Then 
    \begin{align*}
        \norm{f_{h,\sigma}'} \leq \frac{2}{\mu}, \quad \norm{f_{h,\sigma}''} \leq \frac{\sqrt{2/\pi}}{\mu \sigma }, \quad \text{ and } \quad \norm{f_{h,\sigma}'''} \leq \frac{2}{\mu \sigma^2}.
    \end{align*}
    As a consequence, $\abs{f_{h,\sigma}(x)} \leq 2\abs{x}/\mu$ for all $x \in \R$.
\end{lemma}

We are ready to prove Theorem~\ref{thm:gminf:main}.
\startproof{Proof of Theorem~\ref{thm:gminf:main}}
Since $\norm{f_h'}< \infty$ and $\abs{f_h(x)} \leq 2\abs{x}/\mu$ for any $h \in \lipone$, one readily verifies that the BAR \eqref{eq:BAR3:mginf} holds with $f_h(X - \delta \lambda R_{a})$. Applying  the Stein factor bounds of Lemma~\ref{lem:normal:stein:factors} and the moment bound in \eqref{eq:mombound:gminf} of Lemma~\ref{lem:relationships:gminf} to the upper bounds on $\abs{\epsilon_{A}(f_h)}$ and $\abs{\epsilon_{D}(f_h)}$ in Proposition~\ref{prop:extraction:gminf} yields $\abs{\epsilon_{A}(f_h)} \leq K_A$ and $\abs{\epsilon_{D}(f_h)} \leq K_D$. 
 \finishproof
 
 \rr{The following appendix treats the two-station tandem queue, the only example whose diffusion approximation is multidimensional. }
\section{The tandem queue --- a multidimensional example}
\label{app:tandem}
In higher dimensions, the mechanics are unchanged: we derive the diffusion generator from the BAR and
decompose the approximation error into interior and boundary contributions. Our approach also extends
naturally to generalized Jackson networks. The added difficulty
lies in establishing Stein factor bounds for the multidimensional RBMs;
bounds are not known even for the two-dimensional RBM approximating the tandem queue. 

Consider two  first-come-first-served single-server stations in tandem. Customers arrive to station one according to a  renewal
process with generic interarrival time $U$; after service at station one they move to station two and depart after completing service at station two. Service times at station $i$ are i.i.d.\ $S_i$, independent of
the arrival process and of service times at the other station. Assume $\E U^3<\infty$, $\E S_1^3,\E S_2^3<\infty$,
and no simultaneous events almost surely. 
Set
\begin{align*}
    &\lambda = 1/\E U,\quad \mu_i = 1 / \E S_{i},\quad \rho_i = \lambda/\mu_{i}, \\
    &c_{U}^{2} = \lambda^2\text{Var}(U), \quad c_{S,i}^{2} = \mu_{i}^2\text{Var}(S_i).
\end{align*}
For $t\ge0$, let $Q_i(t)$ be the queue length at station $i$, $R_a(t)$ the residual interarrival
time, and $R_{s,i}(t)$ the remaining service at station $i$ (if $Q_i(t)=0$, the service time of the
next customer). Define $X(t)=(\delta_1 Q_1(t),\delta_2 Q_2(t))$ with $\delta_i=1-\rho_i$,
$R_s(t)=(R_{s,1}(t),R_{s,2}(t))$, and $Z(t)=(X(t),R_a(t),R_s(t))$.  Assume $\{Z(t)\}$ is positive
Harris recurrent (cf.\ \cite{DaiMeyn1995}) and let $Z=(X,R_a,R_s)$ have the stationary distribution.
Let $A(t)$ and $D_i(t)$ denote arrivals and station-$i$ departures on $[0,t]$. 

The procedure of applying the generator approach to the tandem queue is similar to the three examples already considered. We therefore  highlight only the key steps, leaving the formal derivations to the interested reader. For sufficiently regular functions, the BAR for $Z$ is 
\begin{align}
    & - \E  (\partial_{r_a} f(Z))  - \sum_{i=1}^{2}\E (1(X_i>0)\partial_{r_{s,i}} f(Z) )   \notag \\
    &+  \E \int_{0}^{1} \Delta f(Z(t-)) d A(t)+ \sum_{i=1}^{2} \E \int_{0}^{1} \Delta f(Z(t-)) d D_{i}(t) = 0.\label{eq:full:BAR:4}
\end{align} 
Let $\TX$ have the stationary distribution of the   compensated queue-length vector $\TX(t)=(\TX_1(t),\TX_2(t))$, where  
\begin{align*}
    \TX_1(t) = X_1(t) - \delta_1\lambda R_{a}(t) + \delta_1\mu_1 R_{s,1}(t), \quad \TX_2(t) = X_2(t) - \delta_2\mu_1 R_{s,1}(t) + \delta_2\mu_2 R_{s,2}(t).
\end{align*} 
The BAR for $\TX$ is 
\begin{align*}
    &   \delta_1 \E \big((\lambda - \mu_1 1(Q_1>0))  \partial_{x_1} f(\TX) \big) + \delta_2 \E \big((  \mu_1 1(Q_1>0) - \mu_2 1(Q_2>0))  \partial_{x_2} f(\TX) \big)    \notag \\
    &+  \E \int_{0}^{1} \Delta f(\TX(t-)) d A(t)+ \sum_{i=1}^{2} \E \int_{0}^{1} \Delta f(\TX(t-)) d D_{i}(t) = 0,
\end{align*}
Omitting the higher-order error terms, it follows that
\begin{align*}
     &\E \int_{0}^{1} \Delta f(\TX(t-)) d A(t) \approx   \frac{1}{2} \delta_1^2  \lambda c_{U}^{2}\E \partial_{x_1}^{2} f(X), \\
    &\E \int_{0}^{1} \Delta f(\TX(t-)) d D_{1}(t)  \approx  \frac{1}{2}   \mu_1 c_{S,1}^{2} \E\big(\delta_1^2  \partial_{x_1}^{2} f(X) - 2\delta_{1} \delta_{2}   \partial_{x_1} \partial_{x_2} f(X) + \delta_{2}^{2}  \partial_{x_2}^{2} f(X) \big), \\
     &\E \int_{0}^{1} \Delta f(\TX(t-)) d D_{2}(t) \approx    \frac{1}{2} \delta_{2}^{2}   \mu_2c_{S,2}^{2}\E \partial_{x_2}^{2} f(X).
\end{align*}
As in our prior examples, the omitted error terms can  be bounded using the first three moments of the primitives as well as the second- and third-order derivatives of $f(x)$. Furthermore, omitting the boundary error terms, we have  
\begin{align*}
    & \delta_1 \E \big((\lambda - \mu_1 1(Q_1>0))  \partial_{x_1} f(\TX) \big) + \delta_2 \E \big((  \mu_1 1(Q_1>0) - \mu_2 1(Q_2>0))  \partial_{x_2} f(\TX) \big) \\
    \approx&\  -\mu_1\delta_1^2  \E \partial_{x_1} f(X) + \delta_2 (\mu_1 \delta_1  - \mu_2 \delta_2 )\E \partial_{x_2} f(X)  \\
    &+ \mu_{1} \big(\delta_{1} 1(Q_1=0)\partial_{x_1} f(X) - \delta_2   1(Q_1=0) \partial_{x_2} f(X) \big) + \mu_2 \delta_2 1(Q_2=0)\partial_{x_2} f(X).
\end{align*} 
The second line on the right-hand side is tied to the reflection structure of the approximating RBM, which we introduce shortly. It is similar to the $f'(0)$ found in \eqref{eq:e0:workload} of Proposition~\ref{prop:extraction:workload}  and is not simply an error term. 

Our expansions suggest the diffusion generator
\begin{align*}
   G_{Y}f(x) =&\ -\mu_1\delta_1^2    \partial_{x_1} f(x) + \delta_2 (\mu_1 \delta_1  - \mu_2 \delta_2 )  \partial_{x_2} f(x)  + \frac{1}{2}\delta_1^2\big(  \lambda c_{U}^{2}  +   \mu_1 c_{S,1}^{2} \big)  \partial_{x_1}^{2} f(x) \\
    &- \delta_1\delta_2  \mu_1 c_{S,1}^{2}  \partial_{x_1} \partial_{x_2} f(x)  + \frac{1}{2} \delta_{2}^{2} \big(   \mu_1 c_{S,1}^{2} +   \mu_2 c_{S,2}^{2}\big) \partial_{x_2}^{2} f(x)  \\
    &+ \mu_1( 1(x_1=0)(\delta_1\partial_{x_1} f(x) - \delta_2 \partial_{x_2} f(x))) + \mu_2 \delta_2 1(x_2=0)\partial_{x_2} f(x),\quad x \in \R^{2}_{+},
\end{align*}
which  corresponds to a two-dimensional RBM $\{Y(t) \in \R^{2}_+ : t \geq 0\}$ on the nonnegative orthant  defined as follows. Set
\begin{align*}
    \mu = \begin{pmatrix}
        \mu_1 \\ \mu_2 
    \end{pmatrix}, \quad 
    \Sigma = 
    \begin{pmatrix}
        \lambda c_{U}^{2}  +   \mu_1 c_{S,1}^{2}  & -  \mu_1 c_{S,1}^{2} \\
        -  \mu_1 c_{S,1}^{2} & \mu_1 c_{S,1}^{2} +   \mu_2 c_{S,2}^{2}
    \end{pmatrix}, \quad 
    R = 
    \begin{pmatrix}
        1 & 0 \\
        -1 & 1
    \end{pmatrix},
\end{align*}
let $\{\xi(t) : t \geq 0 \}$ be a two-dimensional Brownian motion with drift $b = -R \mu$ and covariance matrix $\Sigma$, and let 
\begin{align*}
    Y(t) = \begin{pmatrix}
        \delta_1 & 0 \\
        0 & \delta_2
    \end{pmatrix}  \widetilde{Y}(t),
\end{align*}
where 
\begin{align*}
    \widetilde{Y}(t) =   \xi(t) + R I(t) ,\quad t \geq 0, 
\end{align*}
 and $I: R \to \R^{2}$ is the unique nondecreasing process with $I(0) = 0$ and 
\begin{align*}
    \int_{0}^{\infty} \widetilde{Y}_{i}(t) d I_{i}(t) = 0, \quad i = 1,2.
\end{align*}
For the existence and uniqueness of $\{\widetilde{Y}(t) : t \geq 0\}$  see \cite{HarrReim1981}. 

The final ingredient is the Poisson equation. Fix  $h: \R^2 \to \R$ with $\E \abs{h(Y)} < \infty$, where $Y$ has the stationary distribution of $\{Y(t): t \geq 0\}$, and consider 
\begin{align}
    &G_{Y} f_h(x) =  \E h(Y) - h(x), \quad & x \in \R^{2}_+, \notag \\
    &\delta_1\partial_{x_1} f(x) - \delta_2 \partial_{x_2} f(x) =  0, \quad & x = (0,x_2) \in \R^2_{+}, \label{eq:oblique:1} \\
    &\partial_{x_2} f(x) =  0, \quad & x = (x_1,0) \in \R^{2}_{+}. \label{eq:oblique:2} 
\end{align}
This is known as an oblique derivative problem \citep{Lieb2013} with \eqref{eq:oblique:1} and \eqref{eq:oblique:2} arising from the reflection structure of the RBM, which is driven by the matrix $R$. To bound the diffusion approximation error for the tandem queue, we require bounds on the partial derivatives, up to the third order, of $f_h(x)$. These Stein factor bounds  remain an open problem.

\section{The $G/G/1$ workload: supporting proofs}
\label{app:GG1:workload}

\subsection{Proof of the inversion formula (Lemma~\ref{lem:palm:inversion:2})}
\label{app:workload:inversion}
\startproof{}
Let $\tau_{m}$ be the time of the $m$th arrival   and observe that 
\begin{align*}
    \int_{0}^{1} f(X(t)) dt    =&\  \int_{0}^{\tau_{1}} f(X(t)) dt +  \sum_{m=1}^{A(1)}   \int_{\tau_{m}}^{\tau_{m+1}} f(X(t)) dt     - \int_{1}^{\tau_{A(1)+1}} f(X(t)) dt.
\end{align*}
Since $\tau_{m+1}-\tau_{m} = U(\tau_{m})$, it follows that  
\begin{align*}
    \sum_{m=1}^{A(1)} \int_{\tau_{m}}^{\tau_{m+1}} f(X(t)) dt =&\ \sum_{m=1}^{A(1)} \int_{0}^{U(\tau_{m})} f\big(X(\tau_{m}+u) \big) du \\
    =&\ \int_{0}^{1} \int_{0}^{U(t)} f\big(X(t+u) \big) du d A(t).
\end{align*}
Furthermore, since $\tau_1 = R_{a}(0)$ and $\tau_{A(1)+1} = 1 + R_{a}(1)$, 
\begin{align*}
  \int_{0}^{\tau_{1}} f(X(t)) dt    - \int_{1}^{\tau_{A(1)+1}} f(X(t)) dt =   \int_{0}^{R_{a}(0)} f(X(t)) dt    - \int_{0}^{ R_{a}(1)} f(X(1+t)) dt.
\end{align*}
Therefore, 
\begin{align*}
     \E f(X) =&\ \E  \int_{0}^{1} f(X(t)) dt \\
     =&\ \E \int_{0}^{1} \int_{0}^{U(t)} f\big(X(t+u) \big) du d A(t) + \E \bigg(\int_{0}^{R_{a}(0)} f(X(t)) dt    - \int_{0}^{ R_{a}(1)} f(X(1+t)) dt \bigg)
\end{align*}
and the result follows by noting that since $(X,R_{a})$ is stationary, the joint distribution of $(X(t),R_{a}(t))$ is the same as $(X(1+t),R_{a}(1+t))$, so the two integrals cancel in expectation.
\finishproof

\subsection{Proof of Theorem~\ref{thm:workload:1}}
\label{app:workload:theorem}
\startproof{}
Given $h \in \lipone$, consider the Poisson equation \eqref{eq:stein:exponential} with $\theta = \delta^2$ and $\sigma^2 = \delta^2 \rho \E S (c_{U}^{2}+c_{S}^{2})$. Assuming that the BAR \eqref{eq:BAR2:workload} holds for $f_h(\TX)$, we set  $x = X$ and take expected values to arrive at 
\begin{align*}
    \E h(Y) - \E h(X) =&\ \E G_{Y} f_{h}(X) \\
    =&\  \E G_{Y} f_{h}(X) - \delta  \E  \big((\rho -  1(X>0))f_h'(\TX )\big) -  \E \int_{0}^{1}   \Delta f_h(\TX(t-))    d A(t).
\end{align*}
Taking the supremum over all $h \in \lipone$ and using  Proposition~\ref{prop:extraction:workload} yields 
\begin{align*}
     d_{W}(X,Y) \leq  \abs{\epsilon_{0}(f_h)} + \abs{\epsilon_{A}(f_h)}.
\end{align*}
We now   bound both  $\abs{\epsilon_{0}(f_h)}$ and $\abs{\epsilon_{A}(f_h)}$, and then show that  the BAR \eqref{eq:BAR2:workload} holds for $f_h(\TX)$. Recalling that 
\begin{align*}
\epsilon_{0}(f_h) =&\ \delta^3 \rho \E(R_{a} f_h''(\xi)) - \delta^2 \rho \E (1(X=0)R_{a} f_h''(\xi)),
\end{align*}
the bound on $\abs{\epsilon_{0}(f_h)}$ follows from $\norm{f_h''} \leq \delta^{-2}$ (Lemma~\ref{lem:exp:stein:factors}) and $\E R_{a} = \lambda \E U^2/2$  (Lemma~\ref{lem:relationships:workload}). Similarly, recall that 
\begin{align*}
        \epsilon_{A}(f_h) =&\ \frac{1}{6} \delta^3\E \int_{0}^{1} (S(t) - \rho U(t))^3 f_h'''(\xi(t))d A(t) \\
        &- \frac{1}{2} \delta^2\lambda \E (S -\rho  U)^2  \E \int_{0}^{1} \int_{0}^{U(t)} \big( X(t+u)-X(t-)\big) f_h'''(\xi(t+u))du d A(t).
\end{align*}
The first term is handled using the Stein factor bound $\norm{f_h'''}\leq 4/\sigma^2$, while the second term requires controlling increments of $X$ over an interarrival time. Namely,   
\begin{align*}
    \frac{1}{6} \delta^3\E \abs{\int_{0}^{1} (S(t) - \rho U(t))^3 f_h'''(\xi(t))d A(t)}  \leq  \frac{4 \delta^3\E \abs{S-\rho U}^{3}}{ 6\delta^2 \rho \E S (c_{U}^{2}+c_{S}^{2})} \E A(1) =  \frac{2\delta \E \abs{S-\rho U}^{3}}{ 3  (\E S)^2 (c_{U}^{2}+c_{S}^{2})},
\end{align*}
where the   equality is true because $\E A(1)=\lambda$, and 
\begin{align*}
    &\frac{1}{2} \delta^2\lambda \E (S -\rho  U)^2  \E \abs{\int_{0}^{1} \int_{0}^{U(t)} \big( X(t+u)-X(t-)\big) f_h'''(\xi(t+u))du d A(t)} \\
    \leq&\ 2\E \int_{0}^{1} \int_{0}^{U(t)} \abs{X(t+u)-X(t-)} du d A(t) \\
    =&\ 2 \E \int_{0}^{1} \int_{0}^{U(t)} \abs{(X(t-) + \delta S(t) - \delta u)^{+} - X(t-)} du d A(t)\\
    \leq&\ 2 \E \int_{0}^{1} \int_{0}^{U(t)} \delta (S(t) + U(t)) du d A(t)\\
   =&\ 2 \delta \E (US +U^2) \E A(1)  = 2 \delta (\E S + \lambda \E U^2),
\end{align*}
where the last inequality follows from $\abs{(x+y)^{+}-x^{+}} \leq y^{+}$.

It remains to show that \eqref{eq:BAR2:workload} holds with $f(\TX)=f_h(\TX)$.
By assumption \eqref{eq:moments:workload}, $U$ and $S$ have finite third moments. Since $\E R_{a} = \lambda \E U^2 /2 < \infty$ (Lemma~\ref{lem:relationships:workload}) and $\E X^2 < \infty$ \citep[Theorems~X.3.4 and X.2.1]{Asmu2003}, the fact that  $\abs{f_h(x)} \leq x^2/(2\theta)$ (Lemma~\ref{lem:exp:stein:factors}) yields $\E |f_h(\TX)|<\infty$; the fact that $\E |f_h'(\TX)|<\infty$ is argued similarly. To show that the jump term  in \eqref{eq:BAR2:workload} is well defined, 
note that 
\begin{align*}
    \E \abs{\int_{0}^{1}   \Delta f(\TX(t-))    d A(t)} =&\ \E \abs{\int_{0}^{1}   \int_{0}^{\delta S(t)} f_h'( X(t-)+u) du  d A(t)} \\
    \leq&\ \E \abs{\int_{0}^{1} \int_{0}^{\delta S(t)} \delta^{-2} \abs{X(t-)+u} du d A(t)}\\
    \leq&\ \E  \int_{0}^{1}  \delta^{-1} S(t) (X(t-)+S(t))  d A(t)  \\
    =&\ \delta^{-1}(\E S) \E  \int_{0}^{1}  X(t-)  d A(t)  + \delta^{-1} \E S^2 \E   A(1) 
\end{align*}
The second term on the right-hand side is finite because  $\E A(1) = \lambda$ (Lemma~\ref{lem:relationships:workload}). Letting $\tau_{m}$ be the time of the $m$th jump, the first term is bounded by 
\begin{align*}
  \E  \int_{0}^{1}  X(t-)  d A(t) \leq  \E  \big(A(1) \sup_{0 \leq t \leq 1} X(t)\big)\leq&\  \sqrt{\E (A(1))^2} \sqrt{\E  \Big(  X(0) + \delta \sum_{m=1}^{A(1)} S(\tau_{m}) \Big)^2}. 
\end{align*}
The right-hand side is finite because    $\E X^2 < \infty$,  $\E A^2(1) < \infty$ (since $\E U^2 < \infty$), and that $\E(\sum_{m=1}^{A(1)} S(\tau_{m}))^2 < \infty$ (since $A(1)$ is independent of  $S(\tau_{m})$).  
\finishproof

\section{The prelimit approach: supporting proofs}
\label{app:prelim:proofs}
\startproof{Proof of Lemma~\ref{lem:finite:busy:period} }
The first inequality follows from monotonicity: forcing an arrival to occur immediately
can only increase the initial workload and hence can only delay the first time the workload
empties.

For the second inequality, consider the system initialized with workload $x/\delta+S$
and residual interarrival time $U$. Treat the initial amount $x/\delta$ as low-priority work
and all remaining work, including the initial $S$ and all future arrivals, as high-priority
work. The low-priority work is processed only when no high-priority work is present, that
is, during idle periods of the high-priority system. Therefore, the low-priority work is
cleared once the cumulative idle time of the high-priority system exceeds $x/\delta$.
If $N_x$ denotes the number of idle periods required for this to occur, then
$\E N_x<\infty$ by \cite{Gut1974}. Since the high-priority system has i.i.d. cycles with
finite mean $\E(\bar B+\bar I)$, Wald's identity implies that the expected time until the
end of the $N_x$th idle period is $\E N_x\,\E(\bar B+\bar I)<\infty$. This upper-bounds
$\E(\E_{x+\delta S,U}B_0)$ and proves \eqref{eq:busy:finite}.
\finishproof

\subsection{Proving Proposition~\ref{prop:Mpoisson}}
\label{app:mhorproofs} 
We require several auxiliary lemmas.
\begin{lemma}
\label{lem:dx}
For any $h \in \lipone$ and $M > 0$,  
\begin{align*} 
 \partial_{x} F_{h}^{M}(z)  =&\  \E_{z} \int_{0}^{B_{0} \wedge M  }    h'(X(t))  dt, \quad z = (x,r_{a}) \in \mathbb{S}. 
\end{align*}
\end{lemma} 
\startproof{Proof of Lemma~\ref{lem:dx}} 
The proof is identical to that of Lemma~\ref{lem:abuse:d1}; see Appendix~\ref{app:poissonderive}.
\finishproof
\begin{lemma}
\label{lem:renewalidentity}
For any differentiable $g: \R_{+} \to \R$ with $\E |g(U)|, \E |g(R_{a})|$, $\E |g'(R_{a})| < \infty$, 
\begin{align*}
 \E g'(R_{a}) = \lambda \big(\E g(U) - \lim_{\epsilon \to 0} g(\epsilon)\big).
\end{align*}
\end{lemma}
\startproof{Proof of Lemma~\ref{lem:renewalidentity}}
Since $\E A(1) = \lambda$ (Lemma~\ref{lem:relationships:workload}), the result follows from the BAR in Lemma~\ref{lem:full:BAR:workload} with $f(Z(t)) = g(R_{a}(t))$ there. 
\finishproof
\noindent The next two lemmas are proved in Appendix~\ref{app:propaux}. 
\begin{lemma}
\label{lem:dr}
For any $h \in \lipone$ and almost all $M > 0$,   
\begin{align*}
-\delta \partial_{x} F_h^{M}(z) -  \partial_{r_{a}} F_h^{M}(z) = \E_{z} h(X(M)) - h(x), \quad z=(x,r_{a}) \in \mathbb{S}.
\end{align*}
In particular, $\partial_{r_{a}} F_h^{M}(z)$ is well defined for $z \in \mathbb{S}$. 
\end{lemma} 
\begin{lemma}
\label{lem:jump}
For any $h \in \lipone$ and $M > 0$, 
\begin{align*}
\lim_{\epsilon \to 0} F_h^{M}(x,\epsilon) = \E F_h^{M}(x+\delta S, U).
\end{align*}
\end{lemma}
\startproof{Proof of Proposition~\ref{prop:Mpoisson}}
For almost all $M > 0$, Lemma~\ref{lem:dr} says that
\begin{align}
-\delta \partial_{x} F_h^{M}(z) - \partial_{r_{a}} F_h^{M}(z) = \E_{z} h(X(M)) - h(x), \quad z \in \mathbb{S}. \label{eq:lem16}
\end{align}
Observe that $f(r_{a}) = F_h^{M}(x,r_{a})$ satisfies the conditions of Lemma~\ref{lem:renewalidentity}. Indeed, $\abs{\E f(U)}, \abs{\E f(R_{a})}< \infty$  because $M < \infty$. Furthermore, $\abs{\E f'(R_{a})} < \infty$ follows from the expression for $\partial_{r_{a}} F_h^{M} (z)$ in Lemma~\ref{lem:dr}, together with the observation that $\abs{\partial_{x} F_h^{M}(z)} \leq M$, which follows from Lemma~\ref{lem:dx}. Setting $r_{a} = R_{a}$ in \eqref{eq:lem16} and taking expected values yields 
\begin{align*}
 \E \big( \E_{z,R_{a}} h(X(M)) - h(x) \big) =&\ -\delta  \E \partial_{x} F_{h}^{M}(x,R_{a}) - \E \partial_{r_{a}} F_h^{M}(x,R_{a}) \\
=&\ -\delta   \E \partial_{x} F_{h}^{M}(x,R_{a}) - \lambda \Big(   \E F_{h}^{M}(x,U) - \lim_{\epsilon \to 0} F_{h}^{M}(x,\epsilon )\Big) \\
=&\  -\delta   \E \partial_{x} F_{h}^{M}(x,R_{a}) - \lambda \Big(  \E F_{h}^{M}(x,U) - \E F_{h}^{M}(x+\delta S,U) \Big),
\end{align*}
where the second  and third equalities  follow from Lemmas~\ref{lem:renewalidentity} and \ref{lem:jump}, respectively.  
\finishproof

\subsubsection{Auxiliary lemma proofs}
\label{app:propaux}  
\startproof{Proof of Lemma~\ref{lem:dr}}
We define $\tilde h(x) = h(x) - \E h(X)$ for convenience, in which case 
\begin{align*}
F_h^{M}(z) = \int_{0}^{M} \E_{z} \tilde h(X(t)) dt, \quad z \in \mathbb{S}.
\end{align*}
Our goal is to prove that 
\begin{align}
-\delta \partial_{x} F_h^{M}(z) -  \partial_{r_{a}} F_h^{M}(z) = \E_{z} h(X(M)) - h(x), \quad z=(x,r_{a}) \in \mathbb{S}. \label{eq:drgoal}
\end{align}
Fix $z  = (x,r_{a}) \in \mathbb{S}$ and suppose first that $x = 0$. On one hand, 
\begin{align*}
F_h^{M+\epsilon}(0,r_{a}+\epsilon) =&\  F_{h}^{M}(0,r_{a}+\epsilon)+ \int_{M}^{M+\epsilon} \E_{0,r_{a}+\epsilon} \tilde h(X(t)) dt,
\end{align*}
and on the other, 
\begin{align*}
F_h^{M+\epsilon}(0,r_{a}+\epsilon) = \int_{0}^{\epsilon} \E_{0,r_{a}+\epsilon} \tilde h(X(t)) dt +  \int_{0}^{M} \E_{0,r_{a}}  \tilde h(X(t)) dt  = \epsilon \tilde h(0) + F_{h}^{M} (0,r_{a}),
\end{align*}
  Equating the two expressions and dividing both sides by $\epsilon$ yields
\begin{align*}
\lim_{\epsilon \to 0} \frac{1}{\epsilon}\big(F_{h}^{M}(0,r_{a}+\epsilon) - F_{h}^{M}(0,r_{a})\big) =&\  \tilde h(0) - \lim_{\epsilon \to 0} \frac{1}{\epsilon} \int_{M}^{M+\epsilon} \E_{0,r_{a}+\epsilon} \tilde h(X(t)) dt.
\end{align*}
The left-hand side equals $\delta \partial_{x} F_h^{M}(z) +\partial_{r_{a}} F_h^{M}(z) =   \partial_{r_{a}} F_h^{M}(z)$, since $\partial_{x} F_h^{M}(0,r_{a}) = 0$ by Lemma~\ref{lem:dx}. Thus, to prove \eqref{eq:drgoal} when $x = 0$, it suffices to show that 
\begin{align}
& \frac{1}{\epsilon} \int_{M}^{M+\epsilon} \big(\E_{0,r_{a}+\epsilon} \tilde h(X(t))  - \E_{0,r_{a}} \tilde h(X(t)) \big)dt \notag \\
=&\ \frac{1}{\epsilon} \int_{M}^{M+\epsilon}  \E_{0,r_{a}} \big( \tilde h(X(t-\epsilon))  -\tilde h(X(t)) \big)dt \to 0 \label{eq:intMeps}
\end{align}
as $\epsilon \to 0$, which implies that 
\begin{align*}
\lim_{\epsilon \to 0} \frac{1}{\epsilon} \int_{M}^{M+\epsilon} \E_{0,r_{a}+\epsilon} \tilde h(X(t)) dt = \lim_{\epsilon \to 0} \frac{1}{\epsilon} \int_{M}^{M+\epsilon} \E_{0,r_{a}} \tilde h(X(t)) dt = \E_{0,r_{a}} \tilde h(X(M)).
\end{align*}
Observe that   $\abs{X(t-\epsilon) - X(t)}$ is bounded by  the workload processed during $[t-\epsilon,t]$, which is at most $\delta \epsilon$,  plus any new work that arrives during $[t-\epsilon,t]$. Letting $A([t_1,t_2])$ denote the number of customers arriving during $[t_1,t_2]$, Wald's identity says that the expected workload to arrive during $[t_1,t_2]$ equals $\E S \E A([t_1,t_2])$. Thus, to prove \eqref{eq:intMeps}, we observe that for any $h \in \lipone$ and  for all $t \in [M,M+\epsilon]$, 
\begin{align*}
\E_{0,r_{a}} \big| \tilde h(X(t-\epsilon))  -\tilde h(X(t)) \big| \leq&\  \E_{0,r_{a}} \abs{X(t-\epsilon) - X(t)}  \\
\leq&\ \delta \epsilon + \E_{0,r_{a}} \big( \delta \E S \E (A([t - \epsilon, t])) \big) \\
\leq&\ \delta \epsilon +  \delta \E S \E_{0,r_{a}} \big( A([M - \epsilon, M+\epsilon])\big).
\end{align*}
 It suffices to argue that  the right-hand side goes to zero as $\epsilon \to 0$.   By the dominated convergence theorem, 
\begin{align*}
\lim_{\epsilon \to 0} \E_{0,r_{a}} \big( A([M - \epsilon, M+\epsilon])\big) = \E_{0,r_{a}} \big( A([M, M])\big),
\end{align*}
which equals the expected number of arrivals at time $M$. The right-hand side may be non-zero if the distribution   of $U$ has  point masses.  However, since the number of point masses is at most countable,   then $ \E_{0,r_{a}} \big( A([M, M])\big) = 0$ for all but at most countably many $M$. This proves \eqref{eq:drgoal} when $x = 0$. 

The case when $x> 0$ follows similarly. We repeat the arguments, highlighting the differences. Given $z = (x,r_{a})$, fix $\epsilon < x/\delta$. Then 
\begin{align*}
F_h^{M+\epsilon}(x,r_{a}+\epsilon) =&\  
 F_{h}^{M}(x,r_{a}+\epsilon) + \int_{M}^{M+\epsilon}  \E_{x,r_{a}+\epsilon}\tilde h(X(t)) dt  
\end{align*}
and 
\begin{align*}
F_h^{M+\epsilon}(x,r_{a}+\epsilon) = \int_{0}^{\epsilon} \E_{x,r_{a}+\epsilon} \tilde h(X(t)) dt + F_h^{M}(x-\delta\epsilon,r_{a}).
\end{align*}
Equating both expressions,  subtracting $F_h^{M}(x,r_{a})$ from each side, and dividing by $\epsilon$ yields 
\begin{align*}
& \frac{1}{\epsilon}	 \big( F_{h}^{M}(x,r_{a}+\epsilon) - F_h^{M}(x,r_{a})\big) \\
=&\  \frac{1}{\epsilon}	 \big(F_h^{M}(x-\delta\epsilon,r_{a})- F_h^{M}(x,r_{a})\big) + \frac{1}{\epsilon}\int_{0}^{\epsilon} \E_{x,r_{a}+\epsilon} \tilde h(X(t)) dt - \frac{1}{\epsilon}\int_{M}^{M+\epsilon}  \E_{x,r_{a}+\epsilon}\tilde h(X(t)) dt.
\end{align*}
We now argue that each of the terms on the right-hand side has a well-defined limit as $\epsilon \to 0$, implying that the left-hand side converges to $\partial_{x} F_h^{M}(z)$, which is itself well defined. 
The first term on the right-hand side converges to  $-\partial_{x} F_h^{M}(z)$, which we know exists for all $z \in \mathbb{S}$ by Lemma~\ref{lem:dx}. Furthermore, 
\begin{align*}
    \lim_{\epsilon \to 0} \frac{1}{\epsilon}\int_{0}^{\epsilon} \E_{x,r_{a}+\epsilon} \tilde h(X(t)) dt  = \tilde h(x) \quad \text{ and } \quad  \lim_{\epsilon \to 0} \frac{1}{\epsilon}\int_{M}^{M+\epsilon}  \E_{x,r_{a}+\epsilon}\tilde h(X(t)) dt = \E_{x,r_{a}} \tilde h(X(M)).
\end{align*}
The first equality is straightforward because no arrival occurs during $[0,\epsilon]$, while the second equality is proved the same way as \eqref{eq:intMeps}.
\finishproof

\startproof{Proof of Lemma~\ref{lem:jump}}
Define $\tilde h(x) = h(x) - \E h(X)$  and consider first the case when $x = 0$. Then 
\begin{align*}
F_h^{M}(0,\epsilon) =&\ \int_{0}^{\epsilon} \E_{0,\epsilon} \tilde h(X(t)) dt + \int_{\epsilon}^{M} \E_{0 ,\epsilon} \tilde h(X(t)) dt  \\
=&\ \int_{0}^{\epsilon} \E_{0,\epsilon} \tilde h(X(t)) dt + \E \int_{0}^{M-\epsilon}   \E_{\delta S, U}  \tilde h(X(t)) dt,
\end{align*}
where the outer expectation is with respect to $U$ and $S$. Taking $\epsilon \to 0$, the left-hand side converges to $\lim_{\epsilon \to 0} F_h^{M}(0,\epsilon)$ while the right-hand side converges to 
\begin{align*}
    \lim_{\epsilon \to 0} \E \int_{0}^{M-\epsilon}   \E_{\delta S, U}  \tilde h(X(t)) dt = \E \int_{0}^{M}   \E_{\delta S, U}  \tilde h(X(t)) dt - \lim_{\epsilon \to 0}\E \int_{M-\epsilon}^{M}  \E_{\delta S, U}  \tilde h(X(t)) dt.
\end{align*}
The first term equals $\E F_h^{M}(\delta S, U)$ while the second term is zero because $h \in \lipone$. 
Now suppose that $x >  0$ and take $\epsilon < x/\delta$. Arguing as before, 
\begin{align*}
F_h^{M}(x ,\epsilon) =&\  \int_{0}^{\epsilon} \E_{x,\epsilon} \tilde h(X(t)) dt +\E F_h^{M-\epsilon}(x-\delta \epsilon + \delta S, U).  
\end{align*}
To conclude, we use the fundamental theorem of calculus to write
\begin{align*}
\E F_h^{M-\epsilon}(x-\delta \epsilon + \delta S, U) 
=&\ \E F_h^{M-\epsilon}(x + \delta S, U) + \E \int_{0}^{-\delta \epsilon} \partial_{x} F_h^{M-\epsilon}(x+v + \delta S, U) dv.
\end{align*}
The second term on the right-hand side converges to zero because $\big|\partial_{x}F_h^{M-\epsilon}(z)\big|\leq M$ due to Lemma~\ref{lem:dx}. The first term  converges to $\E F_h^{M}(x + \delta S, U)$ because 
\begin{align*}
   \lim_{\epsilon \to 0} \E \int_{0}^{M-\epsilon}   \E_{x+\delta S, U}  \tilde h(X(t)) dt = \E \int_{0}^{M}   \E_{x+\delta S, U}  \tilde h(X(t)) dt - \lim_{\epsilon \to 0}\E \int_{M-\epsilon}^{M}  \E_{x+\delta S, U}  \tilde h(X(t)) dt,
\end{align*}
and the second term equals zero since $h \in \lipone$. 
\finishproof
\subsection{Proofs of Lemmas~\ref{lem:abuse:d1} and \ref{lem:infpoisson}}
\label{app:poissonderive}
We recall the synchronous coupling $\{Z^{(\epsilon)}(t): t \geq 0\}$ defined in Section~\ref{sec:mhorizon}.
\startproof{Proof of Lemma~\ref{lem:abuse:d1}}
First, observe that 
\begin{align*}
   &\frac{1}{\epsilon} \E \Big( \int_{0}^{\infty} \big( \E_{x+\epsilon,T}  h(X(t)) - \E_{x,T}  h(X(t)) \big) dt \Big) \\
   =&\  \frac{1}{\epsilon} \E \Big(\E_{x,T} \int_{0}^{B_{0}  }  \big( h(X^{(\epsilon)}(t)) - h(X(t)) \big) dt \Big)  +  \frac{1}{\epsilon} \E \Big(\E_{x,T} \int_{B_{0} }^{B^{(\epsilon)}_{0}  }  \big( h(X^{(\epsilon)}(t)) - h(X(t)) \big) dt \Big).
\end{align*}
Note that $\abs{h(X^{(\epsilon)}(t)) - h(X(t))}/\epsilon \leq \norm{h'} \leq 1$  and $B^{(\epsilon)}_{0} \to B_{0}$ as $\epsilon \to 0$. Also note that  for all $\epsilon < 1$, 
\begin{align*}
    \E_{x,T} B^{(\epsilon)}_{0} =  \E_{x+\epsilon,T} B_{0} \leq \E_{x+1,T} B_{0} \leq \E \big(\E_{x+1+\delta S,U} B_{0}\big) < \infty,
\end{align*}
where the second-last inequality follows from the fact that the busy period starting at state $(x+1,T)$ is  made longer if the next arrival happens immediately. The DCT then implies that  
\begin{align*}
    \frac{1}{\epsilon} \E \Big(\E_{x,T} \int_{0}^{B_{0}  }  \big( h(X^{(\epsilon)}(t)) - h(X(t)) \big) dt \Big) \to&\ \E \Big(\E_{x,T} \int_{0}^{B_{0}}    h'(X(t))  dt \Big),\\
    \frac{1}{\epsilon} \E \Big(\E_{x,T} \int_{B_{0} }^{B^{(\epsilon)}_{0}  }  \big( h(X^{(\epsilon)}(t)) - h(X(t)) \big) dt \Big) \to&\ 0.
\end{align*}
\finishproof
 \startproof{Proof of Lemma~\ref{lem:infpoisson}}
Fix $h \in \lipone$. Let  $\mu_{R_a}$ and $\mu_{X\mid r_a}$ denote the law of
$R_a$ and conditional law of $X$ given $R_a=r_a$, respectively. 
We first prove \eqref{eq:lim1}.  Note that 
\begin{align*}
    &\E h(X)   =   \int_{0}^{\infty}\int_{0}^{\infty} \E_{y,r_{a}}h(X(M)) d \mu_{X|r_{a}} (y) d \mu_{R_{a}}(r_{a}), \\
    &\E \big( \E_{x,R_{a}} h(X(M)) \big)  = \int_{0}^{\infty} \E_{x,r_a} h(X(M))   d \mu_{R_{a}}(r_{a})= \int_{0}^{\infty}\int_{0}^{\infty} \E_{x,r_a} h(X(M))  d \mu_{X|r_{a}} (y) d \mu_{R_{a}}(r_{a}).
\end{align*}
The first equality is true because $\E h(X)$ coincides with $\E h(X(M))$ if $Z(0)$ is initialized according to $Z$. It follows that 
\begin{align*}
\E \big( \E_{x,R_{a}} h(X(M)) \big) - \E h(X) =&\  \int_{0}^{\infty}\int_{0}^{\infty} \big(\E_{x,r_a} h(X(M))   - \E_{y,r_a} h(X(M)) \big) d \mu_{X|r_{a}} (y) d \mu_{R_{a}}(r_{a}).
\end{align*}
Since $h \in \lipone$, we can use the synchronous coupling defined in \eqref{eq:couplinggap} to compare two workload processes initialized at $(x,r_a)$ and $(y,r_a)$. Under this coupling, their workloads differ by at most $\abs{x-y}$ until they couple. Coupling occurs no later than the first time the larger initial workload process empties. Therefore,
\begin{align*}
\abs{\E_{x,r_a} h(X(M))  - \E_{y,r_a} h(X(M)) }
\leq
\abs{x-y}\Prob_{x \vee y, r_a} \big(  B_{0} > M\big),
\end{align*}
where the probability on the right-hand side is the probability that the coupled processes have not yet coupled by time $M$. Thus,
\begin{align*}
\lim_{M \to \infty}  \abs{\E \big( \E_{x,R_{a}} h(X(M)) \big) - \E h(X)}  \leq   \lim_{M \to \infty}  \E \Big( \abs{x-X}\Prob_{x \vee X, R_{a}} \big(   B_{0} > M \big) \Big) = 0,
\end{align*}
where the expectation on the right-hand side is taken with respect to the joint law of $(X,R_{a})$.
The last equality follows from the DCT because $\E X < \infty$, and because $\lim_{M \to \infty} \Prob_{x \vee x',r_{a}}\big(   B_{0} > M \big) = 0$ for any $x,x',r_{a} > 0$ by \eqref{eq:busy:finite}. To prove \eqref{eq:lim2}, one can reuse the arguments used to prove Lemma~\ref{lem:abuse:d1} to show that 
\begin{align*}
    \partial_{x} F_{h}^{M}(x,r_{a}) =    \E_{x,r_{a}} \int_{0}^{B_{0} \wedge M  }    h'(X(t))  dt \to \partial_{x} F_{h}(x,r_{a}) \quad \text{ as $M \to \infty$ for all $(x,r_{a}) \in \mathbb{S}$,} 
\end{align*}
and also that $\lim_{M \to \infty}  \E \partial_{x} F_{h}^{M}(x,R_{a})  = \E \lim_{M \to \infty}\partial_{x} F_{h}^{M}(x,R_{a})$.  Lastly, we prove \eqref{eq:lim3}. Similar to the way we argued \eqref{eq:findiff}, 
\begin{align*}
    F_{h}^{M}(x+\epsilon,r_{a}) -  F_{h}^{M}(x,r_{a}) =&\ \E_{x,r_{a}} \int_{0}^{B^{(\epsilon)}_{0}  \wedge M }  \big( h(X^{(\epsilon)}(t)) - h(X(t)) \big) dt \\
    \to&\ F_{h}(x+\epsilon,r_{a}) -  F_{h}(x,r_{a}), \quad \text{ as $M \to \infty$ for all $(x,r_{a}) \in \mathbb{S}$}.
\end{align*}
Let $\hat h(x) = x$ and observe that   since  $h \in \lipone$ and $X^{(\epsilon)}(t) \geq X(t)$, then
\begin{align*}
 \abs{F_{h}^{M}(x+\epsilon,r_{a}) -  F_{h}^{M}(x,r_{a})} \leq&\   \E_{x,r_{a}} \int_{0}^{B^{(\epsilon)}_{0} \wedge M  }   \big(X^{(\epsilon)}(t) - X(t) \big)  dt \leq F_{\hat h}(x+\epsilon,r_{a}) -  F_{\hat h}(x,r_{a}).
\end{align*}
It remains to show   that $\E \big(F_{\hat h}(x+\delta S,U) -  F_{\hat h}(x,U)\big) < \infty$, because then we can use the DCT to conclude \eqref{eq:lim3}. The finiteness of this expectation follows from
\begin{align*}
   \E \big(F_{\hat h}(x+\delta S,U) -  F_{\hat h}(x,U)\big)  =&\  \lim_{M \to \infty}  \E \big(F_{\hat h}^{M}(x+\delta S,U) -  F_{\hat h}^{M}(x,U) \big) \\
   =&\ \lim_{M \to \infty} \E \big( \E_{x,R_{a}} X(M) - x \big) + \lim_{M \to \infty}\delta \E \partial_{x} F_{\hat h}^{M}(x,R_{a}),
\end{align*}
where the first equality is due to the monotone convergence theorem, since   $F_{\hat h}^{M}(x+\epsilon,r_{a}) -  F_{\hat h}^{M}(x,r_{a})$ is increasing in $M$ and is nonnegative for any $(x,r_{a}) \in \mathbb{S}$, and the second equality is due to  \eqref{eq:Mpoisson}. The right-hand side is finite by \eqref{eq:lim1} and \eqref{eq:lim2}.
\finishproof

\subsection{Section~\ref{sec:taylor} Proofs}
\label{app:taylorproofs}
\startproof{Proof of Lemma~\ref{lem:relate}}
Recall that $J(x,r_{a})  = - (x \wedge \delta r_{a}) + \delta S'$. It follows that 
\begin{align*}
 F_h(x+ \delta s,r_{a}) - F_h(x,r_{a})  =&\ \int_{0}^{\infty} \big(\E_{x+\delta s, r_{a}} h(X(t)) - \E_{x,r_{a}} h(X(t)) \big) dt \\
=&\ \int_{0}^{r_{a}} \Big(  h\big((x + \delta s-\delta t)^{+}\big) -  h\big((x  -\delta t)^{+}\big) \Big) dt \\
&+ \E \Big( F_h\big( x+\delta s +J(x+\delta s,r_{a}), U\big)-   F_h\big( x +J(x,r_{a}), U\big) \Big).
\end{align*}
To conclude, note that 
\begin{align*}
 & \E \Big( F_h\big( x+\delta s +J(x+\delta s,r_{a}), U\big)-   F_h\big( x +J(x,r_{a}), U\big) \Big) \\
 =&\  \E \Big(   F_h\big( x+\delta s +J(x,r_{a}), U\big)  -   F_h\big( x +J(x,r_{a}), U\big) \Big) \\
&+ \E \Big(    F_h\big( x+\delta s +J(x+\delta s,r_{a}), U\big) -   F_h\big( x+\delta s +J(x,r_{a}), U\big)\Big).
\end{align*}
Using the fundamental theorem of calculus, together with Lemma~\ref{lem:abuse:d1}, which shows that $\partial_{x}  \E  F_{h}(x+ \delta S,U) =  \E  \partial_{x}  F_{h}(x+ \delta S,U)$, we arrive at
\begin{align*}
& \E \Big(    F_h\big( x+\delta s +J(x+\delta s,r_{a}), U\big) -   F_h\big( x+\delta s +J(x,r_{a}), U\big)\Big) \\ 
=&\ \E \Big(  \int_{-x \wedge(\delta r_{a})}^{-(x+\delta s) \wedge(\delta r_{a})} \partial_{x}  F_h(x+\delta s +v +\delta S',U) dv\Big)\\
=&\ \E^{S'} \Big(  \int_{-x \wedge(\delta r_{a})}^{-(x+\delta s) \wedge(\delta r_{a})} \E^{U} \partial_{x}  F_h(x+\delta s +v +\delta S',U) dv\Big)\\
=&\ \E^{S'} \Big(  \int_{-x \wedge(\delta r_{a})}^{-(x+\delta s) \wedge(\delta r_{a})} \partial_{x}   \E^{U} F_h(x+\delta s +v +\delta S',U) dv\Big).
\end{align*}
Interchanging $\E^{U}$ with the integral in the second equality is justified by the Fubini-Tonelli theorem because $\E^{S'}\E^{U} \abs{\partial_{x}  F_h(x+\delta S',U)}  \leq \E^{S'}\E^{U}\E_{x +\delta S',U} B_{0} < \infty$ for all $x \geq 0$ by Lemma~\ref{lem:abuse:d1} and \eqref{eq:busy:finite}. 

\finishproof

\subsubsection{Proving Lemma~\ref{lem:taylor:final}}
\label{app:taylorlemma}
We recall that $\bar F_h'(x) = \partial_{x}  \E F_h(x,R_{a})$ and that   $\bar F_h''(x)$ and $\bar F_h'''(x)$ are assumed to exist. We recall  \eqref{eq:backpoisson}, or 
\begin{align}
& \E h(X) - \E h(x + J(x,R_{a}'))  \notag \\
=&\ -\delta \E \bar F_h'(x + J(x,R_{a}')) + \lambda \E \big(F_h(x+ \delta S,R_{a}') - F_h(x,R_{a}') \big) - \E \big(\epsilon(x,R_{a}',S)\big), \label{eq:backrecall}
\end{align}
where $J(x,r_{a})  = - (x \wedge \delta r_{a}) + \delta S'$. The following lemma expands the first two terms on the right-hand side of \eqref{eq:backrecall}. We prove it after proving Lemma~\ref{lem:taylor:final}.
\begin{lemma}
For any $x \geq 0$, 
\label{lem:interchangeexpansion}
\begin{align*} 
&\bar F_h'(x + J(x,R_{a}')) = \bar F_h'(x) + \delta(S'-R_{a}')\bar F_h''(x) + 1(\delta R_{a}' < x) \int_{0}^{\delta(S'-R_{a}')} \int_{0}^{v} \bar F_h'''(x+u) du  dv \\
& \hspace{3cm} + 1(\delta R_{a}' \geq x) \big( \bar F_{h}'(\delta S')  - \bar F_h'(x) - \delta(S'-R_{a}')\bar F_h''(x)\big)  \\ 
&\E \big(F_h(x+ \delta S,R_{a}') - F_h(x,R_{a}') \big) = \delta  \E S \bar F_h'(x) + \frac{1}{2} \delta^2 \E S^2 \bar F_h''(x) + \E \int_{0}^{\delta S} (\delta S - v)\int_{0}^{v} \bar F_h'''(x+u) du dv, 
\end{align*}
\end{lemma}
\startproof{Proof of Lemma~\ref{lem:taylor:final} }
Recall that  $\lambda \E S= \rho$.   Combining Lemma~\ref{lem:interchangeexpansion} with \eqref{eq:backrecall} yields 
\begin{align*}
 &\E h(X) - \E h(x + J(x,R_{a}')) \\
 =&\ - \delta  \big( \bar F_h'(x) + \delta \E (S'-R_{a}')\bar F_h''(x) \big)  + \lambda \big( \delta  \E S \bar F_h'(x) + \frac{1}{2} \delta^2 \E S^2 \bar F_h''(x) \big)  \\
 &-\delta \E \Big( 1(\delta R_{a}' \geq x) \big( \bar F_{h}'(\delta S')  - \bar F_h'(x) - \delta(S'-R_{a}')\bar F_h''(x)\big) \Big) \\
 & -\delta \E \Big(  1(\delta R_{a}' < x) \int_{0}^{\delta(S'-R_{a}')} \int_{0}^{v} \bar F_h'''(x+u) du  dv \Big) \\
 &+ \lambda \E \int_{0}^{\delta S} (\delta S - v)\int_{0}^{v} \bar F_h'''(x+u) du dv .
\end{align*}
Using the facts that $\lambda \E S = \rho$, $\lambda \E U = 1$, and that $\E R_{a}' = \lambda \E U^2 /2$, we see that the first line on the right-hand side equals 
\begin{align*}
    -\delta (1-\rho) \bar F_h'(x) + \frac{1}{2} \delta^2 \big( \lambda \E S^2 - 2 \lambda \E U \E S' + \lambda \E U^2  \big) \bar F_h''(x).
\end{align*}
Let us call this term $G_{Y_2} \bar F_h(x)$.
Since $ \bar F_h'(0) = 0$ due to Lemma~\ref{lem:abuse:d1} (because $B_0=0$ if the initial workload $X(0)=0$), our assumptions that $\E \abs{ \bar F_h'(Y_2)}, \E \abs{ \bar F_h''(Y_2)} < \infty$ and integration by parts yield   $\E G_{Y_2} \bar F_h(Y_2) = 0$.  
\finishproof

\startproof{Proof of Lemma~\ref{lem:interchangeexpansion}}
The expression for $\bar F_h'(x + J(x,R_{a}'))$ follows from the facts that
\begin{align*}
    \bar F_h'(x + J(x,R_{a}')) =&\  1(\delta R_{a}' \geq x) \bar F_h'(\delta S') + 1(\delta R_{a}' < x) \bar F_h'(x -\delta R_{a}' + \delta S')
\end{align*}
and, for all $x > \delta R_{a}'$,
\begin{align*}
     \bar F_h'(x -\delta R_{a}' + \delta S')  =   \bar F_h'(x) + \delta(S'-R_{a}')\bar F_h''(x) +  \int_{0}^{\delta(S'-R_{a}')} \int_{0}^{v} \bar F_h'''(x+u) du  dv.
\end{align*}
Next, we argue that 
\begin{align}
\E \big(F_h(x+ \delta S,R_{a}') - F_h(x,R_{a}') \big)   = \E \int_{0}^{\delta S}  \bar F_h'(x + v) dv, \label{eq:somecare}
\end{align}
so that the expression for $\E \big(F_h(x+ \delta S,R_{a}') - F_h(x,R_{a}') \big) $ also follows from Taylor expansion of the integrand around $x$. 
To prove \eqref{eq:somecare}, note that 
\begin{align*}
\E \big(F_h(x+ \delta S,R_{a}') - F_h(x,R_{a}') \big) =&\ \E \int_{0}^{\delta S} \partial_{x}F_{h}(x + v,R_{a}') dv \\
=&\  \E^{S} \int_{0}^{\delta S} \E^{R_{a}'} \partial_{x}F_{h}(x + v,R_{a}') dv  \\
=&\  \E^{S} \int_{0}^{\delta S}  \partial_{x}\E F_{h}(x + v,R_{a}') dv \\
=&\  \E \int_{0}^{\delta S}  \bar F_h'(x + v) dv.
\end{align*}
The first and  second-last equalities follows from Lemma~\ref{lem:abuse:d1}.   Once we justify the interchange of the integral and expectation  in the second equality  using the Fubini-Tonelli theorem, \eqref{eq:somecare} will follow.  Let $\hat h(x) = x$. Using the form of $\partial_{x} F_{h}(x,r_{a})$ from Lemma~\ref{lem:abuse:d1}, it follows that for any $h \in \lipone$,  
\begin{align*}
\abs{\partial_{x} F_{h}(x,r_{a})} \leq \E_{x,r_{a}} B_{0} = \E_{x,r_{a}} \int_{0}^{B_{0}}     \hat h'(X(t))   dt = \partial_{x} F_{\hat h}(x,r_{a}).
\end{align*}
Thus, 
\begin{align*}
 \E \int_{0}^{\delta S}  \abs{\partial_{x}F_{h}(x + v,R_{a}')} dv \leq \E \int_{0}^{\delta S}   \partial_{x}F_{\hat h}(x + v,R_{a}') dv = \E \big(F_{\hat h}(x+ \delta S,R_{a}') - F_{\hat h}(x,R_{a}') \big),
\end{align*}
and the right-hand side is finite because the right-hand side of \eqref{eq:lim3} in Lemma~\ref{lem:infpoisson} is finite.
\finishproof 
 \section{$G/G/1$ workload Stein factor bounds: supporting proofs}
 For the entirety of this section we fix   $h(x) = x$ and assume that  $\overline{\eta} = \sup\{\eta(x) : x \geq 0\}< \infty$.
 \subsection{Second-order bounds (Proof of Lemma~\ref{lem:stein2})}
 \label{app:stein2:proofs}
Recall that $G(x) = \Prob(U \leq x)$. The following auxiliary lemma is needed to prove Lemma~\ref{lem:stein2}. 
\begin{lemma}
\label{lem:threebounds}
For any $\epsilon > 0$ and  $ (x,r_{a}) \in \mathbb{S}$ with $r_{a} < x/\delta$, 
\begin{align}
\frac{1}{\epsilon}\Prob_{x,r_{a}}(R_{a}(B_{0}) < \epsilon/\delta) \leq&\ \overline{\eta}/\delta, \label{eq:bound} \\
\lim_{\epsilon \to 0}  \frac{1}{\epsilon}\Prob_{x,r_{a}}(R_{a}(B_{0}) < \epsilon/\delta) =&\ \frac{1}{\delta} \E_{x,r_{a}} \eta(\alpha(B_{0})). \label{eq:limit}
\end{align} 
\end{lemma}  
\startproof{Proof of Lemma~\ref{lem:threebounds}}
Let $U_{n}$  denote the interarrival time of the $n$th customer, let $W_{0} = V(0)$, and let $W_{n} = V(U_{1} + \cdots + U_{n})$ be the workload in the system right after the $n$th customer arrives, which includes the workload brought  by the $n$th customer. Let 
\begin{align*}
\sigma = \min \{ n \geq 1 : U_{n} > W_{n-1} \}
\end{align*}
be the number of customers served in the first busy period $[0,B_{0}]$. Now assuming that $Z(0) = (x,r_{a}) \in \mathbb{S}$ with $r_{a} < x/\delta$, it must be that $\sigma > 1$, because $W_{0} = x/\delta$ and $U_{1} = r_{a}$. Since $\{R_{a}(B_{0}) \leq \epsilon/\delta\} = \{U_{\sigma} \leq W_{\sigma -1} + \epsilon/\delta\}$, it follows that 
\begin{align*}
 \frac{1}{\epsilon} \Prob_{x,r_{a}}(R_{a}(B_{0}) \leq \epsilon/\delta)  =&\ \frac{1}{\epsilon} \sum_{n=2}^{\infty} \Prob_{x,r_{a}}\big( U_{n} \leq W_{n-1} + \epsilon/\delta \big| \sigma = n\big) \Prob_{x,r_{a}}( \sigma = n) \\ 
=&\ \frac{1}{\epsilon} \sum_{n=2}^{\infty} \E_{x,r_{a}}\bigg[ \Prob_{x,r_{a}}\big( U_{n} \leq W_{n-1} + \epsilon/\delta \big|  \sigma = n, W_{n-1}\big) \bigg| \sigma = n \bigg] \Prob_{x,r_{a}}( \sigma = n).
\end{align*}  
To proceed, note that $\{\sigma = n \} = \{ U_{1} \leq W_{0}, \ldots, U_{n-1} \leq W_{n-2}, U_{n} > W_{n-1}\}$ for any $n \geq 1$, implying that  for any $n \geq 2$, 
\begin{align*}
&\Prob_{x,r_{a}}\big( U_{n} \leq W_{n-1} + \epsilon/\delta \big|  \sigma = n, W_{n-1}\big) \\
=&\ \Prob\big( U_{n} \leq W_{n-1} + \epsilon/\delta \big| W_{0}=x/\delta, U_{1}=r_{a}, \sigma = n, W_{n-1}\big) \\
=&\ \Prob\big( U_{n} \leq W_{n-1} + \epsilon/\delta \big|  W_{0}=x/\delta, U_{1}=r_{a}, U_{1} \leq W_{0}, \ldots, U_{n-1} \leq W_{n-2}, U_{n} > W_{n-1}, W_{n-1}\big) \\
=&\ \Prob\big( U_{n} \leq W_{n-1} + \epsilon/\delta \big|   U_{n} > W_{n-1}, W_{n-1}\big)\\
=&\ \Prob\big( U \leq W_{n-1} + \epsilon/\delta \big|   U  > W_{n-1}, W_{n-1}\big) \\
=&\ \frac{G(W_{n-1}+\epsilon/\delta) - G(W_{n-1})}{1  - G(W_{n-1}) },
\end{align*} 
and therefore 
\begin{align}
\frac{1}{\epsilon} \Prob_{x,r_{a}}(R_{a}(B_{0}) \leq \epsilon/\delta)  =&\ \frac{1}{\epsilon} \sum_{n=1}^{\infty} \E_{x,r_{a}}\bigg[ \frac{G(W_{n-1}+\epsilon/\delta) - G(W_{n-1})}{1  - G(W_{n-1}) } \bigg| \sigma = n \bigg] \Prob_{x,r_{a}}( \sigma = n)  \notag \\
=&\ \frac{1}{\epsilon}  \E_{x,r_{a}}\bigg[ \frac{G(W_{\sigma-1}+\epsilon/\delta) - G(W_{\sigma-1})}{1  - G(W_{\sigma-1}) } \bigg]. \label{eq:wsigma}
\end{align}
To prove \eqref{eq:bound}, observe that the right-hand side of \eqref{eq:wsigma} is bounded by $\overline{\eta}/\delta$ because by the mean value theorem,
\begin{align*}
 \frac{G(w+\epsilon/\delta) - G(w)}{1  - G(w) }  =  \frac{\epsilon}{\delta}\frac{G'(\xi)}{1-G(w)} = \frac{\epsilon}{\delta} \eta(\xi) \frac{1-G(\xi)}{1-G(w)} \leq \frac{\epsilon}{\delta} \overline{\eta}
\end{align*}
for some $\xi \in [w,w+\epsilon/\delta]$, where the last inequality follows from  $\xi \geq w$  and $\eta(x) \leq \overline{\eta}$.  Once we observe that   $W_{\sigma -1} = \alpha(B_{0})$, then  \eqref{eq:limit} follows from taking $\epsilon \to 0$ in \eqref{eq:wsigma} and applying the dominated convergence theorem.
\finishproof

\rr{ \startproof{Proof of Lemma~\ref{lem:stein2}}
Fix $h \in \M_{2}$,  $x \geq 0$, and $\epsilon > 0$, and consider 
\begin{align*}
 &\frac{1}{\epsilon}\big( \partial_{x} \E F_{h}(x+ \epsilon,T)  - \partial_{x} \E F_{h}(x,T)  \big) \notag  \\
 =&\ \frac{1}{\epsilon} \E \bigg(  \E_{ x,T}  \int_{0}^{B_{0}} \big(h'(X^{(\epsilon)}(t)) - h'(X(t))\big) dt   \bigg)   +  \frac{1}{\epsilon}  \E \bigg( \E_{x,T} \int_{B_{0}}^{B_{0}^{(\epsilon)}} h'(X^{(\epsilon)}(t)) dt\bigg).
\end{align*} 
Repeating the proof of Lemma~\ref{lem:abuse:d1} yields  
\begin{align*}
    \frac{1}{\epsilon} \E \bigg(  \E_{ x,T}  \int_{0}^{B_{0}} \big(h'(X^{(\epsilon)}(t)) - h'(X(t))\big) dt   \bigg)  \to \partial_{x} \E F_{h'}(x,T) \quad \text{ as $\epsilon \to 0$}.
\end{align*}
Recall that $R_{a}(B_{0})$ is the residual interarrival time at the end of the initial busy period (which also equals the length of the first idle period $I_{0}$). If $R_{a}(B_{0}) \geq \epsilon/\delta$, then there is no arrival during the interval $[B_{0},B_{0}^{(\epsilon)})$. Since $X^{(\epsilon)}(B_{0}) = \epsilon$, this implies that 
\begin{align*}
    &\frac{1}{\epsilon}  \E \bigg( \E_{x,T} \Big( 1(R_{a}(B_{0}) \geq \epsilon/\delta) \int_{B_{0}}^{B_{0}^{(\epsilon)}} h'(X^{(\epsilon)}(t)) dt\Big)\bigg) \\
    =&\  \E \Big( \Prob_{x,T} \big(R_{a}(B_{0}) \geq \epsilon/\delta\big)\Big) \frac{1}{\epsilon}  \int_{0}^{\epsilon/\delta} h'(\epsilon - \delta t) dt.
\end{align*}
As $\epsilon \to 0$, the right-hand side converges to 
\begin{align*}
      \E \Big( \Prob_{x,T} \big(R_{a}(B_{0}) >0 \big)\Big)  \frac{1}{\delta} h'(0) = \frac{1}{\delta} h'(0).
\end{align*}
To justify the last equality, we observe that $R_{a}(B_{0}) = 0$ would imply that an arrival occurs precisely at the instant that the workload hits zero. Since the workload process is right-continuous, this would imply that $X(B_{0}) > 0$, which contradicts the definition of $B_{0}$. It remains to show that 
 \begin{align*}
       &\frac{1}{\epsilon}  \E \bigg( \E_{x,T} \Big( 1(R_{a}(B_{0}) < \epsilon/\delta) \int_{B_{0}}^{B_{0}^{(\epsilon)}} h'(X^{(\epsilon)}(t)) dt\Big)\bigg) \\
       \to&\   \frac{1}{\delta} \Big( \theta(x/\delta) + \E \big(1(T<x/\delta)\E_{x,T} \eta\big(\alpha(B_{0})\big) \big)\Big) \E \big( \partial_{x} F_{h}(\delta S, U) \big).
 \end{align*}
 Since $R_{a}(B_{0}) < \epsilon/\delta$ implies that an arrival occurs in $[B_{0},B_{0}^{(\epsilon)})$, then 
 \begin{align*}
 & \frac{1}{\epsilon}  \E \bigg( \E_{x,T} \Big( 1(R_{a}(B_{0}) < \epsilon/\delta) \int_{B_{0}}^{B_{0}^{(\epsilon)}} h'(X^{(\epsilon)}(t)) dt\Big)\bigg) \\
   =&\  \frac{1}{\epsilon}  \E \bigg( \E_{x,T} \Big( 1(R_{a}(B_{0}) < \epsilon/\delta)\int_{0}^{R_{a}(B_{0})} h'(\epsilon - \delta t) dt\Big)\bigg) \\
   &+ \frac{1}{\epsilon}  \E \bigg( \E_{x,T} \Big( 1(R_{a}(B_{0}) < \epsilon/\delta)\E\big( \partial_{x} F_{h}(\epsilon - \delta R_{a}(B_{0}) + \delta S, U)\big)\Big)\bigg).
 \end{align*}
 The first term on the right-hand side converges to zero as $\epsilon \to 0$. To analyze the second term, note that $R_{a}(B_{0}) < \epsilon/\delta$ implies that $T < x/\delta + \epsilon/\delta$, and if $x/\delta \leq T < x/\delta + \epsilon/\delta$ then $R_{a}(B_{0}) = T-x/\delta$. Therefore, the second term equals
 \begin{align*}
     &\frac{1}{\epsilon}  \E \bigg(1(T < x/\delta)  \E_{x,T} \Big( 1(R_{a}(B_{0}) < \epsilon/\delta)\E\big( \partial_{x} F_{h}(\epsilon - \delta R_{a}(B_{0}) + \delta S, U)\big)\Big)\bigg) \\
     &+ \frac{1}{\epsilon}  \E \bigg(1( x/\delta \leq T < x/\delta + \epsilon/\delta )  \E_{x,T} \Big(  \E\big( \partial_{x} F_{h}(\epsilon - (\delta T - x)  + \delta S, U)\big)\Big)\bigg).
 \end{align*}
It is straightforward to check  that $\sup_{0 \leq x' \leq \epsilon} \abs{\E\big( \partial_{x} F_{h}(x'+\delta S, U)\big) - \E\big( \partial_{x} F_{h}( \delta S, U)\big)} \to 0$ as $\epsilon \to 0$.  Lemma~\ref{lem:threebounds} and the DCT then yield
 \begin{align*}
     &\frac{1}{\epsilon}  \E \bigg(1(T < x/\delta)  \E_{x,T} \Big( 1(R_{a}(B_{0}) < \epsilon/\delta)\E\big( \partial_{x} F_{h}(\epsilon - \delta R_{a}(B_{0}) + \delta S, U)\big)\Big)\bigg) \\
     \to&\  \frac{1}{\delta}  \E \big(1(T<x/\delta)\E_{x,T} \eta\big(\alpha(B_{0})\big) \big)\E\big( \partial_{x} F_{h}( \delta S, U)\big).
 \end{align*}
 Similarly, using the fact that $\theta(x)$ is bounded, 
\begin{align*}
    &\frac{1}{\epsilon}  \E \bigg(1( x/\delta \leq T < x/\delta + \epsilon/\delta )  \E_{x,T} \Big(  \E\big( \partial_{x} F_{h}(\epsilon - (\delta T - x)  + \delta S, U)\big)\Big)\bigg) \\
    \to&\ \frac{1}{\delta} \theta(x/\delta) \E\big( \partial_{x} F_{h}( \delta S, U)\big).
\end{align*}
\finishproof}

 \subsection{Third-order bounds}

\subsubsection{The renewal process driven by idle times} 
\label{sec:lem4}
In this section we prove Lemma~\ref{lem:mixing:time}.   Recall that $\alpha(t)$ is the age of the interarrival process at time $t\geq 0$. Further recall from the discussion following \eqref{eq:BI} that   $B_n$ and $I_n$, $n \geq 0$,  are the lengths of the $n$th busy and idle periods, respectively. While $B_0$ and $I_0$ are the initial busy and idle period lengths, which depend on the initial condition $Z(0)$, the pairs $(B_n, I_n)$, $n \geq 1$,  are i.i.d.\ with distributions $\bar B$ and $\bar I$, respectively.   Given $Z(0-) = (0,0)$ (an arrival occurs to an empty system at time $t=0$), consider $\alpha(B_0)$, the age at the end of the first busy period, and let  $\bar \alpha$ be the random variable having its   distribution.

Central to our analysis is the following continuous-time Markov process that behaves like the   interarrival age during idle periods of the workload process. Formally, we let $\{\Phi(t): t \geq 0\}$ be defined by the infinitesimal generator 
    \begin{align}
G_{\Phi} f(\gamma) = \lim_{t \to 0} \frac{\E_{\gamma} f(\Phi(t)) - f(\gamma)}{t}   = f'(\gamma) + \eta(\gamma) \big(\E f(\bar \alpha)   - f(\gamma)\big), \quad \gamma \geq 0. \label{eq:gengamma}
\end{align} 
This process grows at unit rate and, at its jump times, $\Phi(t)\stackrel{d}{=}\bar\alpha$. Jumps occur with rate $\eta(\Phi(t))$.

Let $\{\Phi(t): t \geq 0\}$ and $\{ \widetilde \Phi(t): t \geq 0\}$ be two copies (possibly coupled) of the Markov process defined by \eqref{eq:gengamma}. The following technical lemma is proved in Appendix~\ref{app:renewal:idle:phi}. The conditioning on the right-hand side of \eqref{eq:phi:equation} is non standard and means: pick a single $\bar \alpha$-distributed random variable $\xi$ and run $\{\Phi(t), t \geq 0\}$ from time zero starting at $\xi$, and run $\{\widetilde \Phi(t), t \geq 0\}$ from time $s$ starting at $\xi$.
\begin{lemma}
\label{lem:renewal:idle:phi}
Let $v = v(x) = x/\delta$. For any $(x,r_{a}) \in \mathbb{S}$ with $r_{a} < v$ and any $s \geq 0$, 
\begin{align}
     &\E_{x+\delta s,r_{a}} \eta (\alpha(B_{0})) - \E_{x,r_{a}} \eta(\alpha(B_{0})) \notag \\
     =&\   \E\big( \eta (\Phi(v-r_{a}+s)) -   \eta (\widetilde \Phi(v-r_{a}+s)) |  \Phi(0) = \widetilde \Phi(s),  \Phi(0)\stackrel{d}{=} \bar \alpha  \big), \label{eq:phi:equation}
\end{align}
Furthermore if $\Phi(0) \stackrel{d}{=} \bar \alpha$, then the  inter-jump times of $\{\Phi(t): t \geq 0\}$ are i.i.d.\ $\bar I$. 
\end{lemma}
The final ingredient for the proof of Lemma~\ref{lem:mixing:time} is a  coupling $\{(\widetilde \Phi(t), \Phi(t)) :  t \geq 0\}$, which allows us to bound the expression in \eqref{eq:phi:equation}. Defining
\begin{align*}
\eta_m(x,y) = \min\{\eta(x),\eta(y)\}, \quad \eta_{\Delta}(x,y) =  \max\{\eta(x),\eta(y)\} -  \min\{\eta(x),\eta(y)\},
\end{align*}
we let $\{(\widetilde \Phi(t), \Phi(t)),\ t \geq 0\}$ have the same distribution as the Markov process  defined by the generator 
\begin{align}
G_{J} f(\tilde \gamma, \gamma)  =&\ \partial_{\tilde \gamma} f(\tilde \gamma, \gamma) + \partial_{\gamma} f(\tilde \gamma, \gamma) \notag \\
&+ \eta_{m}(\tilde \gamma, \gamma) \big( \E f(\bar \alpha, \bar \alpha) - f(\tilde \gamma, \gamma) \big)   \notag \\
&+  \eta_{\Delta}(\tilde \gamma, \gamma)     1(\eta(\gamma) < \eta(\tilde \gamma)) \big( \E f(\bar \alpha, \gamma)- f(\tilde \gamma, \gamma)\big)  \notag \\
&+ \eta_{\Delta}(\tilde \gamma, \gamma)  1(\eta(\gamma) > \eta(\tilde \gamma)) \big(\E f(\tilde \gamma, \bar \alpha) - f(\tilde \gamma, \gamma)\big). \label{eq:joint:Markov}
\end{align} 
Note that the marginal law of either component of this process is equivalent to the Markov process defined in \eqref{eq:gengamma}. Furthermore, when this process jumps it either couples (with rate $\eta_{m}(\tilde \gamma, \gamma)$), or only one of the components jumps.  
\startproof{Proof of Lemma~\ref{lem:mixing:time}}
Fix $v = x/\delta \geq 0$ and $r_{a} < v$, and  let $ \tau_{C} = \inf\{t \geq s : \widetilde \Phi(t) = \Phi(t)\}$. It follows that 
\begin{align}
&  \abs{\E_{x+\delta s,r_{a}} \eta (\alpha(B_{0})) - \E_{x,r_{a}} \eta(\alpha(B_{0}))} \notag  \\
    =&\ \abs{ \E\big( \eta (\Phi(v-r_{a}+s)) -   \eta (\widetilde \Phi(v-r_{a}+s)) |  \Phi(0) = \widetilde \Phi(s),  \Phi(0)\stackrel{d}{=} \bar \alpha  \big)} \notag \\
    \leq&\ \overline{\eta}   \Prob\big(\tau_{C} > v-r_{a}+s \big|  \Phi(0) = \widetilde \Phi(s),  \Phi(0)\stackrel{d}{=} \bar \alpha    \big). \label{eq:prob:couple}
\end{align}
If  $\widetilde \Phi(s) > \Phi(s)$, then $\eta_{m}(\Phi(t),\widetilde \Phi(t)) = \eta(\widetilde \Phi(t))$ for all $t \geq s$  until the first jump of $\{\widetilde \Phi(t): t \geq 0\}$, at which  point coupling occurs.  Since    $\widetilde \Phi(s)  \stackrel{d}{=} \bar \alpha $, we know by Lemma~\ref{lem:renewal:idle:phi} that the first jump after $s$ happens after $\bar I$ amount of time, implying that the probability of no jump on $(s,v-r_{a}+s]$ is at most
\begin{align}
   \Prob(\bar I > v-r_{a}). \label{eq:bari1}
\end{align}
Alternatively, if $\widetilde \Phi(s) < \Phi(s)$, then $\tau_{C}$ corresponds to the first jump time after $s$ of $\{ \Phi(t): t \geq 0\}$. Unlike $\widetilde \Phi(s)$, it is not true that $\Phi(s)$ is distributed like $\bar \alpha$, because the latter process has been running since time zero with $\Phi(0)\stackrel{d}{=} \bar \alpha$. We will shortly prove that the probability that $\{ \Phi(t): t \geq 0\}$ does not jump on $(s,v-r_{a}+s]$ is at most  
 \begin{align}
     \Prob(\bar I > v-r_{a}) \E  (N_{\Phi}(s) | \Phi(0) \stackrel{d}{=} \bar \alpha) , \label{eq:ibar:to:prove}
 \end{align}
where  $ N_{\Phi}(t)$ is the number of jumps made by $\{ \Phi(t): t \geq 0\}$ on $[0,t]$.  Adding \eqref{eq:bari1} and \eqref{eq:ibar:to:prove} yields 
\begin{align*}
  \overline{\eta}   \Prob\big(\tau_{C} > v-r_{a}+s \big|  \Phi(0) = \widetilde \Phi(s) \stackrel{d}{=} \bar \alpha    \big) \leq \overline{\eta} (1 + \E (N_{\Phi}(s) | \Phi(0) \stackrel{d}{=} \bar \alpha))  \Prob(\bar I > v-r_{a}).
\end{align*}
Let $\{ N_{U}(t) : t \geq 0 \}$ be the counting process of a zero-delayed renewal process whose inter-event times have the interarrival distribution $U$.  It follows that 
\begin{align*}
    \E (N_{\Phi}(s) | \Phi(0) \stackrel{d}{=} \bar \alpha) \leq \E  N_{U}(s) \leq \lambda s + \lambda^2 \E U^2, \quad s \geq 0,
\end{align*}
where the first inequality is due to the  hazard rate $\eta(x)$ being nonincreasing, and the second inequality follows from Lorden's inequality \citep{Lord1970}.

It remains to  verify \eqref{eq:ibar:to:prove}. Let $T_0 = 0$ and $T_n$, $n \geq 1$, be the $n$th jump time of $\{\Phi(t): t \geq 0\}$. Also let $J_n = T_{n}-T_{n-1}$, $n \geq 1$, be the duration between the $n$th and $(n-1)$st jump. We know by Lemma~\ref{lem:renewal:idle:phi} that $J_{n} \stackrel{d}{=} \bar I$ whenever   
$\Phi(0) \stackrel{d}{=} \bar \alpha$. It follows that  (to simplify notation, all probabilities and expectations that follow are conditional on $\Phi(0) \stackrel{d}{=} \bar \alpha$) 
 \begin{align*}
 \Prob(N_{\Phi}(v-r_{a}+s) - N_{\Phi}(s) = 0)  =&\ \Prob( J_{N_{\Phi}(s)+1} > v-r_{a}+s - T_{N_{\Phi}(s)} ) \\
=&\ \sum_{n=0}^{\infty} \Prob(J_{n+1} > v-r_{a}+s - T_{n}, N_{\Phi}(s) = n ).
 \end{align*}
Since $\{N_{\Phi}(s)=n\}  = \{T_n < s, J_{n+1}>s-T_{n} \}$, the right-hand side equals 
\begin{align*}
 \sum_{n=0}^{\infty} \Prob(J_{n+1} > v-r_{a}+s - T_{n},T_n < s)   \leq&\  \sum_{n=0}^{\infty} \Prob(J_{n+1} > v-r_{a},T_n < s ) \\
=&\ \Prob(\bar I > v-r_{a})  \sum_{n=0}^{\infty} \Prob(T_n < s ) \\
=&\ \Prob(\bar I > v-r_{a})  \sum_{n=0}^{\infty} \Prob(N_{\Phi}(s) \geq n  ) \\
=&\ \Prob(\bar I > v-r_{a})  \E N_{\Phi}(s),
\end{align*}
where the first equality is true because $J_{n+1} \stackrel{d}{=}\bar I$  and   $J_{n+1}$ is independent of $T_{n}$.
\finishproof
\begin{remark}
\label{rem:general:eta}
Another family of hazard rates for which the coupling probability can be controlled 
consists of those uniformly bounded away from zero: $\eta(x) \geq \uline \eta > 0$. 
In this case, $\eta_{m}(\tilde \gamma, \gamma) \geq \uline \eta$, so coupling occurs after 
an exponentially distributed time with rate $\uline \eta$, leading to the bound 
\begin{align*}
       \Prob\big(\tau_{C} > v-r_{a}+s \big| \Phi(0) = \widetilde \Phi(s), \ \Phi(0)\stackrel{d}{=} \bar \alpha\big) 
   \leq \overline{\eta}\, e^{-\uline \eta (v-r_{a})}.
\end{align*}
\end{remark}

\subsubsection{Proof of Lemma~\ref{lem:d3:form}.}
\label{app:d3:form:proof}
\begin{lemma} \label{lem:d12}
Suppose that $U$ has a bounded density. Then for any \rr{$h \in \M_{2}$ and any } $x \geq 0$, 
\begin{align*}
    \partial_{x}^{2} \E \big(F_h(x + \delta S,U) - F_h(x,U) \big) = \E^{S}\big(\partial_{x}^{2} \E^{U}  F_h(x + \delta S,U) - \partial_{x}^{2} \E^{U}  F_h(x,U)\big).  
\end{align*}
\end{lemma} 
\startproof{Proof of  Lemma~\ref{lem:d3:form}}
Observe that 
\begin{align*}
    \delta \bar F_{h}'''(x) =&\ \lambda \partial_{x}^{2} \E \big(F_h(x + \delta S,U) - F_h(x,U) \big) \\
    =&\  \lambda \E^{S}\big(\partial_{x}^{2} \E^{U}  F_h(x + \delta S,U) - \partial_{x}^{2} \E^{U}  F_h(x,U)\big),
\end{align*}
where the first equality is due to  \eqref{eq:d3:initial} and the second is due to Lemma~\ref{lem:d12}. Applying the expression for $\partial_{x}^{2} \E^{U}  F_h(\cdot,U)$ from Lemma~\ref{lem:stein2} to the right-hand side yields the result.
\finishproof

\startproof{Proof of Lemma~\ref{lem:d12}}
Though we do not assume $F_{h}(z)$ to be well defined, note that 
\begin{align*}
    &\partial_{x} \E \big(F_h(x + \delta S,U) - F_h(x,U) \big) \\
    =&\ \lim_{\epsilon \to 0}  \bigg(\frac{1}{\epsilon} \E \big(F_h(x + \epsilon + \delta S,U) - F_h(x+ \delta S,U) \big) - \frac{1}{\epsilon} \E \big(F_h(x + \epsilon,U) - F_h(x,U) \big) \bigg) \\
    =&\  \partial_{x}  \E  F_h(x + \delta S,U) - \partial_{x}\E F_h(x,U),
\end{align*}
so that by differentiating both sides with respect to $x$, we arrive at 
\begin{align*}
    \partial_{x}^{2} \E \big(F_h(x + \delta S,U) - F_h(x,U) \big) =  \partial_{x}^{2} \E  F_h(x + \delta S,U) - \partial_{x}^{2} \E F_h(x,U).
\end{align*}
By repeating the arguments used to prove Lemma~\ref{lem:abuse:d1}, one can check that 
\begin{align*}
    \partial_{x}^{2} \E  F_h(x + \delta S,U) = \partial_{x} \E^{S}\partial_{x} \E^{U}  F_h(x + \delta S,U).
\end{align*}
Similarly, 
\begin{align*}
    \partial_{x} \E^{S}\partial_{x} \E^{U}  F_h(x + \delta S,U) =   \E^{S}\partial_{x}^2 \E^{U}  F_h(x + \delta S,U)
\end{align*}
follows   from repeating the proof of Lemma~\ref{lem:stein2}.
\finishproof

\subsubsection{Proof of Lemma~\ref{lem:renewal:idle:phi}.}
\label{app:renewal:idle:phi}
\startproof{} 
Define $\ell_{n}$, $n \geq 0$, to be the end of the $n$th busy period. Namely, $\ell_{0} = B_{0}$ and
\begin{align*}
      \ell_{n} = \ell_{n-1} + I_{n-1} + B_{n}, \quad n \geq 1.
\end{align*} 
Letting $T_{0} = 0$ and $T_{n+1} = T_{n} + I_{n}$, $n \geq 0$, so that the $n$th idle period occurs on the interval $[\ell_{n},\ell_{n} + T_{n+1}-T_{n} )$, we   define 
\begin{align}
    \Gamma(t) = \alpha(\ell_{n} + t-T_{n}), \quad  t \in [T_{n},T_{n+1}),\ n \geq 0,
   \label{eq:gamdef}
\end{align}
to be the age of the interarrival process at the instant when the server has idled for  exactly $t$ time units.

Recalling our synchronous coupling $\{Z^{(x)}(t): t \geq 0\}$ introduced in \eqref{eq:couplinggap} and the corresponding initial busy period    $B_{0}^{(x)}$,  it follows that  
\begin{align*}
      \alpha(B_{0}^{(x)}) = \Gamma(v), \quad x \geq 0,
\end{align*}
because $B_{0}^{(x)}$ is precisely the instant when    $\{Z(t): t \geq 0\}$ idles for   $v$ time units. Therefore,   for any $(x,r_{a}) \in \mathbb{S}$ with $r_{a} < v$ and any $s \geq 0$, 
\begin{align}
\E_{x+\delta s,r_{a}} \eta (\alpha(B_{0})) - \E_{x,r_{a}} \eta(\alpha(B_{0})) =&\ \E_{0,r_{a}} \eta (\alpha(B_{0}^{(x+\delta s)})) - \E_{0,r_{a}} \eta (\alpha(B_{0}^{(x)}))  \notag \\
=&\ \E_{0,r_{a}}\big( \eta (\Gamma(v+s)) -   \eta (\Gamma(v))  \big),   \label{eq:coming:back}
\end{align}
where we recall that  $\E_{x,r_{a}}(\cdot)$ is   the expectation conditioned on $Z(0) = (x,r_{a})$. 
To proceed, we now show that although $\{\Gamma(t) : t \geq 0\}$ is defined via \eqref{eq:gamdef}, it is also a Markov process with generator \eqref{eq:gengamma}.

 We know from   \eqref{eq:gamdef} that   $\{\Gamma(t): t \geq 0\}$  increases at a unit rate. It  jumps at times $T_{n+1}$, $n \geq 0$, and its distribution at jump times satisfies  $\Gamma(T_{n+1}) = \alpha(\ell_{n+1}) \stackrel{d}{=} \bar \alpha$ for all $n \geq 0$, because $\ell_{n+1}$ marks the end of a busy period initiated by an arrival to an empty system. Lastly,   conditioned on $\Gamma(t)$,  the probability that a jump occurs on the interval $(t,t+dt)$ equals
\begin{align*}
    \frac{\Prob(\Gamma(t) < U < \Gamma(t)+dt)}{\Prob(U > \Gamma(t))} =  \eta(\Gamma(t))dt  + o(dt),
\end{align*}
where $o(dt) \to 0$ as $dt \to 0$;  the Markov property is simple to verify. Thus  $\{\Gamma(t): t \geq 0\}$  is a Markov process described by the generator \eqref{eq:gengamma} and, given the same initial condition,  both $\{\Phi(t): t \geq 0\}$ and $\{\Gamma(t): t \geq 0\}$ are equal in distribution. We now argue  that given $\Phi(0) \sim \bar \alpha$, the inter-jump times of $\{\Phi(t): t \geq 0\}$ are i.i.d.\ $\bar I$.

By construction, the value of $\Gamma(\cdot)$ at jump times $T_{n+1}$ satisfies $\Gamma(T_{n+1}) \stackrel{d}{=} \bar \alpha$, and the inter-jump durations $T_{n+2}-T_{n+1} \stackrel{d}{=} \bar I$  for $n \geq 0$. Since $\{\Gamma(t): t \geq 0\}$ is a Markov process, it follows that conditioned on $\Gamma(0) \stackrel{d}{=} \bar \alpha$, the time until the first jump has the same distribution as $\bar I$.

Coming back to right-hand side of \eqref{eq:coming:back}, suppose that   $Z(0)=(0,r_{a})$. Then $\Gamma(r_{a}) = \alpha(\ell_{1})\stackrel{d}{=}\bar \alpha$, because $B_0=0$ and $I_0=r_{a}$, which in turn yields $r_{a} \in [T_{1},T_{2})$ (since $T_1=r_{a}$). Furthermore, conditioned on both $Z(0)=(0,r_{a})$ and   $\Gamma(r_{a})$, the value of $\Gamma(r_{a}+t)$, $t \geq 0$, is independent of $Z(0)$.  It follows that  for any  $v$ and $r_{a} < v$,
\begin{align*}
       \E_{0,r_{a}}\big( \eta (\Gamma(v+s)) -   \eta (\Gamma(v))  \big)  
    =&\   \E\big( \eta (\Gamma(v+s)) -   \eta (\Gamma(v)) | \Gamma(r_{a})  \stackrel{d}{=} \bar \alpha  \big) \\
    =&\  \E\big( \eta (\Phi(v+s)) -   \eta (\Phi(v)) | \Phi(r_{a}) \stackrel{d}{=} \bar \alpha   \big)\\
    =&\  \E\big( \eta (\Phi(v-r_{a}+s)) -   \eta (\Phi(v-r_{a})) | \Phi(0) \stackrel{d}{=} \bar \alpha   \big),
\end{align*} 
where the third equality comes from the Markov property. Finally,  if  $\{ \widetilde \Phi(t) : t \geq 0\}$ is a copy of the Markov process defined by \eqref{eq:gengamma}, then 
\begin{align*}
    & \E\big( \eta (\Phi(v-r_{a}+s)) -   \eta (\Phi(v-r_{a})) | \Phi(0) \stackrel{d}{=} \bar \alpha   \big) \\
    =&\   \E\big( \eta (\Phi(v-r_{a}+s)) -   \eta (\widetilde \Phi(v-r_{a})) |   \Phi(0) = \widetilde \Phi(0),  \Phi(0)\stackrel{d}{=} \bar \alpha \big)\\
    =&\   \E\big( \eta (\Phi(v-r_{a}+s)) -   \eta (\widetilde \Phi(v-r_{a}+s)) |   \Phi(0) = \widetilde \Phi(s),  \Phi(0)\stackrel{d}{=} \bar \alpha \big).
\end{align*}
\finishproof

\subsection{Proof of Lemma~\ref{lem:loulou}}
\label{sec:busy:period:bound}
\startproof{Proof of Lemma~\ref{lem:loulou}}
We first prove \eqref{eq:loulou}.  Let $\hat h(x) = x$, recall that $\bar F_{\hat h}(0) = 0$ due to \eqref{eq:d1}, and consider the Poisson equation \eqref{eq:the:poisson:abuse} with $h(x) = \hat h(x)$ evaluated at $x = 0$, which results in 
\begin{align*}
 \lambda \E \big( F_{\hat h}(\delta S, U)  - F_{\hat h}(0,U)\big) =  \delta \E V,
\end{align*}
Using the synchronous coupling  $\{Z^{(\epsilon)}(t)\}$ introduced in Section~\ref{sec:mhorizon}, it follows that
\begin{align}
\delta \E V =&\ \lambda \E \Big(\int_{0}^{B^{(\delta S)}} \E_{0,U} \big( X^{(\delta S)}(t) - X(t) \big) dt \Big). \label{eq:xrel}
\end{align}
From \eqref{eq:couplinggap}  we know that  the difference $X^{(\delta S)}(t) - X(t)$ decays at rate $\delta$ only during the idle periods of $\{X(t) :  t \geq 0\}$.
Recall that the idle and busy period durations  are $I_{0}, I_{1}, \ldots$, and $B_{0}, B_{1}, B_{2}, \ldots$, respectively, and that $B_{0} = 0$  because $X(0) = 0$. It follows that for all times $t$ corresponding to the busy period  $B_{k}$, $k \geq 1$, 
\begin{align*}
X^{(\delta S)}(t) - X(t) = \delta \Big( S -  \sum_{i=0}^{k-1} I_{i} \Big)^{+}.
\end{align*}
Letting $\mathcal{I} = \{t \in [0,B^{(\delta S)}]: X(t) = 0\}$,  it follows that  
\begin{align*}
  \int_{0}^{B^{(\delta S)}} \big( X^{(\delta S)}(t) - X(t) \big) dt  =&\   \int_{\mathcal{I} } \big( X^{(\delta S)}(t) - X(t) \big) dt +   \int_{[0,B^{(\delta S)}] \setminus \mathcal{I}} \big( X^{(\delta S)}(t) - X(t) \big) dt\\
=&\  \delta S^2/2 +     \sum_{k=1}^{\infty}  B_{k} \delta \Big( S -   \sum_{i=0}^{k-1} I_{i} \Big)^{+}.
\end{align*}
We conclude \eqref{eq:loulou} by combining this equation  with \eqref{eq:xrel}, noting that $B_k $ is independent of $( S - I_{0} - I_1 \ldots - I_{k-1})^{+}$ and $B_{k}\stackrel{d}{=} \bar B$,  and that  $I_{0} \stackrel{d}{=} U$ since $Z(0) = (0,U)$. The equality in \eqref{eq:barBbound} is true because  $\E V =  \lambda \E S^2/2 + \rho \E W $  due to Corollary~X.3.5 of \cite{Asmu2003}. The inequality follows from the well-known bound in (7') of \cite{King1962}, which says that
\begin{align*}
     \E W \leq \frac{\text{Var}(S-U)}{2\E(S-U)} = \frac{\rho \text{Var}(S-U)}{2(1-\rho)\E S }.
\end{align*}
 \finishproof

\subsection{Proof of Lemma~\ref{lem:final:aux}}
\label{app:final:aux:proof}
\startproof{}
For convenience, let $\nu = 1/\E Y_2$. We  use the following inequalities throughout the proof, which follow from the fact that $Y_2$ has density $\nu e^{-\nu y}$ and is independent of $R_{a}$. 
\begin{align}
& \Prob(\delta R_{a} \geq Y_2 ) = \E (1-e^{-\nu\delta R_{a}}) \leq \nu \delta \E R_{a},  \label{eq:mom:0} \\
    &\E Y_2 1(\delta R_{a} \geq Y_2) = \frac{1}{\nu} \E  \big( (1-e^{-\nu\delta R_{a}}) - \nu\delta R_{a} e^{-\nu\delta R_{a}}  \big) \leq 2 \nu \delta^2 \E  R_{a}^2  ,\label{eq:mom:1} \\
    &\E Y_2^2 1(\delta R_{a} \geq Y_2) = \frac{1}{\nu^2} \E  \big( 2(1-e^{-\nu\delta R_{a}}) - 2\nu\delta R_{a} e^{-\nu\delta R_{a}} - (\nu\delta R_{a})^2 e^{-\nu\delta R_{a}} \big) \leq 3 \nu \delta^3 \E  R_{a}^3.  \label{eq:mom:2}
\end{align}
We also use the facts that $U,S,S',R_{a},\bar I$, and $Y_2$ are independent, and that $S$ and $S'$ have the same distribution. We first prove \eqref{eq:ineq:1}--\eqref{eq:ineq:3} and  then \eqref{eq:aineq:1}--\eqref{eq:aineq:2}.

\subsubsection*{Proof of \eqref{eq:ineq:1}--\eqref{eq:ineq:3}.}
Inequality \eqref{eq:ineq:1} follows from $J(x,r_{a})  = - (x \wedge \delta r_{a}) + \delta S'$. Next, we prove \eqref{eq:ineq:2}. We first note that the definition of $\epsilon(x,r_{a},s)$  in Lemma~\ref{lem:relate} and the fact that $h(x) = x$ both imply that 
\begin{align*}
    \int_{0}^{r_{a}} \abs{ h\big((x + \delta s-\delta t)^{+}\big) -  h\big((x  -\delta t)^{+}\big) } dt \leq r_{a} \delta s.
\end{align*}
Furthermore, 
\begin{align*}
& 1(\delta r_{a} > x) \E^{S'} \Big(  \int_{-x}^{-(x+\delta s)  } \partial_{x}   \E^{U} F_h(x+\delta s +v +\delta S',U) dv\Big) \\
\leq&\ 1(\delta r_{a} > x)  \delta s \E^{S'} \Big( \sup_{0 \leq w \leq \delta s}\partial_{x}   \E^{U} F_h(w +\delta S',U) \Big) \\ 
\leq&\ 1(\delta r_{a} > x)   s (\delta s + \delta \E S')( 1 + 2 \overline{\eta}  \E \bar B ),
\end{align*}
where in the second inequality we used \eqref{eq:1bounds} of Lemma~\ref{lem:2bounds:stein:factors}. Combining the bounds and using \eqref{eq:mom:0} yields
\begin{align*}
    \E \abs{\epsilon(Y_2,R_{a},S)} \leq&\ \delta \E R_{a} \E S + \Prob(\delta R_{a} \geq Y_2) \delta (\E S^2 + (\E S)^2)( 1 + 2 \overline{\eta}  \E \bar B ) \\
    \leq&\ \delta \E R_{a} \E S + \delta^2 \nu \E R_{a}(\E S^2 + (\E S)^2)( 1 + 2 \overline{\eta}  \E \bar B ).
\end{align*}
Next we prove \eqref{eq:ineq:3}. Recall from Lemma~\ref{lem:2bounds:stein:factors} that
\begin{align*}
        \abs{\delta \bar F_{h}'(x)} \leq&\   x( 1 + (\lambda + \overline{\eta}) \E \bar B ),\\ 
        \abs{\delta \bar F_{h}''(x)} \leq&\   (1+x)\big( 1 + (\lambda + \overline{\eta}) \E \bar B  \big), \quad x \geq 0. 
    \end{align*}
    Therefore, 
\begin{align*}
     &\E \abs{ 1(\delta R_{a} \geq Y_2) \big(  \delta\bar F_{h}'(\delta S )  -  \delta\bar F_h'(Y_2) - \delta(S-R_{a}) \delta\bar F_h''(Y_2)\big)} \\
    \leq&\   \Big(\E \big(1(\delta R_{a} \geq Y_2)( \delta S + Y_2)\big) + \E \big(1(\delta R_{a} \geq Y_2)\delta(S-R_{a})(1+ Y_2) \big)\Big) \big( 1 + (\lambda + \overline{\eta}) \E \bar B  \big) \\
    \leq&\ \delta^2 \big(  \nu \big( 2\E S \E R_{a} + 5 \E R_{a}^2  \big) +   \E  R_{a}^2 \big) \big( 1 + (\lambda + \overline{\eta}) \E \bar B  \big).
\end{align*} 
The last inequality follows from using  \eqref{eq:mom:0} and \eqref{eq:mom:1} to show that 
\begin{align*}
    &\E \big(1(\delta R_{a} \geq Y_2)( \delta S + Y_2)\big) \leq \delta^2  \nu \E S \E R_{a} +   2 \nu \delta^2  \E R_{a}^2, \\
    & \delta \E S\E \big(1(\delta R_{a} \geq Y_2) (1+ Y_2) \big) \leq \delta \E S  \big(\nu \delta \E R_{a} + 2 \nu \delta^2 \E R_{a}^2\big), \\
    & \delta \E \big(1(\delta R_{a} \geq Y_2)R_{a}(1+ Y_2) \big) \leq \delta   \big(\nu \delta  \E R_{a}^2 + \delta  \E R_{a}^2\big), 
\end{align*} 
 where in the final inequality we used the fact that $ \E \big(1(\delta R_{a} \geq Y_2)R_{a} Y_2 \big) \leq \delta \E R_{a}^2$ instead of \eqref{eq:mom:1}. Using the latter would have resulted in a term involving $\E R_{a}^3$. 

\subsubsection*{Proof of \eqref{eq:aineq:1}--\eqref{eq:aineq:2}.} 
Recall that 
\begin{align*}
    C =    (3   + (\E U + \E \bar I) \lambda  (1 + \rho + \lambda^2 \E U^2))   \nu  \overline{\eta} \delta \E \bar B.
\end{align*}
We  claim that 
\begin{align}
    &\E  \big|\delta^2 \bar F_{ h}'''(Y_2+u) \big|   \leq      \overline{C},   &u > 0, \label{eq:upos}   \\
    & \E^{Y_2} \Big(1(\delta R_{a}  \leq Y_2)  \big|\delta^2 \bar F_{ h}'''(Y_2+u) \big| \Big)   \leq      \overline{C},  &u \in [-\delta R_{a},0]. \label{eq:uneg}
\end{align}
Then  \eqref{eq:aineq:1} follows by applying the bound in \eqref{eq:upos} to 
\begin{align*}
    \lambda    \E \Big(\int_{0}^{\delta S} (\delta S - v)\int_{0}^{v}  \abs{\bar F_{ h}'''(Y_2+u) }     du dv \Big)  =&\ \lambda   \E^{S} \Big(\int_{0}^{\delta S} (\delta S - v)\int_{0}^{v}  \E^{Y_2} \abs{\bar F_{ h}'''(Y_2+u) }     du dv\Big).
\end{align*}
To prove \eqref{eq:aineq:2}, observe that both \eqref{eq:upos} and \eqref{eq:uneg} imply that 
\begin{align*}
    &\delta \frac{1}{\delta^2} \E  \Big(  \int_{0}^{\delta(S -R_{a} )} \int_{0}^{v}   \big| 1(\delta R_{a}  < Y_2) \delta^2\bar F_h'''(Y_2+u) \big| du  dv \Big)  \\
    =&\ \delta \frac{1}{\delta^2} \E^{S,R_{a}}  \Big(  \int_{0}^{\delta(S -R_{a} )} \int_{0}^{v}   \E^{Y_2} \big| 1(\delta R_{a}  < Y_2) \delta^2\bar F_h'''(Y_2+u) \big|  du  dv \Big) \\
    \leq&\ \delta \E (S-R_{a})^2  \overline{C}.
\end{align*}
It remains to prove \eqref{eq:upos} and \eqref{eq:uneg}. 
Recall  from  Lemma~\ref{lem:3bounds:stein:factors} that   for any $x,u$ such that $x + u \geq 0$,   
\begin{align}
      \big|\delta^2 \bar F_{ h}'''(x+u) \big|   \leq&\  \big( \Prob(\delta U > x+u)  3\lambda    +  \Prob(\delta U < x +u  < \delta\bar I +\delta U) \lambda   (1 + \rho + \lambda^2 \E U^2)\big)  \overline{\eta} \E \bar B. \label{eq:recall:3bound}
\end{align}  
Then \eqref{eq:upos} follows once we observe that for   any $u > 0$,  
\begin{align}
     & \Prob(\delta U > Y_2+u) \leq \Prob(\delta U > Y_2)   = \E (1-e^{-\nu \delta U})   \leq \nu \delta \E U  ,\notag  \\
     &\Prob(\delta U < Y_2+u < \delta\bar I +\delta U)\leq \nu \delta (\E U + \E \bar I), \label{eq:aset:1}
\end{align} 
where the  last inequality is true because 
\begin{align*}
    &\Prob(\delta U < Y_2+u < \delta\bar I +\delta U) \\
    =&\ \Prob( Y_2 \leq \delta U <  Y_2+u < \delta\bar I +\delta U) + \Prob( \delta U < Y_2,  Y_2+u < \delta\bar I +\delta U) \\
    \leq&\ \Prob(\delta U \geq Y_2) + \Prob(\delta U < Y_2 < \delta U + \delta \bar I) \\
    =&\ \E (1-e^{-\nu \delta U})  + \E e^{-\nu \delta U} (1-e^{-\nu \delta \bar I})  \leq \nu \delta (\E U + \E \bar I).
\end{align*}
Similarly, \eqref{eq:uneg} follows from the fact that for  any $u \in [-\delta R_{a} , 0]$,  
\begin{align}
     & \E^{Y_2,U} \big( 1(\delta R_{a} < Y_2)1(\delta U > Y_2 + u  )   \big) \leq  \nu \delta \E U, \notag \\
     &\E^{Y_2,U,\bar I} \big(1(\delta U < Y_2+u < \delta\bar I +\delta U)   \big)  \leq \nu \delta ( \E U+ \E \bar I ) \E \bar B. \label{eq:aset:2}
\end{align}
The first inequality in \eqref{eq:aset:2} follows from
\begin{align*}
     &\E^{Y_2,U} \big( 1(\delta R_{a} < Y_2)1(\delta U > Y_2 + u  )   \big)  \\
     \leq&\  \E^{Y_2,U} \big( 1(\delta R_{a} < Y_2 < \delta U +\delta R_{a} )   \big) = e^{-\nu\delta R_{a}} \E (1-e^{-\nu \delta U}) \leq \nu \delta \E U
\end{align*} 
and the second from 
\begin{align*}
E^{Y_2,U,\bar I} \big(1(\delta U < Y_2+u < \delta\bar I +\delta U)   \big)  \leq&\ E^{Y_2,U,\bar I} \big(1(\delta U -u < Y_2 < \delta\bar I +\delta U-u)   \big) \\
    =&\ \E \big( e^{-\nu (\delta U -u)} (1-e^{-\nu\delta \bar I})\big) \leq \nu \delta \E \bar I.
\end{align*}  
\finishproof

\section{The $G/M/\infty$ system: supporting proofs}
\label{app:gminf}
 To prove Proposition~\ref{prop:extraction:gminf} we require the  following inversion formula, which  is proved exactly like Lemma~\ref{lem:palm:inversion:2}.
\begin{lemma}
    \label{lem:palm:inversion:mginf}
   Initialize $Z(0) \sim Z$ and fix  $f: \R \to \R$ with
    \begin{align*}
     \E\abs{f(X)}< \infty \quad \text{ and  } \quad \E \abs{\int_{0}^{R_{a}(0)} f(X(t)) dt} < \infty,
 \end{align*}
 which holds, in particular, when $f(x)$ is bounded.   Then 
    \begin{align}
        \E f(X) =&\ \E   \int_{0}^{1} \int_{0}^{U(t)} f(X(t+u))du  d A(t)  \label{eq:palm:inversion:arrive:mginf}
    \end{align}
\end{lemma}
We now have what we need to prove  Proposition~\ref{prop:extraction:gminf}. 

\startproof{Proof of Proposition~\ref{prop:extraction:gminf}}
Repeating the arguments used in the proof of Proposition~\ref{prop:extraction:workload} yields
\begin{align*}
    \E \int_{0}^{1} \Delta  f(\TX(t-))  d A(t) =&\ \frac{1}{2}\delta^2  \lambda c_{U}^{2} \E f''(X) + \epsilon_{A}(f) = \frac{1}{2} \mu c_{U}^{2} \E f''(X) + \epsilon_{A}(f),
\end{align*}
where
\begin{align} 
        \epsilon_{A}(f) =&\ \frac{1}{6} \delta^3 \E \int_{0}^{1}(1 - \lambda U(t))^3 f'''(\xi(t))d A(t) \notag  \\
    &- \frac{1}{2} \delta^2 \lambda c_{U}^{2} \E \int_{0}^{1} \int_{0}^{U(t)} \big( X(t+u)-X(t-)\big) f'''(\xi(t+u))du d A(t).\label{eq:eA:verbose:gmfin}
     \end{align}
To bound the first term on the right-hand side, we use the fact that  $\E A(1) = \lambda$ and $\delta^2 \lambda = \mu$ to get  
\begin{align*}
    \frac{1}{6} \delta^3 \E \int_{0}^{1}\abs{(1 - \lambda U(t))^3 f'''(\xi(t))}d A(t)  \leq \frac{1}{6} \delta  \mu \norm{f'''} \E \abs{1-\lambda U}^{3} 
\end{align*}
To bound the second term in \eqref{eq:eA:verbose:gmfin},  let $D[t,t+u]$ be the number of departures on $[t,t+u]$ and recall that $X(t) = \delta(Q(t)-R)$. It follows that  
\begin{align*}
    & \E \int_{0}^{1} \int_{0}^{U(t)} \abs{\big( X(t+u)-X(t-)\big) f'''(\xi(t+u))} du d A(t) \\
    \leq&\ \norm{f'''} \delta \E \int_{0}^{1} \int_{0}^{U(t)} \abs{Q(t+u)-Q(t-)} du dA(t)\\
    \leq&\ \norm{f'''} \delta \E \int_{0}^{1} \int_{0}^{U(t)} (1 + D[t,t+u]) du  dA(t)\\
    \leq&\ \norm{f'''} \delta \E \int_{0}^{1} \int_{0}^{U(t)} (1 + \mu Q(t) u) du  dA(t)\\
    =&\ \norm{f'''} \delta \E \int_{0}^{1}  (U(t) + \mu Q(t) U^2(t)/2)   dA(t)\\
    =&\ \norm{f'''} \delta \Big(1 + \frac{1}{2} \E U^2 \mu  \E \int_{0}^{1}   Q(t)     dA(t)\Big)\\
    =&\ \norm{f'''} \delta \Big(1 + \frac{1}{2} \E U^2 \mu  \E \int_{0}^{1}  (1+Q(t-))     dA(t)\Big).
\end{align*}
The third inequality is true because  the number of departures on $[t,t+u]$ is bounded by a Poisson process with rate $\mu Q(t)$. Using $\E A(1) = \lambda$ and \eqref{eq:qt-:gminf} of Lemma~\ref{lem:relationships:gminf}, the right-hand side equals 
\begin{align*}
     & \norm{f'''} \delta \Big(1 + \frac{1}{2} \E U^2 \mu (\lambda +  \mu \E (Q-R)^2 + \lambda R -\lambda)\Big) \\
     =&\ \norm{f'''} \delta \Big(1 + \frac{1}{2} \lambda^2\E U^2  ( \E (Q-R)^2/R^2 + 1  ) \Big).
\end{align*}
We arrive at the stated  bound on $\abs{\epsilon_{A}(f)}$.

We now prove \eqref{eq:eD:gminf}.   Taylor expansion yields 
\begin{align*}
      \mu Q \big( f(\TX -\delta)-f(\TX)\big)  =&\ -\mu \delta Q f'(\TX) + \frac{1}{2} \mu \delta^2 Q f''(\TX) + \frac{1}{6} \mu \delta^3 Q f'''(\xi).
\end{align*}
Since $\mu R = \lambda$, $X = \delta (Q-R)$, and $\TX - X = -\delta \lambda R_{a}$, the first two terms on the right-hand side satisfy
\begin{align*}
    -\mu \delta Q f'(\TX) =&\  -\mu \delta R f'(\TX)  -\mu \delta (Q-R) f'(\TX)  \\
    =&\   -\mu \delta R f'(\TX) -\mu \delta (Q-R) f'(X) -\mu \delta (Q-R) ( f'(\TX) -f'(X)) \\
    =&\ -\delta \lambda f'(\TX) - \mu X f'(X) + \mu X \delta  \lambda R_{a} f''(\xi),
\end{align*}
and 
\begin{align*}
    \frac{1}{2} \mu \delta^2 Q f''(\TX) =&\ \frac{1}{2} \mu \delta^2 R f''(X) + \frac{1}{2} \mu \delta^2 (Q-R) f''(X)  + \frac{1}{2} \mu \delta^2 Q (f''(\TX) - f''(X))    \\
    =&\ \frac{1}{2} \mu   f''(\TX) + \frac{1}{2} \mu \delta X f''(X)- \frac{1}{2} \mu \delta^3 Q  \lambda R_{a}  f'''(\xi),
\end{align*}
implying that \eqref{eq:eD:verbose:gminf} holds with 
\begin{align}
        \epsilon_{D}(f) =&\  \mu \delta \E \big(X \lambda R_{a} f''(\xi)\big) - \frac{1}{2} \mu \delta^3  \E \big(Q   \lambda R_{a} f'''(\xi)\big)+\frac{1}{6} \mu \delta^3  \E \big(Q f'''(\xi)\big). \label{eq:eD:verbose:gminf}
\end{align}
To conclude, we recall that $\E Q = R$ and  bound each of the three terms on the right-hand side as follows:
\begin{align*}
    \mu \delta \E \abs{X \lambda R_{a} f''(\xi)} \leq&\ \mu \delta \norm{f''} \sqrt{ \E X^2} \sqrt{\lambda^2\E R_{a}^2} = \mu \delta \norm{f''} \sqrt{ \E X^2} \sqrt{\lambda^3\E U^3/3}, \\
    \frac{1}{2} \mu \delta^3  \E \big(Q  \lambda R_{a} f'''(\xi)\big) \leq&\ \frac{1}{2} \mu \norm{f'''} \delta^3  \E \big(Q \lambda R_{a} \big) \\
    \frac{1}{6} \mu \delta^3  \E \big(Q f'''(\xi)\big) \leq&\ \frac{1}{6} \mu \norm{f'''} \delta^3  \E Q = \frac{1}{6} \mu \norm{f'''} \delta.
\end{align*}  
\finishproof

\subsection{Proof of Lemma~\ref{lem:relationships:gminf}}
\label{sec:relationships:proof:gminf}
\startproof{}
The first two  equalities in \eqref{eq:standard:gminf} are argued as in Lemma~\ref{lem:relationships:workload} and the third is obtained from the BAR \eqref{eq:full:BAR:gminf} with $f(z) = q \wedge M$. For the relationship in \eqref{eq:mombound:gminf}, we invoke the BAR  with    $f(z) = ((q \wedge M)  - \lambda r_{a})^2$ to get  
\begin{align*}
    0 =&\ 2\lambda  \E    \big(Q \wedge M - \lambda R_{a}\big) + \mu \E \big(Q 1(Q \leq M)  \big( ( Q-1   - \lambda R_{a})^2 - ( Q  - \lambda R_{a})^2 \big) \big) \\
    &+ \E \int_{0}^{1} 1(Q(t-) \leq M-1) \big( ( Q(t-)+1 - \lambda U(t))^2 - ( Q(t-))^2  \big) d A(t)\\
    =&\ 2\lambda  \E    \big(Q \wedge M - \lambda R_{a}\big) + \mu \E\big( Q 1(Q \leq M)  \big( -2(Q- \lambda R_{a}) + 1 \big) \big) \\
    &+ \E \int_{0}^{1} 1(Q(t-) \leq M-1) \big( 2 Q(t-)(1- \lambda U(t)) + (1-\lambda U(t))^2 \big) d A(t)
\end{align*}
Taking $M \to \infty$ yields 
\begin{align*}
    0 =&\ 2\lambda \E Q - 2\lambda^2 \E R_{a} -2\mu \E Q^2 + 2\mu \lambda  \E Q R_{a} + \mu \E Q + \E(1-\lambda U)^2 \E A(1).
\end{align*}
Divide by $2\mu$, move $\E Q^2$ to the left-hand side, and use $\E Q = R$ to get 
\begin{align*}
    \E Q^2 =&\ \E (Q-R)^2 + R^2 \\
    =&\ R^2 -  R \lambda^2 \E U^2/2 + \lambda \E Q R_{a} + R/2 +  R \E(1-\lambda U)^2/2,
\end{align*}
implying that
\begin{align*}
    \E (Q-R)^2  =&\  -   R \lambda^2 \E U^2/2 + \lambda \E Q R_{a} + R/2 +   R (1 - 2 + \lambda^2 \E U^2)/2 \\
    =&\  \lambda \E Q R_{a}  \\
    =&\ \lambda \E (Q-R) R_{a} + R \lambda^2 \E U^2/2 \\
    \leq&\ \lambda \sqrt{ \E (Q-R)^2} \sqrt{\E R_{a}^2} + R \lambda^2 \E U^2/2.
\end{align*}
Since  $x^2 \leq bx + c$ implies that $x^2 \leq (b^2 + x^2)/2 + c$ and, therefore,  $x^2 \leq b^2  +2 c$, and since $\E R_{a}^{2} = \lambda \E U^3/3$, we obtain \eqref{eq:mombound:gminf}.

To prove \eqref{eq:qt-:gminf} we use  $f(z) = (q \wedge m)^2$ in the BAR to get 
\begin{align*}
    0 =&\ \mu \E \big( Q  1(Q \leq M) ( -2Q + 1)\big)   +  \E  \int_{0}^{1} 1(Q(t-) \leq M-1) (2 Q(t-)+1) d A(t) \\
    =&\ - 2 \mu \E(Q^2 1(Q \leq M)) + \mu E(Q 1(Q \leq M)) \\
    &+ 2\E  \int_{0}^{1} 1(Q(t-) \leq M-1)   Q(t-)  d A(t) + \E  \int_{0}^{1} 1(Q(t-) \leq M-1)  d A(t).
\end{align*}
Taking $M \to \infty$ and using $\E Q = R$ we arrive at 
\begin{align*}
    2 \E Q^2 =&\ \E Q +  \frac{2}{\mu} \E  \int_{0}^{1}    Q(t-)  d A(t) + \frac{1}{\mu} \E A(1) \\
    =&\ R + \frac{2}{\mu} \E  \int_{0}^{1}    Q(t-)  d A(t) + R.
\end{align*}
Using   $\E A(1) = \lambda$ and rearranging terms, we get 
\begin{align*}
    \frac{1}{\mu }\E  \int_{0}^{1}    Q(t-)  d A(t) = \E Q^2 - R = \E (Q-R)^2 + R^2 -R.
\end{align*} 
\finishproof

\subsection{Stein factor bounds for the Normal distribution}
\label{app:stein:normal}
\startproof{Proof of Lemma~\ref{lem:normal:stein:factors}}
Given $h \in \lipone$, define $\bar h(x) =   h(\sigma x)/\sigma$ and $\bar f(x) = f_{\bar h, 1}(x)$. Let us verify that  $f_{h,\sigma}(x) = \sigma \bar f(x/\sigma)$ solves the Poisson equation. Since 
\begin{align*}
    f_{h,\sigma}'(x) = \bar f'(x/\sigma) \quad \text{ and } \quad f_{h,\sigma}''(x) = \bar f''(x/\sigma)/\sigma,
\end{align*}
it follows that 
\begin{align*}
    - x f_{h,\sigma}'(x) + \frac{1}{2} \sigma^2 f_{h,\sigma}''(x)  =&\ \sigma \big( - \frac{x}{\sigma} \bar f'(x/\sigma) + \frac{1}{2}   \bar f''(x/\sigma) \big)  \\
    =&\ \frac{\sigma}{\mu }   ( \E \bar h( Y/\sigma) - \bar h(  x/\sigma))  = \frac{1}{\mu } ( \E h(Y) - h(x)),
\end{align*}
and our Stein factor bounds follow from the standard result \citep[Lemma 2.4]{ChenGoldShao2011} that:
\begin{align*}
        \norm{\bar f'} \leq \frac{2}{\mu}, \quad \norm{\bar f''} \leq \frac{\sqrt{2/\pi}}{\mu }, \quad \text{ and } \quad \norm{\bar f'''} \leq \frac{2}{\mu  }.
    \end{align*}
\finishproof

\section{The JSQ system: supporting proofs}
\label{app:jsq}
Let $A(t)$ and $U(t)$ be defined as in Section~\ref{sec:model:workload}. Similarly, let  $D_{i}(t)$ denote the number of departures from server $i$  on $[0,t]$ and if a departure occurs from server $i$ at time $t$, we let $S_{i}(t)$ be the service time of the next customer to be served by that server. We also let $\Lambda_{i}(t)$ denote the (random) remaining time until a customer gets routed to server $i$. The following BAR is proved just like Lemma~\ref{lem:full:BAR:workload}.
\begin{lemma}
\label{lem:full:BAR:JSQ}
Initialize $Z(0) \sim Z$. If $f(Z(s))$ satisfies the FTC conditions with probability one under $Z(0) \sim Z$ and if the integrability condition 
\begin{align*}
    \E \abs{f(Z)}, \  \E \abs{\partial_{r_{a}} f(Z)}, \  \E \abs{\partial_{r_{s,i}} f(Z)},\  \E \int_{0}^{t} \abs{\Delta  f(Z(s-))} dD_{i}(s),\  \E \int_{0}^{t} \abs{\Delta  f(Z(s-))} dA(s) < \infty
\end{align*} 
holds, then  
\begin{align}
    &- \E (\partial_{r_a} f(Z)) -   \sum_{i=1}^{n}\E (1(Q_i>0) \partial_{r_{s,i}} f(Z))   \notag \\
    &+  \E  \int_{0}^{1} \Delta f(Z(t-)) d A(t) + \sum_{i=1}^{n} \E  \int_{0}^{1} \Delta f(Z(t-)) d D_{i}(t) = 0.\label{eq:full:BAR:JSQ}
\end{align}
\end{lemma}
Define the compensated total customer count 
\begin{align*}
    \TX(t) = X(t) - \delta \lambda R_{a}(t) + \sum_{i=1}^{n} \delta \mu R_{s,i}(t) \quad \text{ and } \quad \TX = X  - \delta \lambda R_{a}  + \sum_{i=1}^{n} \delta \mu R_{s,i}.
\end{align*}
Specialized to $\TX$, the BAR \eqref{eq:full:BAR:JSQ} becomes 
\begin{align}
    &\delta \E \Big(\Big(\lambda - \sum_{i=1}^{n} 1(Q_i>0) \mu\Big) f'(\TX)\Big) +  \E  \int_{0}^{1} \Delta f(\TX(t-)) d A(t) \notag \\
    &+ \sum_{i=1}^{n} \E  \int_{0}^{1} \Delta f(\TX(t-)) d D_{i}(t) = 0. \label{eq:BAR3}
\end{align}
The following is an analog of  Proposition~\ref{prop:extraction:workload}, and  is proved in Appendix~\ref{app:extraction:JSQ}.
\begin{proposition}
\label{prop:extraction:JSQ}
If  $f\in C^{2}(\R)$ with $f''(x)$ absolutely continuous and $\norm{f''},\norm{f'''} < \infty$, then, provided that all expectations are well defined,
\begin{align}
    \delta \E \Big(\Big(\lambda - \sum_{i=1}^{n} 1(Q_i>0) \mu\Big) f'(\TX)\Big) =&\ -n\mu \delta^2 \E f'(X) + n \mu \delta^2 f'(0) + \epsilon_{0}(f), \label{eq:e0:JSQ} \\
        \E \int_{0}^{1} \Delta  f(\TX(t-))  d A(t) =&\ \frac{1}{2}\delta^2  \lambda c_{U}^{2} \E f''(X) + \epsilon_{A}(f),  \label{eq:eA:JSQ} \\
         \E \int_{0}^{1} \Delta  f(\TX(t-))  d D_{i}(t) =&\  \frac{1}{2} \delta^2  \mu c_{S}^{2} \E f''(X) + \epsilon_{D,i}(f), \label{eq:eD:JSQ}
    \end{align}
    where 
\begin{align*}
    \abs{\epsilon_{0}(f)} \leq&\ \frac{1}{2} n \mu \delta^3 \norm{f''} \big(\lambda^2 \E U^2/2 +   ( \lambda \mu  \E S^2/2 + \mu \delta\E S)\big) \\
    &+  \delta^2 \mu \norm{f''} \sum_{i=1}^{n} \E \Big(1(Q_i=0)   \Big(\sum_{j=1}^{n} Q_j + \lambda R_{a} + \sum_{j=1}^{n} \mu R_{s,j}\Big)   \Big), \\
    \abs{\epsilon_{A}(f)} \leq&\      \frac{1}{2} \delta^3 \norm{f'''} \Big( \frac{1}{3} \lambda \E\abs{1-\lambda U}^{3} +   c_{U}^{2}  n\mu \lambda(\rho (\lambda/n) \E S^2/2 + \delta \E S+ \lambda \E U^2/2) \\
    & \hspace{2cm} + \lambda c_{U}^{2}   (1 + n + \lambda \mu \E U^2 +  \mu^2 \E S^2 )\Big),\\
    \abs{\epsilon_{D,i}(f)} \leq&\ \frac{1}{2} \delta^3 \norm{f'''} \bigg( \frac{1}{3}   \frac{\lambda}{n} \E\abs{1 - \mu S}^{3} \\
    &+    c_{S}^{2}   \big(   (\lambda^2 \mu \rho + 2(n-1)\mu^2 \rho \lambda)   \E S^2/2n    + \delta(\lambda \mu + 2(n-1)\mu^2)\E S  + \frac{1}{2} \lambda^2 \mu \E U^2 \big)\\ 
    &+   \mu c_{S}^{2}     n \rho  + \mu c_{S}^{2}      (\lambda /n) \big( \E S+\lambda \E S^2 + \lambda^2 \E U^2 \E S \big)  \\
    &+  c_{S}^{2}  \lambda \big( 1 +    (\lambda \mu \E U^2/2 + 1/\rho)  + 2\mu^2 \E S^2  \big)\bigg) +\frac{1}{2} \delta^3  \norm{f''}  c_{S}^{2}\mu\lambda/n.
\end{align*}
 \end{proposition} 
The expansions in \eqref{eq:e0:JSQ}--\eqref{eq:eD:JSQ} suggest the diffusion   generator 
\begin{align*}
  G_{Y} f(x) =  -n\mu \delta^2 f'(x) + \frac{1}{2}\delta^2 ( \lambda c_{U}^{2} + n\mu c_{S}^{2}) f''(x) + n\mu \delta^2 f'(0), \quad x \geq 0,
\end{align*}
which, like in Section~\ref{sec:workload}, corresponds to the  exponential distribution with mean $( \rho c_{U}^{2} +   c_{S}^{2})/2$.  We now combine everything to prove Theorem~\ref{thm:jsq:main}.
\startproof{Proof of Theorem~\ref{thm:jsq:main}} 
Fix $h \in \lipone$. Since  $\E X^2 < \infty$, we can verify that $f_h(\TX)$ satisfies the conditions of the BAR \eqref{eq:full:BAR:JSQ} in the same way as we did in the proof of Theorem~\ref{thm:workload:1}. For the Stein factor bounds, we use  Lemma~\ref{lem:exp:stein:factors} with $\theta = n \mu \delta^2$ and $\sigma^2 = \delta^2 (\lambda c_{U}^{2} + n \mu c_{S}^{2}) / 2$ there to get
\begin{align*}
     \norm{f_h''} \leq \frac{1}{n \mu \delta^2} \quad \text{ and } \quad \norm{f_h'''} \leq \frac{4}{\delta^2(\lambda c_{U}^{2} + n \mu c_{S}^{2})}.
\end{align*}
Combining these with the bounds on   $\abs{\epsilon_{0}(f_{h})}$, $\abs{\epsilon_{A}(f_{h})}$, and $\abs{\epsilon_{D,i}(f_{h})}$ in Proposition~\ref{prop:extraction:JSQ} concludes the proof.
\finishproof

\subsection{Expanding the BAR}
\label{app:extraction:JSQ}
To prove  Proposition~\ref{prop:extraction:JSQ} we require  two auxiliary lemmas. The first lemma contains some useful relationships derived using the BAR.
\begin{lemma}
    \label{lem:relationships:JSQ}
    For any $1 \leq i \neq j  \leq n$ and $m > 1$,  
    \begin{align}
    & \E A(1) = n\E D_{i}(1) = \lambda,  \quad  \Prob(Q_i > 0) = \rho, \label{eq:arrival:departure:JSQ} \\ 
        &\E R_{a}^{m-1} =  \lambda \E U^{m}/m, \quad \E (R_{s,i}^{m-1} 1(Q_i>0)) =   (\lambda/n)\E S^{m}/m, \quad \E (R_{s,i}^{m}|Q_{i}=0)= \E S^{m}, \label{eq:residual:moments:JSQ} \\
        & \E \int_{0}^{1} R_{a}(t)   d D_{i}(t)   \leq \mu(\E R_{s,i} + \E R_{a}), \quad \E \int_{0}^{1} R_{s,i}(t)  d A(t) \leq \lambda(\E R_{s,i} + \E R_{a}), \label{eq:mix1:JSQ}  \\
        & \E \int_{0}^{1} R_{s,i}(t)  d D_{j}(t) \leq \mu(\E R_{s,i} + \E R_{s,j}). \label{eq:mix2:JSQ}
    \end{align}
\end{lemma}
\startproof{Proof of Lemma~\ref{lem:relationships:JSQ}}
We prove each of \eqref{eq:arrival:departure:JSQ}--\eqref{eq:mix2:JSQ} by plugging in truncated test functions into the BAR \eqref{eq:BAR3}. Using $f(z)  = r_{a} \wedge M$, $f(z) = (q_1+\cdots+q_n) \wedge M$ (together with the symmetry of the servers), and $f(z) = r_{s,i} \wedge M$   yields $\E A(1)=\lambda$,   $n\E D_{i}(1) = \E A(1)$, and $\E S \E D_{i}(1) = \Prob(Q_i>0)$, respectively, which proves \eqref{eq:arrival:departure:JSQ}. Using $f(z) = r_{a}^m \wedge M$ yields the first equality in \eqref{eq:residual:moments:JSQ} and using $f(z) = r_{s,i}^{m} \wedge M$ yields $m\E (R_{s,i}^{m-1} 1(Q_i>0)) = \E S^{m} \E D(1) = (\lambda/n) \E S^m$. The  third equality in \eqref{eq:residual:moments:JSQ} follows  once we note that $\E (R_{s,i}^{m} 1(Q_i=0)) = \E S^{m} \Prob(Q_i=0)$ because when $Q_i(t) = 0$, then $R_{s,i}(t)$ equals the service time of the next customer to arrive.  Using $f(z) = (r_a \wedge M) (r_{s,i} \wedge M)$  yields
\begin{align*}
    &\E \int_{0}^{1} R_{a}(t) \big(1(Q_{i}(t)=0) \Lambda_{i}(t) + S_{i}(t)\big) d D_{i}(t) + \E \int_{0}^{1} R_{s,i}(t) U(t) d A(t)  \\
    =&\ \E R_{s,i} + \E 1(Q_i>0)R_{a},
\end{align*}
from which \eqref{eq:mix1:JSQ}  follows. The inequality  \eqref{eq:mix2:JSQ} is proved similarly using $f(z) = (r_{s,i} \wedge M) (r_{s,j} \wedge M)$.
\finishproof
The second lemma is the Palm inversion formula, analogous to Lemma~\ref{lem:palm:inversion:2}. For simplicity, we only prove it for bounded functions $f(x)$, though this assumption can be relaxed as in Lemma~\ref{lem:palm:inversion:2}.
\begin{lemma}
    \label{lem:palm:inversion:JSQ}
    For any bounded $f: \R \to \R$, 
    \begin{align}
        \E f(X) =&\ \E \Big( \int_{0}^{1} \int_{0}^{U(t)} f(X(t+u))du  d A(t)\Big)  \quad \text{ and }  \label{eq:palm:inversion:arrive:JSQ}\\
        \E f(X)  =&\ \E \Big( \int_{0}^{1}   \int_{0}^{1(Q_i(t)=0)\Lambda_{i}(t) + S_{i}(t)} f(X(t+u))du   d D_{i}(t)\Big). \label{eq:palm:inversion:depart:JSQ}
    \end{align}
\end{lemma}
\startproof{Proof of Lemma~\ref{lem:palm:inversion:JSQ}}
Since $\E R_{a} < \infty$ (Lemma~\ref{lem:relationships:JSQ}), the proof of \eqref{eq:palm:inversion:arrive:JSQ} is identical to that of \eqref{eq:palm:inversion:arrive:workload} in Lemma~\ref{lem:palm:inversion:2}. Similarly, to prove \eqref{eq:palm:inversion:depart:JSQ}, let $R_{d,i}(t) = 1(Q_{i}(t)=0)\Lambda_{i}(t)   + R_{s,i}(t)$ be the remaining time until the next departure from server $i$  and let $R_{d,i}$ denote the steady-state version. We only need to argue that $\E R_{d,i} < \infty$. This is implied by the facts that 1)  the expected time until the next customer arrives to server $i$ is finite because $\E R_{a} < \infty$ and because of our uniform tie-breaking rule (meaning that an arrival will be allocated an idling server with probability at least $1/n$) and 2) $\E R_{s,i} < \infty$ (Lemma~\ref{lem:relationships:JSQ}).
\finishproof

\startproof{Proof of Proposition~\ref{prop:extraction:JSQ}}
The proof has two parts. In part one we assume that \eqref{eq:e0:JSQ}--\eqref{eq:eD:JSQ} hold  with
\begin{align}
        \epsilon_{0}(f) =&\ - n \mu \delta^3 \E \Big(\Big(-\lambda R_{a} +   \sum_{i=1}^{n} \mu R_{s,i}\Big)f''(\xi)\Big) \notag \\
        &+ \delta^2  \mu \sum_{i=1}^{n} \E \Big(1(Q_i=0)   \Big(\sum_{j=1}^{n} Q_j - \lambda R_{a} + \sum_{j=1}^{n} \mu R_{s,j}\Big) f''(\xi) \Big),  \label{eq:e0:verbose:JSQ}\\
        \epsilon_{A}(f) =&\ \frac{1}{6} \delta^3 \E \int_{0}^{1}(1 - \lambda U(t))^3 f'''(\xi(t))d A(t) \notag \\
    &+ \frac{1}{2} \delta^3 c_{U}^{2} \E \int_{0}^{1} \Big(\sum_{i=1}^{n}  \mu R_{s,i}(t)\Big) f'''(\xi(t-))d A(t) \notag \\
    &- \frac{1}{2} \delta^2 \lambda c_{U}^{2} \E \int_{0}^{1} \int_{0}^{U(t)} \big( X(t+u)-X(t-)\big) f'''(\xi(t+u))du d A(t), \label{eq:eA:verbose:JSQ}
     \end{align}
    and
    \begin{align}
        \epsilon_{D,i}(f) =&\ \frac{1}{6} \delta^3 \E \int_{0}^{1}(1 - \mu S_{i}(t))^3 f'''(\xi(t))d D_{i}(t) \notag \\
    &+ \frac{1}{2} \delta^3 c_{S}^{2} \E \int_{0}^{1} \Big(-  \lambda R_{a}(t) + \sum_{j \neq i}  \mu R_{s,j}(t)\Big) f'''(\xi(t-))d D_{i}(t) \notag \\
    &- \frac{1}{2} \delta^2  \mu c_{S}^{2} \E \int_{0}^{1}  \int_{0}^{1(Q_{i}(t)=0 ) \Lambda_{i}(t)} f''(X(t+u))  du d D_{i}(t)  \notag  \\
    &- \frac{1}{2} \delta^2  \mu c_{S}^{2} \E \int_{0}^{1} \int_{1(Q_{i}(t)=0 ) \Lambda_{i}(t)}^{1(Q_{i}(t)=0 ) \Lambda_{i}(t) +S_{i}(t)} \big(X(t+u)-X(t-)\big)f'''(\xi(t+u)) du  d D_{i}(t), \label{eq:eD:verbose:JSQ}
     \end{align}
     and we bound the terms on   the right-hand side, with the exception of  the boundary term in $\epsilon_{0}(f)$. In part two, we expand all  terms in the BAR \eqref{eq:BAR3} to verify that \eqref{eq:e0:JSQ}--\eqref{eq:eD:JSQ} does indeed hold with   $\epsilon_{0}(f), \epsilon_{A}(f)$, and $\epsilon_{D,i}(f)$ as in \eqref{eq:e0:verbose:JSQ}--\eqref{eq:eD:verbose:JSQ}.

\textbf{Part one.}  To prove  the bound on $\epsilon_{0}(f)$, we bound the non-boundary term  by noting that
\begin{align*}
    n \mu \delta^3 \E \Big(\Big(\lambda R_{a} +   \sum_{i=1}^{n} \mu R_{s,i}\Big)\norm{f''}\Big) \leq n \mu \delta^3 \norm{f''} \big(\lambda^2 \E U^2/2 +   ( \lambda \mu  \E S^2/2 + \mu \delta\E S)\big),
\end{align*}
where the inequality follows from Lemma~\ref{lem:relationships:JSQ}. Next we bound the three terms composing  $\epsilon_{A}(f)$. For the first term, noting that  $\E A(1)=\lambda$, by Lemma~\ref{lem:relationships:JSQ},  yields 
\begin{align*}
    \frac{1}{6} \delta^3 \E \int_{0}^{1}(1 - \lambda U(t))^3 f'''(\xi(t))d A(t) \leq  \frac{1}{6} \delta^3 \norm{f'''}  \lambda \E\abs{1-\lambda U}^{3}. 
\end{align*}
 For the second term, we use Lemma~\ref{lem:relationships:JSQ} to get
\begin{align*}
     &\frac{1}{2} \delta^3 c_{U}^{2} \E \int_{0}^{1} \Big(\sum_{i=1}^{n}  \mu R_{s,i}(t)\Big) f'''(\xi(t-))d A(t) \\
     \leq&\ \frac{1}{2} \delta^3 c_{U}^{2} \norm{f'''} n\mu \lambda(\E R_{s,i}+ \E R_{a}) \\
     =&\ \frac{1}{2} \delta^3 c_{U}^{2} \norm{f'''} n\mu \lambda(\rho (\lambda/n) \E S^2/2 + \delta \E S+ \lambda \E U^2/2).
\end{align*}
For the third term, note that
\begin{align*}
     &\frac{1}{2} \delta^2 \lambda c_{U}^{2} \E \int_{0}^{1} \int_{0}^{U(t)} \big( X(t+u)-X(t-)\big) f'''(\xi(t+u))du d A(t) \\
     \leq&\ \frac{1}{2}\delta^2 \lambda c_{U}^{2}  \norm{f'''} \E \int_{0}^{1} U(t) \sup_{0 \leq u \leq U(t)}  \abs{X(t+u)-X(t-)}   d A(t) \\
     \leq&\ \frac{1}{2}\delta^3 \lambda c_{U}^{2}  \norm{f'''} (1 + n + \lambda \mu \E U^2 +  \mu^2 \E S^2 ).
\end{align*} 
To justify the second inequality, let $D_{i}[t,t+u]$ be the   number of service completions by server $i$ on $[t,t+u]$ and note that 
\begin{align*}
    \sup_{0 \leq u \leq U(t)}\abs{X(t+u)-X(t-)} \leq \delta + \delta \sum_{i=1}^{n} D_{i}[t,t+U(t)].
\end{align*}
By forcing server $i$ to continue working even if it has no customers in its buffer, we can  construct a zero-delayed counting process $\{\bar D_{i}(t), t \geq 0\}$ keeping track of \emph{potential} service completions, which would be independent of the arrival process, and satisfy
\begin{align*}
    D_{i}[t,t+u] \leq 1 + \bar D_{i}(u).
\end{align*}
Since $\E \bar D_{i}(u) \leq \mu u + \mu^2 \E S^2$ by Lorden's inequality (\cite{Lord1970}), it follows that 
\begin{align}
    \E \int_{0}^{1} U(t) \sup_{0 \leq u \leq U(t)}  \abs{X(t+u)-X(t-)}   d A(t)   \leq&   \delta  \E \Big(U \big(1 + n + \sum_{i=1}^{n} \bar D_{i}(U)\big)\Big) \E A(1) \notag \\
    \leq&\  \delta \E \Big(U \big(1 + n +\mu U + \mu^2 \E S^2 \big)\Big) \E A(1) \notag \\
    =&\  \delta (1 + n + \lambda \mu \E U^2 +  \mu^2 \E S^2  ) \label{eq:lorden:use}
\end{align}
where in the last equality we used $ \E U \E A(1) = 1$.  Lastly, we bound the four terms composing $\epsilon_{D,i}(f)$. First, since $\E D_{i}(1)  = \lambda/n$ by Lemma~\ref{lem:relationships:JSQ}, we have
\begin{align*}
       \frac{1}{6} \delta^3 \E \int_{0}^{1}(1 - \mu S_{i}(t))^3 f'''(\xi(t))d D_{i}(t) \leq \frac{1}{6} \delta^3  \norm{f'''} \frac{\lambda}{n} \E\abs{1 - \mu S}^{3}.
     \end{align*}
Second, using Lemma~\ref{lem:relationships:JSQ}, we have 
\begin{align*}
    & \frac{1}{2} \delta^3 c_{S}^{2} \E \int_{0}^{1} \Big( \lambda R_{a}(t) + \sum_{j \neq i}  \mu R_{s,j}(t)\Big) f'''(\xi(t-))d D_{i}(t) \\
    \leq&\    \frac{1}{2} \delta^3 c_{S}^{2} \norm{f'''}  \Big(  \lambda \mu(\E R_{s,i} + \E R_{a}) + \sum_{j\neq i}  \mu^2(\E R_{s,i} + \E R_{s,j})\Big)  \\
    \leq&\ \frac{1}{2} \delta^3 c_{S}^{2} \norm{f'''}  \Big(  \lambda \mu(\rho (\lambda/n)\E S^2/2 + \delta \E S + \lambda \E U^2/2 ) +  2(n-1) \mu^2(\rho (\lambda/n)\E S^2/2 + \delta \E S)\Big) \\
    =&\ \frac{1}{2} \delta^3 c_{S}^{2} \norm{f'''}  \Big(  \frac{\lambda^2 \mu \rho + 2(n-1)\mu^2 \rho \lambda} {2n}   \E S^2    + \delta(\lambda \mu + 2(n-1)\mu^2)\E S  + \frac{1}{2} \lambda^2 \mu \E U^2 \Big).
 \end{align*}
Third, the Palm inversion formula (Lemma~\ref{lem:palm:inversion:JSQ}) applied to $f(z)=1(q_i=0)$, together with Lemma~\ref{lem:relationships:JSQ} yield  
\begin{align*}
    \E \int_{0}^{1} 1(Q_{i}(t)=0 ) \Lambda_{i}(t) d D_{i}(t) = \Prob(Q_i=0) = \delta, 
\end{align*}
so that 
\begin{align*}
           \frac{1}{2} \delta^2  \mu c_{S}^{2} \E \int_{0}^{1}  \int_{0}^{1(Q_{i}(t)=0 ) \Lambda_{i}(t)} f''(X(t+u))  du d D_{i}(t) \leq  \frac{1}{2} \delta^3  \norm{f''}  c_{S}^{2}\mu\lambda/n.
\end{align*}
Fourth,  
     \begin{align*}
       & \frac{1}{2} \delta^2  \mu c_{S}^{2} \norm{f'''} \E \int_{0}^{1} \int_{1(Q_{i}(t)=0 ) \Lambda_{i}(t)}^{1(Q_{i}(t)=0 ) \Lambda_{i}(t) +S_{i}(t)} \abs{X(t+u)-X(t-)}  du  d D_{i}(t) \\ 
       \leq&\ \frac{1}{2} \delta^2  \mu c_{S}^{2} \norm{f'''} \E \int_{0}^{1}  S_{i}(t) \sup_{0 \leq u \leq S_{i}(t)} \abs{X(t+1(Q_{i}(t)=0 ) \Lambda_{i}(t)+u)-X(t-)}   d D_{i}(t)\\ 
     \end{align*}
     Observe that  
     \begin{align*}
   \sup_{0 \leq u \leq S_{i}(t)} \abs{X(t+u)-X(t-)} \leq&\   \delta A\big[t,t+1(Q_{i}(t)=0 ) \Lambda_{i}(t)\big]  \\
    &+ \delta A\big(t+1(Q_{i}(t)=0 ) \Lambda_{i}(t),t+1(Q_{i}(t)=0 ) \Lambda_{i}(t) + S_i(t)\big] \\
    &+ \delta \sum_{j \neq i} D_{j}\big[t,t+1(Q_{i}(t)=0 ) \Lambda_{i}(t) + S_i(t)\big],
\end{align*}
where $A[a,b]$ and $D_{j}[a,b]$ are the number of arrivals and server $j$ departures, respectively, on $[a,b]$. To complete the bound on the fourth term in $\epsilon_{D,i}(f)$, we now bound  
\begin{align}
    & \frac{1}{2} \delta^2  \mu c_{S}^{2} \norm{f'''} \E \int_{0}^{1}  S_{i}(t) \delta A\big[t,t+1(Q_{i}(t)=0 ) \Lambda_{i}(t)\big]   d D_{i}(t) \label{eq:a} \\ 
       &+ \frac{1}{2} \delta^2  \mu c_{S}^{2} \norm{f'''} \E \int_{0}^{1}  S_{i}(t) \delta A\big(t+1(Q_{i}(t)=0 ) \Lambda_{i}(t),t+1(Q_{i}(t)=0 ) \Lambda_{i}(t) + S_i(t)\big]   d D_{i}(t)\label{eq:b}\\ 
       &+ \frac{1}{2} \delta^2  \mu c_{S}^{2} \norm{f'''} \sum_{j \neq i} \E \int_{0}^{1}  S_{i}(t) \delta  D_{j}\big[t,t+1(Q_{i}(t)=0 ) \Lambda_{i}(t) + S_i(t)\big]  d D_{i}(t).\label{eq:c}
\end{align}
First, we bound   $A\big[t,t+1(Q_{i}(t)=0 ) \Lambda_{i}(t)\big]$, the number of arrivals before server $i$ gets a customer. Recall that our tie-breaking routing rule is to select among the least loaded servers uniformly at random. If server $i$ becomes idle at time $t$, then an arrival will be routed to this server with probability at least $1/n$, and it is therefore possible to construct a geometrically distributed random variable $\Gamma \geq 1$ with mean $n$ (success probability $1/n$) that a) upper bounds the number arrivals needed to route a customer to server $i$ and b) is independent of $S_{i}(t)$ and $\{D_{j}(t), t \geq 0\}_{j=1}^{n}$. It follows that  $A\big[t,t+1(Q_{i}(t)=0 ) \Lambda_{i}(t)\big] \leq \Gamma$. Furthermore, we define the time until the $\Gamma$th arrival as 
\begin{align*}
    \overline\Lambda_{i}(t) = R_{a}(t) + \sum_{m=2}^{\Gamma} U(\tau_{m}),
\end{align*}
where $\tau_{m}$ is the time of the $m$th arrival after time $t$, and note that 
\begin{align*}
    1(Q_{i}(t)=0)\Lambda_{i}(t) \leq  \overline\Lambda_{i}(t).
\end{align*}
We now have the tools needed to bound \eqref{eq:a}--\eqref{eq:c}. First, we observe that 
\begin{align*}
    &\frac{1}{2} \delta^2  \mu c_{S}^{2} \norm{f'''} \E \int_{0}^{1}  S_{i}(t) \delta A\big[t,t+1(Q_{i}(t)=0 ) \Lambda_{i}(t)\big]   d D_{i}(t) \\ 
    \leq&\ \frac{1}{2} \delta^3  \mu c_{S}^{2} \norm{f'''} \E \int_{0}^{1}  S_{i}(t) \Gamma   d D_{i}(t)  =  \frac{1}{2} \delta^3  \mu c_{S}^{2} \norm{f'''} \E S \E \Gamma \E D_{i}(1) = \frac{1}{2} \delta^3   \mu c_{S}^{2} \norm{f'''}   n \rho,
\end{align*}
where in the last equality we used $\E \Gamma = n$, $\E S = 1/\mu$, and $\E D_{i}(1) = \lambda/n$. Second, we argue using Lorden's inequality, as we did in \eqref{eq:lorden:use}, that 
\begin{align*}
    &\frac{1}{2} \delta^2  \mu c_{S}^{2} \norm{f'''} \E \int_{0}^{1}  S_{i}(t) \delta A\big(t+1(Q_{i}(t)=0 ) \Lambda_{i}(t),t+1(Q_{i}(t)=0 ) \Lambda_{i}(t) + S_i(t)\big]   d D_{i}(t) \\
    \leq&\ \frac{1}{2} \delta^3  \mu c_{S}^{2} \norm{f'''} \E \int_{0}^{1}  S  (1+\lambda S + \lambda^2 \E U^2)  d D_{i}(t)\\
    =&\ \frac{1}{2} \delta^3  \mu c_{S}^{2} \norm{f'''}     \E \big(S(1+\lambda S + \lambda^2 \E U^2)\big)   \E D_{i}(1) \\
    =&\ \frac{1}{2} \delta^3 \mu c_{S}^{2} \norm{f'''}     (\lambda /n) \big( \E S+\lambda \E S^2 + \lambda^2 \E U^2 \E S \big).
\end{align*}
Third, we bound \eqref{eq:c} 
\begin{align*}
    & \frac{1}{2} \delta^2  \mu c_{S}^{2} \norm{f'''} \sum_{j \neq i} \E \int_{0}^{1}  S_{i}(t) \delta  D_{j}\big[t,t+1(Q_{i}(t)=0 ) \Lambda_{i}(t) + S_i(t)\big]  d D_{i}(t) \\
    \leq&\ \frac{1}{2} \delta^3  \mu c_{S}^{2} \norm{f'''} \sum_{j \neq i} \E \int_{0}^{1}  S    D_{j}\big[t,t+\overline \Lambda_{i}+ S\big]  d D_{i}(t)\\
    \leq&\ \frac{1}{2} \delta^3  \mu c_{S}^{2} \norm{f'''} (n-1) \E \big(S   \big(1 + \mu(\overline \Lambda_{i}+ S) + \mu^2 \E S^2\big) \big)     \E D_{i}(1)\\
    =&\ \frac{1}{2} \delta^3  \mu c_{S}^{2} \norm{f'''} (n-1)  \big( \E S +   \E \overline \Lambda_{i} + 2\mu \E S^2  \big)   \lambda/n \\
    \leq&\ \frac{1}{2} \delta^3  \mu c_{S}^{2} \norm{f'''}  \lambda \big( \E S +   (\E R_{a} + \E U \E \Gamma)   +2 \mu  \E S^2\big)   \\
    =&\ \frac{1}{2} \delta^3  c_{S}^{2} \norm{f'''}  \lambda \big( 1 +    (\lambda \mu \E U^2/2 + 1/\rho)  + 2\mu^2 \E S^2  \big).
\end{align*}

\textbf{Part two.} We first prove \eqref{eq:e0:JSQ}:
\begin{align*}
    \delta \E \Big((\lambda - \sum_{i=1}^{n} 1(Q_i>0) \mu) f'(\TX)\Big) =&\ -n \mu \delta^2 \E f'(\TX) + \delta  \mu \sum_{i=1}^{n} \E \big(1(Q_i=0) f'(\TX)\big).
\end{align*}
Recalling that   $\TX = \delta \sum_{i=1}^{n} Q_i - \delta \lambda R_{a} + \delta \sum_{i=1}^{n} \mu R_{s,i}$,  we have
\begin{align*}
    -n \mu \delta^2 \E f'(\TX) = -n \mu \delta^2 \E f'(X) - n \mu \delta^3 \E \Big(\Big(-\lambda R_{a} +   \sum_{i=1}^{n} \mu R_{s,i}\Big)f''(\xi)\Big),
\end{align*}
and 
\begin{align*}
&\delta  \mu \sum_{i=1}^{n} \E \big(1(Q_i=0) f'(\TX)\big) \\
=&\  \delta  \mu \sum_{i=1}^{n} \E \Big(1(Q_i=0) \big( f'(0) + \delta \Big(\sum_{i=1}^{n} Q_i - \lambda R_{a} + \sum_{i=1}^{n} \mu R_{s,i}\Big) f''(\xi)\big)\Big),
\end{align*}
which proves \eqref{eq:e0:verbose:JSQ}. The proof of \eqref{eq:eA:verbose:JSQ} is very similar to the proof of \eqref{eq:eA:workload} in Proposition~\ref{prop:extraction:workload}:  if an arrival happens at time $t$, then  $\Delta \TX(t-) = \delta(1-\lambda U(t))$ and, since $\E (1-\lambda U) = 0$, 
\begin{align}
    \E \int_{0}^{1} \Delta  f(\TX(t-))  d A(t)  =&\   \frac{1}{2} \delta^2\E (1 -\lambda  U)^2  \E \int_{0}^{1}f''(\TX(t-))d A(t) \notag \\
    &+ \frac{1}{6} \delta^3\E \int_{0}^{1} (1 - \lambda U(t))^3 f'''(\xi(t))d A(t), \label{eq:analogous:JSQ}
\end{align}
and the first term expands into 
\begin{align*}
    &\frac{1}{2} \delta^2\E (1 -\lambda  U)^2 \E \int_{0}^{1}f''(X(t-))d A(t) \\
    &+ \frac{1}{2} \delta^3\E (1 -\lambda  U)^2 \E \int_{0}^{1} \big(\TX(t-)-X(t-)\big) f'''(\xi(t-))d A(t).
\end{align*}
Noting that $\E (1 -\lambda  U)^2 = c_{U}^{2}$ and $\TX(t-)-X(t-) = \sum_{i=1}^{n}  \mu R_{s,i}(t)$, to conclude \eqref{eq:eA:verbose:JSQ} we use Lemma~\ref{lem:palm:inversion:JSQ} and repeat the argument following \eqref{eq:analogous:2} to get
\begin{align*}
    &\E \int_{0}^{1}  f''(X(t-))  d A(t)\\
    =&\ \lambda \E f''(X) - \lambda \E  \int_{0}^{1} \int_{0}^{U(t)} \big(X(t+u)-X(t-)\big)f'''(\xi(t+u)) du  d A(t).
\end{align*}
Lastly we prove \eqref{eq:eD:verbose:JSQ}. Similar to \eqref{eq:analogous:JSQ} and the display that follows, we  expand the departure jump term  to get
\begin{align*}
    \E \int_{0}^{1} \Delta  f(\TX(t-))  d D_{i}(t)  =&\   \frac{1}{2} \delta^2\E (1 -\mu  S)^2 \E \int_{0}^{1}f''(X(t-))d D_{i}(t) \\
    &+ \frac{1}{2} \delta^3\E (1 -\mu  S)^2 \E \int_{0}^{1} \big(\TX(t-)-X(t-)\big) f'''(\xi(t-))d D_{i}(t)\\
    &+ \frac{1}{6} \delta^3\E \int_{0}^{1} (1 - \mu S_{i}(t))^3 f'''(\xi(t))d D_{i}(t),
\end{align*}
where $\TX(t-)-X(t-)= -\delta \lambda R_{a}(t) + \sum_{j \neq i}\delta \mu R_{s,j}(t)$. We conclude the proof  by using Lemma~\ref{lem:palm:inversion:JSQ} to get
\begin{align*}
    \E f''(X) =&\ \E \int_{0}^{1}  \int_{0}^{1(Q_{i}(t)=0 ) \Lambda_{i}(t)} f''(X(t+u))  du d D_{i}(t)  + \E S \int_{0}^{1}   f''(X(t-))  d D_{i}(t) \\ 
    &+ \E \int_{0}^{1} \int_{1(Q_{i}(t)=0 ) \Lambda_{i}(t)}^{1(Q_{i}(t)=0 ) \Lambda_{i}(t) +S_{i}(t)} \big(X(t+u)-X(t-)\big)f'''(\xi(t+u)) du  d D_{i}(t).
\end{align*}

\finishproof

\end{APPENDICES}

\end{document}